%% file: smid.tex
\documentclass [12pt]{amsart}

\usepackage{amssymb,amsxtra,amsfonts}
\usepackage{graphics}
\usepackage{amsmath, accents}
\usepackage[colorlinks]{hyperref}

\input{macros}

\input{macros-thm1.1.smid}

\renewcommand{\pii}{\pi_{\hbox{\footnotesize\rm irr}}}

\begin{document}
\date{}

\nocite{HL}		%
\nocite{BDV}		%
\nocite{FG}		%

\title{Symplectic Manifolds and Isomonodromic Deformations}

\author{Philip Boalch}

\thispagestyle{empty}

\maketitle

\thispagestyle{empty}

\ 

\ 

\ 

\ 

\ 

\

\centerline{\input{extimdstex.big.pstex_t}}	
$$\text{ Isomonodromic Deformations (Figure \ref{imd fig}, p.\pageref{imd fig})}$$

\newpage

\ 

\

\begin{abstract}We study moduli spaces of meromorphic connections 
(with arbitrary order poles) over Riemann surfaces together with the 
corresponding spaces of monodromy data (involving Stokes matrices). 
Natural symplectic structures are found and described both explicitly
and from an infinite dimensional viewpoint (generalising the Atiyah-Bott
approach).
This enables us to give an intrinsic symplectic description of 
the isomonodromic deformation equations of Jimbo, Miwa and Ueno,
thereby putting the
existing results for the six Painlev\'e equations
and Schlesinger's equations into a uniform framework. 
\end{abstract}

\ 

\

\vspace{0.1cm}
{\footnotesize
\bf CONTENTS}
\begin{centering}
\footnotesize
$$
\begin{array}{lcr}
\ref{sn: intro}. \text{ Introduction} & \ \qquad\qquad\qquad\ &
\pageref{sn: intro} \\
\ref{sn: triv}. \text{ Meromorphic Connections on Trivial Bundles} &&
\pageref{sn: triv} \\
\ref{sn: gm}. \text{ Generalised Monodromy} &&
\pageref{sn: gm} \\
\ref{sn: smooth}. \text{ ${C^\infty}$  Approach to Meromorphic Connections} &&
\pageref{sn: smooth} \\
\ref{sn: ss+mm}. \text{ Symplectic Structure and Moment Map} &&
\pageref{sn: ss+mm} \\
\ref{sn: mm is sp}. \text{ The Monodromy Map is Symplectic} &&
\pageref{sn: mm is sp} \\
\ref{sn: imds}. \text{ Isomonodromic Deformations and Symplectic Fibrations} &&
\pageref{sn: imds} \\
\text{Appendix} && 
\pageref{apx} \\
\text{Acknowledgments and References} &&
\pageref{acknow} 
\end{array}
$$
\end{centering}

\ 

\

\begin{equation*}\fbox{$
\begin{array}{cccccc}
&&&&&\\
&&  \wt\M(\ba) & \mapright{\cong} & \wt\A_{\fl}(\ba)/\G_1 &\\
&&  \bigcup && \mapdown{\cong}&\\
\hphantom{\rightarrow}\wt O_1\times\cdots\times\wt O_m\spq G &\cong&
  \wt\M^*(\ba) & \mapright{\wt\nu} & \wt M_0(\ba) &\hphantom{\rightarrow}\\
&&&&&\\
\end{array}
$}
\end{equation*}

$$\text{ Main commutative diagram (p.\pageref{main diagram})}$$

\newpage

\section{Introduction} \label{sn: intro}

Moduli spaces of representations of fundamental groups of Riemann
surfaces have been intensively studied in recent years and have
an incredibly rich structure:
For example, a theorem of 
Narasimhan and Seshadri \cite{NarSes} identifies 
the space of irreducible unitary representations of the
fundamental group of a compact Riemann surface 
with the moduli space of stable
holomorphic vector bundles on the surface.
In particular, this description puts a K\"ahler structure on the space
of fundamental group representations---it has a 
symplectic structure together with a compatible complex
structure. 
A remarkable fact is that although the complex structure on the space
of representations will depend on the complex structure of the
surface, the {\em symplectic} structure only depends on the topology,
a fact often referred to as `the symplectic nature of the fundamental
group' \cite{Gol84}.

The geometry
is richer still if we consider the moduli space
of {\em complex} fundamental group representations:
Due to results of Hitchin, Donaldson and Corlette,
the K\"ahler structure above now becomes a {\em hyper-K\"ahler}
structure and the symplectic structure becomes a {\em complex
symplectic} structure, which is still topological.
One of the main aims of this paper is to generalise this complex symplectic
structure. 
(Hyper-K\"ahler structures will not be considered here.)

First recall that, over a Riemann surface,
there is a one to one correspondence between complex fundamental group
representations and {\em holomorphic} connections (obtained by taking
a holomorphic connection to its monodromy/holonomy representation).
Then replace the word `holomorphic' by `meromorphic'---we
will study the symplectic geometry of moduli spaces
of meromorphic connections.

In fact, as in the holomorphic case, these moduli spaces may also be
realised in a more topological way, using 
a generalised notion of monodromy data. 
By restricting a meromorphic connection to the complement of its polar
divisor and taking the corresponding monodromy representation, a map
is obtained from the moduli space of meromorphic connections to the
moduli space of representations of the fundamental group of the
punctured Riemann surface.
For connections with only simple poles 
this map is generically a covering map and
so we are essentially in the well-known case of
representations of fundamental groups of punctured Riemann
surfaces.
However in general there are local moduli at the poles---it
is not sufficient to restrict to the complement of the polar divisor
and take the monodromy representation as above. 

Fortunately this extra data---the local moduli of meromorphic
connections---has been studied in the theory of differential equations
for many years and has a monodromy-type description in terms
of `Stokes matrices', which encode (as we will explain) the change in
asymptotics of solutions on sectors at the poles.
The Stokes matrices and the fundamental group representation fit
together in a natural way and
the main question we ask is simply:
``What is the symplectic geometry of
these moduli spaces of generalised monodromy data?''

Recently, Martinet and Ramis \cite{MR91} have constructed a huge group
associated to any Riemann surface, the 
`wild fundamental group', whose set of finite dimensional representations
naturally corresponds to the set of meromorphic connections on the
surface. 
Although we will not directly use this perspective,
the question above can then be provocatively rephrased as asking:
``What is the symplectic nature of the {\em wild} fundamental group?''

The motivation behind these questions is to
understand intrinsically the symplectic geometry of
the full family of isomonodromic deformation
equations of Jimbo, Miwa and Ueno \cite{JMU81}.
The initial impetus was the theorem of B.Dubrovin \cite{Dub95} 
identifying
the local moduli space of semisimple Frobenius manifolds 
with a space of Stokes matrices. 
(In brief, this means certain Stokes matrices parameterise 
certain two-dimensional topological quantum field theories.)
The original aim was to find a more intrinsic approach to the 
intriguing (braid group invariant) Poisson
structure written down by Dubrovin 
on this space of Stokes matrices 
in the rank three case (see \cite{Dub95} appendix F and also the
recent paper \cite{Ugag} of M.Ugaglia for the higher rank formula). 
The key step in the proof of Dubrovin's theorem is that
(in the semisimple case) the WDVV equations 
(of Witten-Dijkgraaf-Verlinde-Verlinde)
are {\em equivalent} to the equations for isomonodromic deformations of
certain meromorphic connections
on the Riemann sphere with just two poles, of
orders one and two respectively (so the space of solutions corresponds
to the moduli of such connections---the Stokes matrices).

More generally
Jimbo, Miwa and Ueno \cite{JMU81} have written down
a vast family of nonlinear differential equations, governing
isomonodromic deformations of meromorphic connections over $\IP^1$
having arbitrarily many poles of arbitrary order 
(on arbitrary rank bundles).
These are of independent interest
and can be thought of as a {\em universal} family of nonlinear 
equations: They are the largest family of nonlinear
differential 
equations known to have the `Painlev\'e property' (that, except on fixed
critical varieties, solutions will only have poles as singularities).
Special cases include the six Painlev\'e equations 
(which arise as the isomonodromic deformation equations for 
connections on rank two bundles over $\IP^1$, with total pole
multiplicity four) and Schlesinger's equations 
(the simple pole case---see below).

In brief, the six Painlev\'e equations were found almost a
hundred years ago, as a means to construct new transcendental
functions (namely their solutions---the Painlev\'e transcendents);
R.Fuchs discovered then that the sixth Painlev\'e equation arises as an
isomonodromic deformation equation. 
The subject then lay more or less
dormant until the late 1970's when (spectacularly) 
Wu, McCoy, Tracey and Barouch
\cite{WMTB,MTW}
found that the correlation functions of certain quantum field
theories satisfied Painlev\'e equations.
Subsequently Jimbo, Miwa, M\^{o}ri and Sato \cite{SMJ, JMMS} 
showed that this was a special
case of a more general phenomenon and developed the theory of
`holonomic quantum fields' which led to \cite{JMU81}.
See for example \cite{Iwetal, Ume98} for more background material.

One expects isomonodromic deformations to lurk underneath most
integrable partial differential equations since
the heuristic `Painlev\'e integrability test'
says that a nonlinear PDE
will be `integrable' if it admits some reduction to an ODE with
the Painlev\'e property
(for example the KdV equation has a reduction to the first Painlev\'e
equation and all six Painlev\'e equations appear as reductions of the
anti-self-dual Yang-Mills equations; 
see \cite{AC,MasWoo}).
They certainly appear in a diverse range of nonlinear problems 
in geometry and theoretical
physics, such as Frobenius manifolds \cite{Dub95}
or in the construction of Einstein metrics \cite{Tod,Hit-tei}.

On the other hand general solutions of isomonodromy 
equations cannot be given explicitly in terms of known special
functions; as mentioned above
general solutions 
are {\em new} transcendental functions (see \cite{Ume89}).
This is the reason we turn to geometry to understand more about 
these equations.
Recent work on isomonodromic deformations 
seems to have focused mainly on particular examples, 
in particular exploring the rich geometry
of the six Painlev\'e equations and searching for the few, very special,
explicit 
solutions that they do admit.
The question we address here is simply 
``What is the symplectic geometry 
of the full family of isomonodromic
deformation equations of Jimbo, Miwa and Ueno?''

Geometrically the isomonodromy equations constitute
a flat (Ehresmann) connection on a fibre bundle
over a space of deformation parameters, the fibres being certain
moduli spaces
of meromorphic connections over $\IP^1$.
Thus the idea is to find natural symplectic structures on such moduli
spaces and then prove they are preserved by the isomonodromy equations.
In certain special cases, such as the Schlesinger or six Painlev\'e
equations the symplectic geometry is well-known (see \cite{HW,Hit95,Oka}).
The main results of this paper are analogous to those of Hitchin \cite{Hit95}
and Iwasaki \cite{Iwa92} who explained intrinsically why
Schlesinger's equations and certain rank two higher genus isomonodromy
equations (respectively) admit a time-dependent Hamiltonian description. 
However for the general isomonodromy equations considered here a Hamiltonian 
description is still not known: this work indicates strongly that such
description should exist. (See also Remark \ref{new features}.)

For example,
understanding the symplectic geometry of the isomonodromic deformation
equations enables us to
ask questions about their quantisation. 
This has been addressed in certain cases by Reshetikhin \cite{Resh92} 
and Harnad \cite{Harn96} and leads to Knizhnik-Zamolodchikov type 
equations.

A key step (Theorem \ref{thm mmap is sp})
is to establish that the transcendental map taking a
meromorphic connection to its (generalised) monodromy data, is 
a symplectic map.
This is the `inverse monodromy theory' version of
the well-known result in inverse scattering
theory, that the map from the set of initial potentials to
scattering data is a symplectic map (see \cite{FadTak} Part 1, Chapter III).

Although apparently not mentioned in the literature, 
a useful perspective (explained in Section \ref{sn: imds}) has 
been to interpret the paper \cite{JMU81}
of Jimbo, Miwa and Ueno, as stating that the 
Gauss-Manin connection in non-Abelian cohomology (in the sense of
Simpson \cite{Sim94ab}) generalises 
to the case of {\em meromorphic} connections.
This offers a fantastic guide for future generalisation.

\subsection*{The Prototype: Simple Poles}

The simplest way to explain the results of this paper is to
first describe the intrinsic symplectic geometry of
Schlesinger's equations, following Hitchin \cite{Hit95}.

Choose matrices 
$A_1,\ldots,A_m\in\End(\IC^n)$, distinct numbers
$a_1,\ldots, a_m\in\IC$ and
consider the following meromorphic
connection on the trivial rank $n$ holomorphic vector bundle over the
Riemann sphere:
\begin{equation} \label{nabla}
\nabla := d - \left(A_1\frac{dz}{z-a_1}+\cdots+A_m\frac{dz}{z-a_m}\right).
\end{equation}
This has a simple pole at each $a_i$ and will have no further
pole at $\infty$ if and only if
$A_1+\cdots+A_m = 0$, which we will assume to be the case.
Thus, on removing a small open disc $D_i$ from around each $a_i$ and 
restricting $\nabla$ to the $m$-holed sphere
${{S}}:=\IP^1\setminus(D_1\cup\cdots\cup D_m)$,
we obtain a (nonsingular)
holomorphic connection.
In particular it is flat and so, taking its monodromy, a
representation of the fundamental group of the $m$-holed sphere is
obtained.
This procedure defines a holomorphic map, which we will call the
{\em monodromy map}, from the set of such connection 
coefficients to the set of complex
fundamental group representations:
\begin{equation} \label{mon map}
\bigl\{ (A_1,\ldots,A_m) \ \bigl\vert \ \text{$\sum A_i= 0$} \bigr\}
\ \mapright{\nu_{\bf a}}\ 
\bigl\{ (M_1,\ldots,M_m) \ \bigl\vert \ M_1\cdots M_m = 1 \bigr\}
\end{equation}
where appropriate loops generating the fundamental group
of ${{S}}$ have been chosen 
and the matrix $M_i\in G:=GL_n(\IC)$ is the
automorphism obtained by parallel translating a basis of solutions
around the $i$th loop.

This map is the key to the whole theory and is generically a local
analytic isomorphism.
It is tempting to think of $\nu_{\bf a}$ as a generalisation of the
exponential function, but note the dependence on the pole positions
$\bf a$ 
is rather complicated since the monodromy map {\em solves} Painlev\'e type
equations (see below).

We can however study the geometry of the monodromy map, particularly
the {\em symplectic} geometry.
First, to remove the base-point dependence, quotient both sides of 
(\ref{mon map}) by the diagonal conjugation action of $G$.
Secondly restrict the matrices $A_i$ to be in fixed adjoint orbits.
(These may be
identified with coadjoint orbits using the trace, and so have natural 
complex symplectic structures.)
Thus we pick $m$ generic (co)adjoint orbits $O_1,\ldots,O_m$
and require $A_i\in O_i$.
Also define $\cC_i\subset G$ to be the conjugacy class
containing $\exp(2\pi \sqrt{-1} A_i)$ for any $A_i\in O_i$.
Fixing $A_i\in O_i$ implies $M_i\in \cC_i$.
The key fact now is that the sum $\sum A_i$ is
a moment map for the diagonal conjugation
action of $G$ on $O_1\times \cdots\times O_m$ and so (\ref{mon map}) becomes
\begin{equation} \label{mon map'}
O_1\times\cdots\times O_m\spq G
\quad \mapright{\nu_{\bf a}}\quad 
\Hom_{\bf\cC}\bigl(\pi_1({{S}}),G\bigr)/
G
\end{equation}
where the subscript ${\bf\cC}$ 
means we restrict to representations 
having local monodromy around $a_i$ in the conjugacy class $\cC_i$.
The symplectic geometry of this set of representations has
been much studied recently.
The primary symplectic description is due to Atiyah
and Bott \cite{AB82,Ati90} 
and involves interpreting it as an
infinite dimensional symplectic quotient, starting with all 
$C^\infty$ connections on the manifold-with-boundary ${{S}}$
(see also \cite{Aud95}).
Alternatively, 
a purely finite dimensional description of the symplectic structure is
given by the cup product in parabolic group cohomology
\cite{BK93} and 
finding a finite dimensional proof of the closedness of this symplectic
form has occupied many people. (See \cite{Kar92,FR,AMR,GHJW,AMM}.) 

By construction the left-hand side of (\ref{mon map'}) 
is a finite dimensional symplectic quotient and
one of the key results of \cite{Hit95} was that, for any
choice of pole positions ${\bf a}$, 
the  monodromy map $\nu_{\bf a}$ is {\em symplectic}; 
it pulls back the Atiyah-Bott symplectic structure on the right to the
symplectic structure on the left, coming from the coadjoint orbits.
This fact is the key to understanding
intrinsically why Schlesinger's equations are symplectic, as
we will now explain.

Observe that 
if we vary the positions 
of the poles slightly 
then the
spaces on the left and the right of  (\ref{mon map'}) do not change.
However the monodromy map $\nu_{\bf a}$ does vary.
Schlesinger \cite{Schles} wondered  
how the matrices $A_i$
should vary with respect to the pole positions $a_1,\ldots,a_m$ 
such that the monodromy
representation $\nu_{\bf a}(A_1,\ldots,A_m)$ stays fixed, and
thereby discovered the beautiful family
of nonlinear differential equations which now bear his name:
$$
\frac{\partial A_i}{\partial a_j}=\frac{[A_i,A_j]}{a_i-a_j}\qquad 
\text{if } i\ne j, \quad \text{ and} \qquad 
\frac{\partial A_i}{\partial a_i}=-\sum_{j\ne
i}\frac{[A_i,A_j]}{a_i-a_j}.
$$

These are the equations for {\em isomonodromic deformations} of the
logarithmic connections $\nabla$ on $\IP^1$ that we began
with in (\ref{nabla}).
Hitchin's observation now is that the local self-diffeomorphisms of the
symplectic manifold $O_1\times\cdots \times O_m\spq G$ induced by
integrating Schlesinger's equations, are clearly symplectic
diffeomorphisms,
because they are of the form 
$\nu^{-1}_{{\bf a'}} \circ \nu_{{\bf a}}$
for two sets of pole positions 
${\bf a}$ and ${\bf a'}$ and the monodromy map is
a local {\em symplectic} isomorphism for all ${\bf a}$.

This is the picture we will generalise to the case of higher order
poles, after rephrasing it in terms of symplectic fibrations. 
The main missing ingredient is the Atiyah-Bott
construction of a symplectic structure on the generalised monodromy
data;
when the discs are removed any local moduli at the poles is lost.
We will work throughout over $\IP^1$ since 
the weight of this paper is to see what to do
locally at a pole of order at least two, and because 
our main interest is the Jimbo-Miwa-Ueno
isomonodromy equations, which are for connections over $\IP^1$. 
However, apart from in Section \ref{sn: triv} and for the explicit form of
the isomonodromy equations, global coordinates on $\IP^1$ are not used
and so most of this work generalises immediately to arbitrary
genus compact Riemann surfaces, possibly with boundary.

The organisation of this paper is as follows.
The next three sections each give a different approach to meromorphic
connections.
In Section \ref{sn: triv} we generalise the left-hand side of 
(\ref{mon map'}) and prove the results we
will need later regarding the symplectic
geometry of these spaces.
Section \ref{sn: gm} describes the generalised monodromy data of a
meromorphic connection on a Riemann surface, both the local data
(the Stokes matrices) and the global data fitting together the local data
at each pole. 
This generalises the spaces of fundamental group
representations above---the notion of fixing the conjugacy class of
local monodromy is replaced by fixing the `formal equivalence class'.
In Section \ref{sn: smooth} we introduce an appropriate notion of
$C^\infty$ singular connections and prove the basic results one might
guess from the non-singular case, 
relating flat singular connections to spaces of monodromy data 
and to meromorphic connections on degree zero holomorphic bundles.
The notion of fixing the formal type
of a meromorphic connection corresponds nicely to the notion of
fixing the `$C^\infty$ Laurent expansion' of a 
flat $C^\infty$ singular connection.
Section \ref{sn: ss+mm} then shows that the Atiyah-Bott symplectic
structure generalises naturally to these spaces of $C^\infty$ singular
connections, and that as in the non-singular case,
the curvature, when defined appropriately, 
is a moment map for the gauge group action. 
Thus the spaces of generalised monodromy data also appear as infinite
dimensional symplectic quotients.
(One should note that the `naive' extension of the Atiyah-Bott symplectic
structure to $C^\infty$ singular connections does not work 
since the two-forms that arise are too singular to be 
integrated.)
Section \ref{sn: mm is sp} summarises the preceding sections in a 
commutative diagram and then proves the key result, that
the monodromy map pulls back the Atiyah-Bott type symplectic structure
on the generalised monodromy data to the explicit symplectic structure of 
Section \ref{sn: triv}, on the spaces of meromorphic connections.
Section \ref{sn: imds} explains geometrically what the
Jimbo-Miwa-Ueno isomonodromy equations are (we will write them
explicitly in the appendix), and 
then proves
the main result, Theorem \ref{thm2}, 
that the isomonodromic deformation equations of
Jimbo-Miwa-Ueno are equivalent to a
flat {\em symplectic} connection on a {\em symplectic} fibre bundle 
having the moduli spaces of Section \ref{sn: triv} as fibre.
Note that in the general case there are
more deformation parameters: we may vary the `irregular types' of the
connections at the poles, as well as the pole positions. 
This produces, in particular, many nonlinear symplectic braid group
representations on the spaces of monodromy data.
Finally, we end by sketching a relationship between
Stokes matrices and Poisson Lie groups.

\section{Meromorphic Connections on Trivial Bundles} \label{sn: triv}

Let $D=k_1(a_1)+\cdots+k_m(a_m) > 0$ be an effective divisor
on $\IP^1$ (so that $a_1,\ldots,a_m\in\IP^1$ are distinct points 
and $k_1,\ldots, k_m > 0$ are positive integers)
and let $V\to\IP^1$ be a rank $n$ holomorphic vector 
bundle.
\begin{defn}	\label{dfn: mero conn}
A {\em meromorphic connection} $\nabla$ on $V$ with poles on $D$ is a map
$\nabla: V\rightarrow V\otimes K(D)$
from the sheaf of holomorphic sections of $V$ to the sheaf of sections of 
$V\otimes K(D)$,
satisfying the Leibniz rule:
$\nabla(fv)=(df)\otimes v+ f\nabla v$, 
where $v$ is a local section of $V$, $f$ is a local holomorphic
function and $K$ is the sheaf of holomorphic one-forms on $\IP^1$.
\end{defn}
If we choose a local coordinate $z$ on $\IP^1$ vanishing
at $a_i$ then in terms of a local trivialisation of $V$, 
$\nabla$ has the form
$\nabla=d-{A}$,
where 
\begin{equation}	\label{eqn: mero one-form}
{A}=A_{k_i}\frac{dz}{z^{k_i}}+\cdots+
	A_1\frac{dz}{z}+A_0dz+\cdots 
\end{equation}
is a matrix of meromorphic one-forms and $A_j\in\End(\IC^n), j\le k_i$.

\begin{defn}	\label{generic defn}
A meromorphic connection $\nabla$ will be said to be {\em generic} if 
at each $a_i$ the leading coefficient $A_{k_i}$ is 
diagonalisable with distinct eigenvalues (if $k_i\ge 2$), or
diagonalisable with distinct eigenvalues mod $\IZ$ (if $k_i=1$).
\end{defn}
This condition is independent of the trivialisation and
coordinate choice.
We will restrict to such generic connections since they are
simplest yet sufficient for our purpose (to describe the symplectic
nature of isomonodromic deformations).

\begin{defn}	\label{cmp framing}
A {\em compatible framing at $a_i$} of a vector bundle $V$ with
generic connection $\nabla$ 
is an isomorphism $g_0:V_{a_i}\to\IC^n$ 
between the fibre $V_{a_i}$ and $\IC^n$
such that 
the leading coefficient of $\nabla$ is {\it diagonal} in any
local trivialisation of $V$ extending $g_0$.
\end{defn}
Given a trivialisation of $V$ in a
neighbourhood of $a_i$ so that $\nabla=d-A$ as above,
then a compatible framing is represented by a constant matrix
(also denoted $g_0$)
that diagonalises the leading coefficient:
$g_0\in G$ such that $g_0A_{k_i}g_0^{-1}$
is diagonal.

At each point $a_i$ choose a germ $d- \iAo$ of a diagonal 
generic meromorphic
connection on the trivial
rank $n$ vector bundle. 
(We use the terminology that {\em a} trivial bundle is just
trivialisable, but {\em the} trivial bundle has a chosen
trivialisation. 
Also pre-superscripts $\iA$, whenever used, will signify local
information near $a_i$.) 
Thus $\iAo$ is a matrix of germs of meromorphic one-forms, which we
require (without loss of generality) to be diagonal.
If $z_i$ is a local coordinate vanishing at $a_i$, write
\begin{equation}	\label{Q and lambda}
\iAo = d(\iQ) + \iLo\frac{dz_i}{z_i}
\end{equation}
where $\iLo$ is a constant diagonal matrix and $\iQ$ is a diagonal
matrix of meromorphic functions.
\begin{defn}
A connection $(V,\nabla)$ with compatible framing $g_0$ at $a_i$ has
{\em irregular type} $\iAo$ if $g_0$ extends to a formal
trivialisation of $V$ at $a_i$, in which
$\nabla$ differs from $d-\iAo$ by a matrix of one-forms with just
simple poles.
\end{defn}
Equivalently this means, if $\nabla=d-A$ in some 
local trivialisation, we require
$gAg^{-1}+(dg)g^{-1}= d(\iQ) + \iL dz_i/z_i$ for some
diagonal matrix $\iL$  not necessarily equal to $\iLo$ 
and some formal bundle automorphism 
$g\in G\flb z_i \frb=GL_n(\IC\flb z_i \frb)$ with
$g(a_i)=g_0$.
The diagonal matrix $\iL$ appearing here will be referred to as the 
{\em exponent of formal monodromy} of $(V,\nabla,g_0)$.  

Let $\ba$ denote the choice of the effective divisor $D$ and all the germs
$\iAo$.
The spaces which generalise those on the left-hand side of 
(\ref{mon map'}) are defined as follows.

\begin{defn}	\label{dfn: pp mfd}
The moduli space 
$\M^*(\ba)$ is the set of isomorphism classes of pairs
$(V,\nabla)$ where $V$ is a
{\em trivial} rank $n$ holomorphic vector bundle over $\IP^1$
and $\nabla$ is a meromorphic connection on $V$
which is formally equivalent
to $d-\iAo$ at $a_i$ for each $i$ and has
no other poles.
\end{defn}
Following \cite{JMU81}, we also define slightly larger moduli spaces:
\begin{defn} \label{dfn: eppm}
The {\em extended moduli space} 
$\wt\M^*(\ba)$ is the set of isomorphism classes of triples
$(V,\nabla,{\bf g})$ consisting of a generic connection
$\nabla$ (with poles on $D$) 
on a {\em trivial} holomorphic vector bundle $V$ over $\IP^1$
with compatible framings 
${\bf g}=({^1\!g_0},\ldots,{^m\!g_0})$
such that
$(V,\nabla,{\bf g})$ has irregular type
$\iAo$
at each $a_i$.
\end{defn}
The term `extended moduli space' is taken from the paper \cite{Jef94}
of L.Jeffrey, since these spaces play a similar role (but are not the same). 

\begin{rmk}	\label{rmk: deg zero spaces}
For use in later sections we also define spaces $\M(\ba)$ and
$\wt\M(\ba)$ simply by replacing the word `trivial' by `degree zero' in
Definitions \ref{dfn: pp mfd} and \ref{dfn: eppm} respectively.
\end{rmk}

Since $\M^*(\ba)$ and $\wt\M^*(\ba)$
are moduli spaces of connections on trivial bundles we can
obtain explicit descriptions of them. First define 
$G_{k}$ to be the group of $(k-1)$-jets of bundle automorphisms:
$$G_{k}:= GL_n\bigl(\IC[\ze]/\ze^{k}\bigr)$$
where $\ze$ is an indeterminate.
Then the main result of this section is:
\begin{prop} \label{prop: main pp}

$\bullet$ $\M^*(\ba)$ is isomorphic to a complex symplectic quotient
\begin{equation} \label{eqn: mod sp}
\quad\M^*(\ba)\cong O_1\times\cdots\times O_m\spq G\quad
\end{equation}
where $G:= GL_n(\IC)$ and $O_i\subset \gkis$ is a coadjoint orbit of
$G_{k_i}$.

$\bullet$ Similarly there are complex symplectic manifolds (extended orbits)
$\wt O_i$ with $\dim(\wt O_i)=\dim( O_i) + 2n$ and (free) Hamiltonian
$G$ actions, such that
\begin{equation} \label{eqn: ext mod sp}
\quad\wt\M^*(\ba)\cong \wt O_1\times\cdots\times \wt O_m\spq G.
\end{equation}

$\bullet$ In this way $\wt\M^*(\ba)$ inherits (intrinsically) the
structure of a complex symplectic manifold, the torus actions changing
the choices of compatible framings are Hamiltonian (with moment maps
given by the values of the $\iL$'s) and $\M^*(\ba)$
arises as a symplectic quotient by these $m$ torus actions.
\end{prop}
Because of the third statement here (and that $\M^*(\ba)$ may not be Hausdorff)
we will mainly work with the extended moduli spaces.
They will be the phase spaces of the isomonodromy equations.
Before proving Proposition \ref{prop: main pp}  
we first collect together all the results we will need regarding the
extended orbits $\wt O_i$.
\subsection*{Extended Orbits}
Fix a positive integer $k\ge 2$.
Let $B_k$ be the subgroup of $G_k$ of elements having constant term
$1$. This is a unipotent normal subgroup and in fact $G_k$ is the
semi-direct product  $G\ltimes B_k$ (where 
$G:=GL_n(\IC)$ acts on $B_k$ by conjugation).
Correspondingly the Lie algebra of $G_k$ decomposes as a vector space
direct sum and dualising we have:
\begin{equation}	\label{eqn: co-sum}
\gks=\lbks\oplus \gs.
\end{equation}
Concretely if we have a matrix of meromorphic one-forms $A$ as in 
(\ref{eqn: mero one-form})  with $k=k_i$ then the principal part
of $A$ can be identified as an element of $\gks$ simply by replacing
the coordinate $z$ by the indeterminate $\ze$:
\begin{equation} 	\label{eqn: pp map}
A_{k}\frac{d\ze}{\ze^{k}}+\cdots+ A_1\frac{d\ze}{\ze}\in
\gks.
\end{equation}
Abusing notation, this element of $\gks$
will also be denoted by $A$. 
Such $A$'s are identified as
elements of $\gks$ via the pairing 
$\langle A,X\rangle := \res_0(\tr(A(\ze)\cdot X))=
\sum_{i=1}^{k}\tr(A_iX_{i-1})$
where 
$X=X_0+X_1\ze+\cdots+X_{k-1}\ze^{k-1}\in \gk$.
Then from (\ref{eqn: co-sum}), 
$\lbks$ is identified with the set of $A$ in (\ref{eqn: pp map})
having zero residue and $\g^*$ with those having only a residue term 
(zero irregular part).
Let $\pir:\gks\to\gs$ and $\pii:\gks\to\lbks$ denote the corresponding
projections.

Now 
choose a diagonal element 
$A^0=
A^0_kd\ze/\ze^k+\cdots+A^0_{2}d\ze/\ze^{2}$
of $\lbk^*$ whose leading coefficient $A^0_k$ has distinct eigenvalues.
For example if $k=k_i$, such $A^0$ arises from the irregular
part $d(\iQ)$ of $\iAo$ in (\ref{Q and lambda}).
Let $O_{\!B}\subset \lbks$ denote the $B_k$ coadjoint orbit containing $A^0$.

\begin{defn}	\label{dfn: ext orb}
The {\em extended orbit} $\wt O\subset G \times \gks$ associated to $O_{\!B}$ is:
$$\wt O := 
\left\{
(g_0,A)\in G \times \gks \ \bigl\vert
\ \pii(g_0 A g_0^{-1})\in O_{\!B} \right\} $$ 
where $\pii:\g_k^*\to\lbks$ is the natural projection removing the residue.
\end{defn}
If $(g_0,A)\in\wt O$ then eventually $A$ 
will correspond to the principal part of a generic
meromorphic connection and $g_0$ to a compatible framing.

\begin{lem}
The extended orbit $\wt O$ is canonically isomorphic to the symplectic
quotient of the product $T^*G_k\times O_{\!B}$ by $B_k$, where both the
cotangent bundle $T^*G_k$ and the coadjoint orbit $O_{\!B}$ have their natural
symplectic structures.
\end{lem}
\pf
$B_k$ acts by the coadjoint action on $O_{\!B}$ and by the 
standard (free) left action on $T^*G_k$ (induced
from left multiplication of the groups).
A moment map is given by
$\mu: T^*G_k\times O_{\!B}\to \lbks; \ 
	(g,A,B)\mapsto -\pii(\Ad^*_g(A))+B$
where $B\in O_{\!B}$ and $(g,A)\in G_k\times\gks\cong T^*G_k$ via the left
trivialisation.
Thus
\begin{equation}	\label{eqn: prezero}
\mu^{-1}(0)=\left\{ (g,A,B) \ \bigr\vert\ \pii(gAg^{-1})=B \right\}.
\end{equation}
It is straightforward to check that the map
\begin{equation}	\label{eqn: chi}
\chi:\mu^{-1}(0)\to \wt O;\qquad (g,A,B)\mapsto (g(0),A)
\end{equation}
is well-defined, surjective and has precisely the $B_k$ orbits as fibres.
\epf 

This gives $\wt O$ the structure of a complex symplectic manifold.
Next we examine the torus action on $\wt O$ corresponding to changing
the choice of compatible framing.
If $(g_0,A)\in\wt O$ then by hypothesis there is some $g\in G_k$ such
that  $gAg^{-1} = A^0 + \Lambda d\ze/\ze$ for some matrix $\Lambda$.
It is easy now to modify $g$ such that $\Lambda$ is in fact
diagonal. (Conjugating by $1+X\ze^{k-1}$ for an appropriate matrix $X$
will remove any off-diagonal part of $\Lambda$.)
It follows that there is a well-defined map:
$$ \mu_T: \wt O \to \lt^*; \qquad (g_0,A)\mapsto -\Lambda\frac{d\ze}{\ze},$$
where, as above, if $R=\Lambda{d\ze}/{\ze}\in\lt^*$ and $\Lambda'\in \lt$ then
$\langle R, \Lambda' \rangle= \tr( \Lambda \Lambda')$.

\begin{lem}	\label{lem: taction}
1) The map $\mu_T$ is a moment map for 
the  free action of $T\cong(\IC^*)^n$ on $\wt O$ defined by
$t(g_0,A)=(tg_0,A)$ where $t\in T$.

2) The symplectic quotient at the value $-R$ of $\mu_T$ is
the $G_k$ coadjoint orbit through the element 
$A^0+R$ of $\gks$.

3) Any tangents $v_1,v_2$ to $\wt O\subset G\times\gks$ at $(g_0,A)$
are of the form
$$v_i=(X_i(0), [A,X_i] + g_0^{-1} \dot R_i g_0)\in\g\times\gks$$
for some $X_1,X_2\in \gk$ and $ \dot R_1,\dot R_2\in\lt^*$ 
(where $\g\cong T_gG$ via 
left multiplication), and the symplectic structure on $\wt O$ is then given
explicitly by the formula:
\begin{equation}	\label{extorb formula}
\omega_{\wt O}(v_1,v_2) =   \langle \dot R_1,\wt X_2 \rangle
- \langle \dot R_2,\wt X_1 \rangle 
+ \langle A,[X_1,X_2] \rangle
\end{equation}
where $\wt X_i:=g_0X_ig_0^{-1}\in \gk\text{ for $i=1,2$.}$
\end{lem}
\pf
There is a surjective `winding' map $w:G_k\times\lt^*\to\wt O$ defined by
$(g,R)\mapsto (g(0),g^{-1}(A^0+R)g)$. It fits into the
commutative diagram:
\begin{equation}	\label{cd: extorb}
\begin{array}{ccccccc}
G_k\times\lt^* & \hookmapright{\iota} & \mu^{-1}(0) & \subset & 
T^*G_k\times O_{\!B} & \mapright{\prsmall} & T^*G_k \\
\shortmapdown{w} &   & \shortmapdown{\chi}  &   &   &   &   \\
\wt O & = & \wt O &   &   &   &   
\end{array}
\end{equation}
where $\chi$ is from (\ref{eqn: chi}), $\pr$ is the projection and 
$\iota(g,R) := (g,g^{-1}(A^0+R)g,A^0)$.  
Since the $O_{\!B}$ component of $\iota$ is constant the pullback of the symplectic
structure on $T^*G_k$ along $\pr\circ\iota$ is the pullback of the
symplectic structure on $\wt O$ along $w$.
Let $T$ act on $T^*G_k$ by the standard left action $t(g,A)=(tg,A)$,
on $O_{\!B}$ by conjugation ($t(B)=tBt^{-1}$) and on $G_k\times\lt^*$ by
left multiplication: $t(g,R)=(tg,R)$. Observe that all the maps
in (\ref{cd: extorb}) are then $T$-equivariant and that a moment map on
$T^*G_k$ is given by 
$\nu: T^*G_k\to\lt^*;\ (g,A)\mapsto
-\delta(\pir(gAg^{-1}))$
since the map $\delta \circ \pir$ (taking the diagonal part of the residue term
of an element of $\gks$) is the dual of the derivative of the
inclusion $T\hookrightarrow G_k$.
Statement 1) now
follows from the fact that the pullback of $\nu$ along $\pr\circ
\iota$ is the pullback of $\mu_T$ along $w$ (both maps pullback to the
projection $G_k\times\lt^*\to \lt^*$).
The third statement is proved by directly calculating the pullback of the 
symplectic structure on $T^*G_k$ along $\pr\circ\iota$. 
(Note $(X_i,\dot R_i)$  is 
just a lift of $v_i$ to $G_k\times \lt^*$.) 
The second statement follows directly from (\ref{extorb formula}).
\epf

Thus, by projecting to $\gks$, we see 
 $\wt O$ is a principal $T$ bundle over an $n$-parameter family of $G_k$
coadjoint orbits. An alternative description will also be useful:
\begin{lem}[Decoupling]		\label{lem: decoupling}
The following map is a symplectic isomorphism:
$$\wt O \cong T^*G \times O_{\!B}; \quad
(g_0,A)\mapsto(g_0,\pir(A), \pii(g_0Ag_0^{-1}) )$$
where $T^*G\cong G\times\gs$ via the left trivialisation.
\end{lem}
\pf
It is an isomorphism as the map
$(g_0,S,B)\mapsto (g_0,g_0^{-1}Bg_0+S)\in \wt O$ 
(where $(g_0,S,B)\in T^*G \times O_{\!B}$) is an inverse.
Under this identification, a section $s$ of the projection $\chi$ in 
(\ref{cd: extorb}) is given by: $s:T^*G\times O_{\!B}\to T^*G_k\times
O_{\!B};$
$(g_0,S,B)\mapsto (g_0, g_0^{-1} B g_0 + S, B)$
where left multiplication is used to trivialise the cotangent bundles.
A straightforward calculation shows $s$ is symplectic.
\epf

This will be important because $O_{\!B}$ admits {\em global} Darboux coordinates.

\begin{cor}	\label{cor: G action}
The free $G$ action $h(g_0,A):=(g_0h^{-1},hAh^{-1})$ on $\wt O$ 
is Hamiltonian with moment map
$\mu_{G}: \wt O \to \gs; \ (g_0,A)\mapsto \pir(A).$
\end{cor}
\pf
After decoupling $\wt O,\ G$ acts only on the $T^*G$ factor and it
does so by the standard action coming from right multiplications,
which has moment map $\mu_G$.
\epf 

Finally in the simple pole case ($k=1$) not yet considered we define 
$$
\wt O := \left\{ (g_0,A)\in G\times \gs
 \ \bigl\vert \ g_0Ag_0^{-1}\in \lt' \right\}\subset G\times \gs
$$
where $\lt'\subset \lt^*$ is the subset containing diagonal
matrices whose eigenvalues are distinct mod $\IZ$.
If we identify $G\times \gs$ with $T^*G$ then $\wt O$ is in fact a
symplectic submanifold (see \cite{GS} Theorem 26.7).
The formula (\ref{extorb formula}) holds unchanged and 
the free $G$ and $T$ actions are still Hamiltonian with the same moment maps
as above (the diagonalisation of $A$ used to define $\mu_T$ is simply
$g_0Ag_0^{-1}$).
Note that
the winding map $w:G\times\lt'\to \wt O;$ 
$(g_0,R) \mapsto(g_0,g_0^{-1}Rg_0)$ is now an isomorphism.

\ 

\pfms[of Proposition \ref{prop: main pp}].\ 
Choose a coordinate $z$ to identify $\IP^1$ with $\IC\cup\infty$
such that each $a_i$ is finite. 
Define $z_i:=z-a_i$.
The chosen meromorphic connection germs $d-\iAo$ determine $G_{k_i}$
coadjoint orbits $O_i$ and extended orbits $\wt O_i$ as above: 
Define $O_i$ to be the coadjoint orbit through the
point of $\gkis$ determined (using the coordinate choice $z_i$)
by the principal part of $\iAo$ in (\ref{Q and lambda}).
Similarly the irregular part of $\iAo$ determines a point of $\lbkis$ 
and $\wt O_i$ is the extended orbit associated to the $B_{k_i}$
coadjoint orbit through this point.

Now suppose $\nabla$ is a meromorphic connection on a holomorphically 
trivial bundle
$V$ over $\IP^1$ with poles on the divisor $D$.
Upon  trivialising $V$ we find $\nabla=d-A$ for a matrix $A$ of
meromorphic one-forms of the form
\begin{equation}	\label{global A}
A=\sum_{i=1}^{m}\Big(
\iA_{k_i}\frac{dz}{(z-a_i)^{k_i}}+\cdots+\iA_{1}\frac{dz}{(z-a_i)}\Big)
\end{equation}
where the $\iA_j$ are $n\times n$ matrices.
The principal part of $A$ at $a_i$ determines an element $\iA\in
\gkis$ as above (replacing $z-a_i$ by $\ze$ in the $i$th term of the sum
(\ref{global A})).

The crucial fact now is that $\nabla$ is formally equivalent to
$d-\iAo$ at $a_i$ if and only if $\iA$ is in $O_i$.
The `only if' part is clear since the gauge action
restricts to the coadjoint action on the principal parts of $A$.
The converse is not true in general (even if formal
{\em meromorphic} transformations are allowed: see \cite{BV83}),
but it does hold in the generic case we are considering here, and is
well-known (see \cite{BJL79}).
Also, using the description of the extended orbits  as principal $T$
bundles, it follows that
if $\nabla$ is generic and has compatible framings 
${\bf g}=({^1\!g_0},\ldots,{^m\!g_0})$ then $(V,\nabla,{\bf g})$ has
irregular type $\iAo$ at $a_i$ if and only if 
$(\igo,\iA)$ is in $\wt O_i$.

Thus any meromorphic connection on the trivial bundle with the correct formal
type determines and is determined by a point of the product 
$O_1\times\cdots\times O_m$. 
Observe however that a general point of $O_1\times\cdots\times O_m$
will give a connection with an additional pole at $z=\infty$ unless we
impose the constraint
\begin{equation}	\label{sum of residues}
{^1\!A_{1}}+\cdots+{^m\!A_{1}}=0.
\end{equation}
Also observe that the choice of global trivialisation of $V$
corresponds to the action 
of $G$ on $O_1\times\cdots\times O_m$ by diagonal conjugation.
The first statement in Proposition \ref{prop: main pp} follows
simply by observing that 
the left-hand side of (\ref{sum of residues})
is a moment map for this $G$ action on $O_1\times\cdots\times O_m$
(since the $G$ action on each $O_i$ factor
is  the restriction of the
coadjoint action to $G\subset G_{k_i}$).

The proof of the second statement is completely analogous. (The $G$
action on $\wt O_i$ is given in Corollary \ref{cor: G action}.)

Lemma \ref{lem: decoupling} makes it transparent that $\wt\M^*(\ba)$
is a smooth complex manifold: the symplectic quotient by $G$ just removes a
factor of $T^*G$ from the product of extended orbits.
It is straightforward to check  
the complex symplectic structure so defined on $\wt\M^*(\ba)$ is
independent of the coordinate choices used above.
(In fact arbitrary {\em local} coordinates $z_i$ may be used.)

Finally the statements concerning the torus actions are immediate from
Lemma \ref{lem: taction}, since the $G$ and $T$ action on each extended
orbit  commute. 
\epfms

\begin{rmk} 
Open subsets of 
the symplectic quotients 
$
O_1\times\cdots\times O_m \spq G
$
in (\ref{eqn: mod sp}) have been previously studied:
They are
algebraically completely integrable Hamiltonian systems \cite{AHH,Bea}.
See also \cite{DonMar} Sections 4.3 and 5.3.
The perspective there is to regard these as spaces 
of meromorphic Higgs fields, rather
than as spaces of meromorphic connections.
\end{rmk}
\input{smid-monod}

\input{smid-smooth}

\input{smid-mmsp}
\input{smid-imds}

\ 

{ \small
{\bf Acknowledgments.}\ 	\label{acknow}
The results of this paper were announced
at the 1998 ICM in Berlin and
appeared in my Oxford D.Phil. thesis \cite{Boa}, which was supported by an
E.P.S.R.C. grant. 
I would like to thank my D.Phil. supervisor
Nigel Hitchin for the guidance and inspiration 
provided throughout this project, 
Boris Dubrovin whose work introduced me to the geometry of Stokes matrices,
and Hermann Flaschka whose lectures at the
1996 isomonodromy conference in Luminy helped germinate the idea that the
main results of this paper might be true.
Many other mathematicians have also influenced this work through helpful
comments, 
and I am grateful to them all, 
especially 
Mich\`ele Audin, John Harnad, Alastair King, Mich\`ele Loday-Richaud, 
Eyal Markman, Nitin Nitsure, Claude Sabbah, Graeme Segal
and Nick Woodhouse.
}

\renewcommand{\baselinestretch}{1}              %
\normalsize
\bibliographystyle{amsplain}    \label{biby}
\bibliography{../thesis/syr} 
\ 

\ 

\noindent
Journal ref: Adv. Math. 163, (2001) 137-205

\end{document}

%% file: macros.tex
\long\def\pb #1*/{}

\def\reE@DeclareMathSymbol#1#2#3#4{%
    \let#1=\undefined
    \DeclareMathSymbol{#1}{#2}{#3}{#4}}
\DeclareSymbolFont{symbolsC}{U}{txsyc}{m}{n}
\SetSymbolFont{symbolsC}{bold}{U}{txsyc}{bx}{n}
\DeclareFontSubstitution{U}{txsyc}{m}{n}
\reE@DeclareMathSymbol{\strictiff}{\mathrel}{symbolsC}{76}

\newcommand\beq{\begin{equation}}
\newcommand\eeq{\end{equation}}
\newcommand\bal{\begin{align*}}
\newcommand\eal{\end{align*}}   
\newcommand\bmx{\left(\begin{matrix}}
\newcommand\emx{\end{matrix}\right)}
\newcommand\bsmx{\left(\begin{smallmatrix}}
\newcommand\esmx{\end{smallmatrix}\right)}
\newcommand\bmxnp{\begin{matrix}}
\newcommand\emxnp{\end{matrix}}
\newcommand\bsmxnp{\begin{smallmatrix}}
\newcommand\esmxnp{\end{smallmatrix}}

\newcommand{\wt}{\widetilde}
\newcommand{\wA}{{\wt \A}}

\DeclareMathSymbol{\widehatsym}{\mathord}{largesymbols}{"62}
\newcommand\lowerwidehatsym{%
  \text{\smash{\raisebox{-1.3ex}{%
    $\widehatsym$}}}}
\newcommand\fixwidehat[1]{%
  \mathchoice
    {\accentset{\displaystyle\lowerwidehatsym}{#1}}
    {\accentset{\textstyle\lowerwidehatsym}{#1}}
    {\accentset{\scriptstyle\lowerwidehatsym}{#1}}
    {\accentset{\scriptscriptstyle\lowerwidehatsym}{#1}}
}

\newcommand{\wh}{\fixwidehat}

\newcommand{\ip}{\int_{\IP^1}}

\newcommand{\pir}{\pi_{\hbox{\footnotesize\rm res}}}
\newcommand{\pii}{\pi_{\hbox{\footnotesize\rm irr}}}

\newcommand{\spq}{/\!\!/}

\providecommand{\deg}{\text{\rm deg}}

\def\part#1{\frac{\partial\phantom{q}}{\partial#1}}

\newcommand {\flb}{\lbrack\!\lbrack}
\newcommand {\frb}{\rbrack\!\rbrack}

\newcommand{\mki}{{k_i}}                

\newcommand{\hh}{\mathop{\rm H}}

\newcommand{\jmu}{{\mathop{\rm\footnotesize JMU}}}

\newcommand{\Lia}{{\mathop{\rm Lie}}}     

\DeclareMathOperator{\ISto}{{\IS}to} 
 
\newcommand{\Syst}{\mathop{\rm Syst}}

\newcommand{\Sect}{\text{\rm Sect}}

\newcommand{\Ssect}{\wh {\mathop{\rm Sect}}}
\newcommand{\Roots}{\mathop{\rm Roots}}

\newcommand{\Mult}{\mathop{\rm Mult}}

\newcommand{\asexp}{{\mathop{\rm \AE}}}   
\newcommand{\const}{{\mathop{\rm Const}}}
\newcommand{\Ad}{{\mathop{\rm Ad}}}

\newcommand{\PP}{{\rm PP}}

\DeclareMathOperator{\pr}{pr}
\newcommand{\prsmall}{\mathop{\footnotesize\rm pr}}
\newcommand{\res}{{\mathop{\rm Res}}}
\newcommand{\Res}{{\mathop{\rm Res}}}

\newcommand{\tr}{{\mathop{\rm Tr}}}

\DeclareMathOperator{\Hom}{Hom}         

\renewcommand{\ker}{\mathop{\rm Ker}}
\DeclareMathOperator{\End}{End}

\newcommand{\diag}{{\mathop{\rm diag}}}

\newcommand{\fl}{{\mathop{\footnotesize\rm {f}l} }}   
\newcommand{\od}{{\mathop{\footnotesize\rm od}}}          
\newcommand{\vrt}{{\mathop{\footnotesize\rm Vert}}}       
\newcommand{\hol}{{\mathop{\footnotesize\rm hol}}}        
\newcommand{\rel}{{\mathop{\footnotesize\rm Rel}}}

\newcommand{\ba}{{\bf a}}

\newcommand{\bap}{{\bf A\!'}}

\newcommand{\bd}{{\bf d}}

\newcommand{\bs}{{\bf S}}

\DeclareSymbolFont{bbold}{U}{bbold}{m}{n}
\DeclareSymbolFontAlphabet{\mathbbold}{bbold}

\newcommand{\IA}{\mathbb{A}}

\newcommand{\IC}{\mathbb{C}}
\newcommand{\ID}{\mathbb{D}}

\newcommand{\IP}{\mathbb{P}}                                     
                           
\newcommand{\IR}{\mathbb{R}}                           
\newcommand{\IS}{\mathbb{S}}

\newcommand{\IZ}{\mathbb{Z}}

\newcommand{\A}{\mathcal{A}}

\newcommand{\cC}{\mathcal{C}}

\newcommand{\wcc}{\widetilde \cC}       

\newcommand{\F}{\mathcal{F}}

\newcommand{\G}{\mathcal{G}}
\newcommand{\cH}{\mathcal{H}}

\newcommand{\M}{\mathcal{M}}

\newcommand{\cO}{\mathcal{O}}

\newcommand{\cS}{\mathcal{S}}
\newcommand{\T}{\mathcal{T}}
\newcommand{\cT}{\mathcal{T}}

\newcommand{\g}{       \mathfrak{g}     }
\newcommand{\gk}{      \mathfrak{g}_k   }

\newcommand{\gs}{      \mathfrak{g}^*       }
\newcommand{\gks}{     \mathfrak{g}_k^*   }
\newcommand{\gkis}{    \mathfrak{g}_{k_i}^* }

\newcommand{\lb}{       \mathfrak{b}            }
\newcommand{\lbk}{      \mathfrak{b}_k          }

\newcommand{\lbks}{     \mathfrak{b}_k^*        }
\newcommand{\lbkis}{    \mathfrak{b}_{k_i}^*    }

\newcommand{\lt}{\mathfrak{t}}

\newcommand{\gl}{       \mathfrak{gl}     } 

\newcommand{\al}{\alpha}

\newcommand{\be}{\beta}
\newcommand{\ga}{\gamma}

\newcommand{\Ga}{\Gamma}

\newcommand{\si}{\sigma}

\newcommand{\Sif}{\Sigma_i(\wh F)}
\newcommand{\Sjf}{\Sigma_j(\wh F)}
\newcommand{\ze}{\zeta}

\renewcommand{\bar}{\overline}

\newcommand{\iA}{{^i\!A}}
\newcommand{\iQ}{{^i\!Q}}
\newcommand{\iL}{{^i\!\Lambda}}
\newcommand{\bL}{{\bf\Lambda}}
\newcommand{\idLj}{{^i\!\dot\Lambda_j}}
\newcommand{\idRj}{{^i\!\dot R_j}}
\newcommand{\iLp}{{^i\!\Lambda'}}
\newcommand{\iLo}{{^i\!\Lambda^{\!0}}}
\newcommand{\Lo}{{\Lambda^{\!0}}}
\newcommand{\iwbj}{{^i\wt b_j}}
\newcommand{\iwbo}{{^i\wt b_0}}
\newcommand{\isj}{{^i\si_j}}
\newcommand{\iso}{{^i\si_0}}
\newcommand{\iAo}{{^i\!A^0}}
\newcommand{\iwAo}{{^i\!\wt A^0}}
\newcommand{\iwAod}{{^i\!\wt A^0_\Delta}}
\newcommand{\owAod}{{^1\!\wt A^0_\Delta}}

\newcommand{\oa}{{^0\!\!\A}}
\newcommand{\oaf}{{^0\!\!\A_{\mathop{\footnotesize\rm {f}l}}}}    
\newcommand{\igo}{{^i\!g_0}}
\newcommand{\iwg}{{^i\wh g}}            

\newcommand{\owg}{{^1\wh g}}            

\newcommand{\og}{{^0\!\G}}

\makeatletter
 \newlength{\typesize}
\makeatother

\newlength{\vvoff}
\newlength{\hhoff}

\def\mapright#1{\smash{
        \mathop{\longrightarrow}\limits^{#1}}}

\def\hookmapright#1{\smash{
        \mathop{\hookrightarrow}\limits^{#1}}}
\def\mapdown#1{\Big\downarrow
        \rlap{$\vcenter{\hbox{$\!\scriptstyle#1$}}$}}
\def\shortmapdown#1{\big\downarrow
        \rlap{$\vcenter{\hbox{$\scriptstyle#1$}}$}}
\def\shortmapup#1{\big\uparrow
        \rlap{$\vcenter{\hbox{$\scriptstyle#1$}}$}}

\def\underset#1#2{\ \smash{\mathop{ #2 }\limits_{#1}}\ }

\newcommand{\pf}{\begin{bpf}}

\newcommand{\pfms}{\begin{bpfms}}
\newcommand{\epf}{\end{bpf}\hfill$\square$\\}           
\newcommand{\epfms}{\end{bpfms}\hfill$\square$\\}       

\newcommand{\idea}{\begin{bidea}}

\newcommand{\eidea}{\end{bidea}\hfill$\square$\\}           

\newcommand{\sk}{\begin{bsk}}    

\newcommand{\esk}{\end{bsk}\hfill$\square$\\}           
\newcommand{\sketch}{\begin{bsketch}}

\newcommand{\esketch}{\end{bsketch}\hfill$\square$\\}



%% file: macros-thm1.1.smid.tex
\theoremstyle{plain}  

\newtheorem{thm}{Theorem}[section]

\newtheorem{lem}{Lemma}[section]
\newtheorem{cor}[lem]{Corollary}
\newtheorem{prop}[lem]{Proposition}
\newtheorem{claim}{Claim}[section]

\theoremstyle{definition} \newtheorem{rmk}{Remark}[section]\newtheorem {defn}{Definition}[section]

%% file: extimdstex.big.pstex_t
\begin{picture}(0,0)%
\includegraphics{extimdstex.big.pstex}%
\end{picture}%
\setlength{\unitlength}{2565sp}%
\begingroup\makeatletter\ifx\SetFigFont\undefined%
\gdef\SetFigFont#1#2#3#4#5{%
  \reset@font\fontsize{#1}{#2pt}%
  \fontfamily{#3}\fontseries{#4}\fontshape{#5}%
  \selectfont}%
\fi\endgroup%
\begin{picture}(5802,5610)(1411,-7159)
\put(4726,-3511){\makebox(0,0)[lb]{\smash{{\SetFigFont{14}{16.8}{\rmdefault}{\mddefault}{\updefault}{\color[rgb]{0,0,0}$\wt\nu$}%
}}}}
\put(1501,-4636){\makebox(0,0)[lb]{\smash{{\SetFigFont{14}{16.8}{\rmdefault}{\mddefault}{\updefault}{\color[rgb]{0,0,0}$\wt\M^*$}%
}}}}
\put(5551,-7036){\makebox(0,0)[lb]{\smash{{\SetFigFont{14}{16.8}{\rmdefault}{\mddefault}{\updefault}{\color[rgb]{0,0,0}$\wt X$}%
}}}}
\put(1951,-7036){\makebox(0,0)[lb]{\smash{{\SetFigFont{14}{16.8}{\rmdefault}{\mddefault}{\updefault}{\color[rgb]{0,0,0}$\wt X$}%
}}}}
\put(5251,-4636){\makebox(0,0)[lb]{\smash{{\SetFigFont{14}{16.8}{\rmdefault}{\mddefault}{\updefault}{\color[rgb]{0,0,0}  $\wt M$}%
}}}}
\put(1426,-2236){\makebox(0,0)[lb]{\smash{{\SetFigFont{14}{16.8}{\rmdefault}{\mddefault}{\updefault}{\color[rgb]{0,0,0}$\wt\M^*(\ba)$}%
}}}}
\put(5251,-2236){\makebox(0,0)[lb]{\smash{{\SetFigFont{14}{16.8}{\rmdefault}{\mddefault}{\updefault}{\color[rgb]{0,0,0}  $\wt M(\ba)$}%
}}}}
\end{picture}%

%% file: smid-monod.tex
\section{Generalised Monodromy}	\label{sn: gm}

This section describes the monodromy
data of a generic meromorphic connection, 
involving both a fundamental
group representation and Stokes matrices,
largely following \cite{BV83, BJL79, JMU81, L-R94, MR91}.
The presentation here is quite nonstandard however and 
care has been taken to keep track of all the choices 
made and thereby see what is intrinsically defined.

Fix the data $\ba$ of 
a divisor $D=\sum k_i(a_i)$ on $\IP^1$ and connection germs
$d-\iAo$ at each $a_i$ as in Section \ref{sn: triv}.
Let $(V,\nabla,{\bf g})$ be a compatibly framed
meromorphic connection on a holomorphic vector bundle $V\to\IP^1$
with irregular type $\ba$.

In brief, the monodromy data arises as follows. 
The germ $d-\iAo$ canonically determines some 
directions at $a_i$ (`anti-Stokes directions') for each $i$ and 
(using local coordinate choices) we can consider the sectors at each
$a_i$ bounded by these directions (and having some small fixed radius). 
Then the key fact is that the framings ${\bf g}$ 
(and a choice of branch of logarithm at each pole)
determine, in a canonical way, a
choice of basis of solutions of the
connection $\nabla$ on each such sector at each pole.
Now along any path in the punctured sphere 
$\IP^1\setminus\{a_1.\ldots,a_m\}$ between
two such sectors we can extend the two corresponding bases of
solutions and obtain a constant $n\times n$ matrix relating these two bases.
The monodromy data of $(V,\nabla,{\bf g})$ is simply the set of all
such constant matrices, plus the exponents of formal monodromy.

Before filling in the details of this procedure we will give a
concrete definition of the monodromy manifolds that store this
monodromy data and so give a clear idea of where we are going.
All of the monodromy manifolds are of the following form.
Suppose $N_1,\ldots,N_m$  are manifolds, we have maps 
$\rho_i:N_i\to G'$ to some group $G'$ for each $i$ and 
that there is an action of $G=GL_n(\IC)$ on $G'$ 
(via group automorphisms) and on each $N_i$
such that $\rho_i$ is $G$-equivariant.
Define a map ${\bf{\rho}}$ to be the (reverse ordered) product of
the $\rho_i$'s:
$$
{\bf{\rho}} : N_1\times\cdots\times N_m \to G';\quad
(n_1,\ldots,n_m)\mapsto \rho_m(n_m)\cdots\rho_2(n_2)\rho_1(n_1).
$$
Since $G$ acts on $G'$ by automorphisms, $\bf{\rho}$ is
$G$-equivariant and ${\bf \rho}^{-1}(1)$ is a $G$-invariant 
subset of the product
$N_1\times\cdots\times N_m$. 
We will write the quotient as:
\begin{equation}	\label{eqn: msq}
N_1\times\cdots\times N_m \spq G := {\bf \rho}^{-1}(1)/G.
\end{equation}
This is viewed simply as a convenient way to write down the
various sets of monodromy data that arise\footnote{
The relationship with \cite{AMM} will be discussed elsewhere.
}.
There is no conflict of notation since 
the symplectic quotients of Section \ref{sn: triv} arise in this way 
by taking $N_i= O_i$ (or $\wt O_i$), $G'=(\g^*,+)$ and the $\rho_i$'s
as the moment maps for the $G$ actions.
All of the examples in this section however will have $G':= G$ acting on
itself by conjugation.

Recall in the simple pole case that we fixed generic coadjoint
orbits $O_1,\ldots,O_m$ of $G$ to define a symplectic space of
connections on trivial bundles over $\IP^1$.
By choosing appropriate generators of the fundamental group of the
punctured sphere we see that the corresponding space of
monodromy data is of the above form:
\begin{equation} \label{eqn: homc}
\Hom_{\bf\cC}
\bigl(\pi_1(\IP^1\setminus\{a_1,\ldots,a_m\}),G\bigr)/G \cong
\cC_1\times\cdots\times\cC_m\spq G
\end{equation}
where $G$ acts on each conjugacy class $\cC_i$ by conjugation 
and each map $\rho_i:\cC_i\to G$ is just the inclusion. 

Considering higher order poles in Section \ref{sn: triv}
amounted to replacing the coadjoint
orbits of $G$ above by coadjoint orbits of $G_{k_i}$ (still denoted $O_i$)
or by extended orbits $\wt O_i$.
By analogy, in this section we now replace each conjugacy class
in the simple pole case by a
larger manifold $\cC_i$ (the multiplicative version of $O_i$), or by
$\wt \cC_i$  (the multiplicative version of $\wt O_i$).
The basic definition is somewhat surprising:

\begin{defn} \label{def: gcc}

$\bullet$ Let  $U_{+/-}$ be the upper/lower triangular unipotent subgroups of
$G$, then
$$\wt \cC_i := G\times (U_+\times U_-)^{k_i-1}\times \lt$$
where $\lt$ is the set of diagonal $n\times n$ matrices and
$k_i$ is the pole order at $a_i$.
(If $k_i=1$ replace $\lt$ by $\lt'$ here; the elements 
with distinct eigenvalues mod $\IZ$.)
A point of $\wt \cC_i$ will be denoted 
$({C_i}, {^i\bs}, \iLp)$ where
$^i\bs=({^i\!S_1},\ldots,{^i\!S_{2k_i-2}})\in U_+\times U_-\times U_+
\times \cdots \times U_-$.

$\bullet$ 
The formula
 $t({C_i}, {^i\bs}, \iLp) =
\Bigl(t\cdot{C_i},\ 
(t\cdot{^i\!S_1}\cdot t^{-1},\ldots,t\cdot{^i\!S_{2k_i-2}}\cdot t^{-1}),
\  \iLp\Bigr)$
defines a free action of the torus $T$ on $\wt \cC_i$
and, given some fixed choice of $\iLp$,
we define $\cC_i$ to be the subset of the
quotient having $\lt$ component fixed to this value:
$\quad\cC_i := (\wcc_i/T)\vert_\iLp.$

$\bullet$ The map $\rho_i:\wcc_i\to G$ is defined by the formula
$$\rho_i(C_i,{^i\bs},\iLp)= 
C_i^{-1}\cdot{^i\!S_{2k_i-2}}\cdots{^i\!S_2}\cdot{^i\!S_1}\cdot
\exp\bigl((2\pi\sqrt{-1})\iLp\bigr)\cdot C_i.
$$%
(This is $T$ invariant so descends to define $\rho_i:\cC_i\to G$.)

$\bullet$ Finally $G$ acts on $\wcc_i$ (and on $\cC_i$) via 
$g(C_i,{^i\bs},\iLp)= (C_i g^{-1},{^i\bs},\iLp)$
(so that $\rho_i$ is clearly $G$-equivariant, where $G$ acts on itself
by conjugation).
\end{defn}
The triangular matrices ${^i\!S_j}$ (with $1$'s on their diagonals)
appearing here are the Stokes
matrices.
Note that in every case 
the dimension of $\wt \cC_i$ is the same as the dimension of
the extended orbit $\wt O_i$ (and similarly $\dim (\cC_i)= \dim (O_i)$). 
Also note that if the pole is simple ($k_i=1$) then $\wcc_i= G\times\lt'$
and that $\cC_i$ can naturally be identified with the conjugacy class
through $\exp(2\pi\sqrt{-1}\cdot \iLp) \in G$.

Our aim in the rest of this section is to define an (abstract) space of
monodromy data $M(\ba)$ and an intrinsic holomorphic map $\nu$ from
the moduli space $\M^*(\ba)$ of Section \ref{sn: triv} to $M(\ba)$,
obtained by taking monodromy data. We will call $\nu$ the monodromy
map,
although the names Riemann-Hilbert map or 
de Rham morphism are also appropriate.
Recall in Proposition \ref{prop: main pp} that after making some
choices (of local coordinates in that case) a concrete description of 
the moduli space $\M^*(\ba)$ was obtained. 
Analogously here, after making some choices
(of some `tentacles'; something like a choice of generators of the
fundamental group---see Definition \ref{def: tentacles}), 
we will see that the quotient $\cC_1\times\cdots\times\cC_m\spq G$ 
is a concrete realisation of the
space of monodromy data.
Thus we will have the diagram:
\begin{equation}	\label{diag: }
\begin{array}{ccc}
O_1\times\cdots\times O_m\spq G &  & 
\cC_1\times\cdots\times\cC_m\spq G \\
\shortmapup{\cong} &   & \shortmapup{\cong}  \\
\M^*(\ba) & \mapright{\nu} & M(\ba).    
\end{array}
\end{equation}
As in Section \ref{sn: triv} we will work mainly with the extended version
(putting tildes on all the spaces in the above diagram) since the spaces
are then manifolds (and again the non-extended version 
may be deduced by considering torus actions).

\begin{lem}	\label{lem: sdim}
The extended monodromy manifold 
$\wt M(\ba)\cong \wcc_1\times\cdots\times\wcc_m\spq G$ is indeed a 
complex manifold and has the same dimension as $\wt \M^*(\ba)$.
\end{lem}
\pf
Remove the $G$ action by fixing $C_1=1$, so
$\wt M(\ba)$ is identified with the subvariety 
$\rho_m\cdots\rho_1 = 1$ of
the product
$\wcc'_1\times\wcc_2\times\cdots\times\wcc_m$
where
$\wcc'_1$ is the subset of $\wcc_1$ having $C_1=1$.
The result now follows from the implicit function theorem since
the map 
$\rho_m\cdots\rho_1: 
\wcc'_1\times\wcc_2\times\cdots\times\wcc_m \to G$
is surjective on tangent vectors (except in the trivial case $m=1,
k_1=1$).
In particular $\dim\wt M(\ba) =\sum\dim(\wcc_i)-2n^2$ 
and, from Proposition \ref{prop: main pp},
this is $\dim\wt\M^*(\ba)$.
\epf
\subsection*{Local Moduli: Stokes Matrices} \label{ssn: local moduli}
First we will set up the necessary auxiliary data.
Let $d- A^0$ be a diagonal generic meromorphic
connection on the trivial rank $n$ vector bundle over the unit disc 
$\ID\subset \IC$ with
a pole of order $k\ge 2$ at $0$ and no others.
Let $z$ be a coordinate on $\ID$ vanishing at $0$.
Thus 
(as in Section \ref{sn: triv})
$A^0 = dQ + \Lo\frac{dz}{z}$
where $\Lo$ is a constant diagonal matrix and $Q$ is a diagonal
matrix of meromorphic functions.
$Q$ is determined by $A^0$ and $z$ by requiring that it has constant term
zero in its Laurent expansion with respect to  $z$.
Write $Q=\diag(q_1,\ldots,q_n)$ and
define $q_{ij}(z)$ to be the leading term of $q_i-q_j$.
Thus if $q_i-q_j=a/z^{k-1}+b/z^{k-2}+\cdots$ then $q_{ij}=a/z^{k-1}$.

Let the circle $S^1$ parameterise rays (directed lines) emanating from 
$0\in\IC$; 
intrinsically one can think of this circle as being the boundary circle
of the real oriented blow up of $\IC$ at the origin.
If $d_1,d_2\in S^1$ then $\Sect(d_1,d_2)$ will denote the (open) 
sector swept out
by rays rotating in a positive sense from $d_1$ to $d_2$.
The radius of these sectors will be taken 
sufficiently small when required later.

\begin{defn}	\label{dfn: asds}
The {\it anti-Stokes directions} $\IA\subset S^1$ 
are the directions $d\in S^1$ such that
for some $i\ne j$:
$q_{ij}(z)\in\IR_{<0} \text{ for $z$ on the ray specified by } d.$
\end{defn} 
These are the directions along which $e^{q_i-q_j}$ decays most rapidly
as $z$ approaches $0$ and it follows that $\IA$ is independent of the
coordinate choice $z$.
(For uniform notation later, define $\IA$ to contain a
single, arbitrary direction if $k=1$; this will be used only as a
local branch cut.)
\begin{defn}	\label{dfn: roots}
Let $d\in S^1$ be an anti-Stokes direction.

$\bullet$
The {\em roots} of $d$ are the ordered pairs
$(ij)$ `supporting' $d$:
$$\Roots(d):=\{ (ij) \ \bigl\vert \ 
q_{ij}(z)\in\IR_{<0} \text{ along } d \}.$$

$\bullet$
The {\em multiplicity} $\Mult(d)$ of $d$ is the number of roots
supporting $d$.

$\bullet$ 
The {\em group of Stokes factors}
associated to $d$ is the group
$$\ISto_d(A^0):=\{ K \in G \ \bigl\vert \ (K)_{ij}=\delta_{ij}
\text{  unless $(ij)$ is a root of $d$} \}. $$	
\end{defn}
It is straightforward to check that $\ISto_d(A^0)$
is a unipotent subgroup of $G=GL_n(\IC)$ of dimension equal to the
multiplicity of $d$.
Beware that the terms `Stokes factors' and `Stokes matrices' are used
in a number of different senses in the literature.
(Our terminology is closest to Balser, Jurkat and Lutz \cite{BJL79}. 
However our approach is perhaps closer to that of Martinet and Ramis
\cite{MR91} but what we call Stokes factors, they
call Stokes matrices, 
and they do not use the things we call Stokes matrices.)

The anti-Stokes directions $\IA$ have $\pi/(k-1)$ rotational symmetry: if
$q_{ij}(z)\in\IR_{<0}$ then
$q_{ji}(z\exp(\frac{\pi\sqrt{-1}}{k-1}))\in\IR_{<0}$.
Thus the number $r:=\#\IA$ of anti-Stokes directions is divisible by
$2k-2$ and we define $l:= r/(2k-2)$.
We will refer to an $l$-tuple $\bd=(d_1,\ldots,d_l)\subset \IA$ 
of {\em consecutive} anti-Stokes directions as a `half-period'.
When weighted by their multiplicities, 
the number of anti-Stokes 
directions in any half-period is $n(n-1)/2 = \Mult(d_1)+\cdots+\Mult(d_l)$.
Also a half-period $\bd$ determines a total ordering of the set 
$\{q_1,\ldots,q_n\}$ defined by:
\begin{equation}	\label{order defn}
q_i\underset{{\bd}}{<}q_j \qquad\Longleftrightarrow\qquad 
(ij) \text{ is a root of some }d\in\bd.
\end{equation}
A simple check shows $(ij) \text{ is a root of some }d\in\bd$ 
precisely if  $e^{q_i}/e^{q_j}\to 0$ as $z\to 0$ along the ray  
$\theta({\bd})\in S^1$ bisecting $\Sect(d_1,d_l)$ (so that 
(\ref{order defn}) is 
the natural dominance ordering along $\theta(\bd)$). 
In turn there is a permutation matrix $P\in G$ associated to $\bd$
given by 
$(P)_{ij} = \delta_{\pi(i)j}$ where $\pi$ is the permutation of
$\{1,\ldots,n\}$ corresponding to (\ref{order defn}):
$q_i\underset{{\bd}}{<}q_j \Leftrightarrow 
\pi(i)<\pi(j)$. 
A key result is then:

\begin{lem} \label{lem: sf to sm}
Let $\bd=(d_1,\ldots,d_l)\subset \IA$ be a half-period (where
$d_{i+1}$ is the next anti-Stokes direction after $d_i$ in 
a positive sense).

$1)$ The product of the corresponding groups of Stokes factors is
isomorphic as a variety, via the product map, to the subgroup of $G$ conjugate
to $U_+$ via $P$:
$$ \prod_{d\in\bd}\ISto_d(A^0) \cong PU_+P^{-1};\quad
(K_1,\ldots,K_l)\mapsto K_l\cdots K_2K_1\in G
$$

$2)$ Label the rest of $\IA$ uniquely as $d_{l+1},\ldots,d_r$
(in order) then the following map from the product of all the groups
of Stokes factors, is an isomorphism of varieties:
$$\prod_{d\in \IA}\ISto_d(A^0) \cong (U_+\times U_-)^{k-1};\quad
(K_1,\ldots,K_r)\mapsto (S_1,\ldots S_{2k-2})$$
where
$S_i:=P^{-1}K_{il}\cdots K_{(i-1)l+1}P\in U_{+/-}$
if $i$ is odd/even.
\end{lem}
\pf
$1)$ holds since the groups of Stokes factors are a set of `direct
spanning' subgroups of $PU_+P^{-1}$; see Borel \cite{Bor91} Section
14. It is also proved directly in Lemma 2 of \cite{BJL79}.
$2)$ follows from $1)$ simply by observing that the orderings
associated to neighbouring half-periods are opposite.
\epf
Now we move on to the local moduli of meromorphic connections.
Let $\Syst(A^0)$ denote the set of germs at $0\in\IC$ 
of meromorphic connections on 
the trivial rank $n$ vector bundle, that are formally
equivalent to $d-A^0$.
Concretely we have
$$\Syst(A^0) = \{ \ d-A \ \bigl\vert \  A = \wh F[A^0] \text{ for some }
\wh F \in G\flb z \frb \ \}
$$
where $A$ is a matrix of germs of meromorphic one-forms, 
$G\flb z \frb = GL_n(\IC\flb z \frb)$ and 
$\wh F[A^0] = (d\wh F)\wh F^{-1} + \wh FA^0\wh F^{-1}$.
The group $G\flb z \frb$ does not act on $\Syst(A^0)$ since in general
$\wh F[A^0]$ will not have convergent entries. The subgroup 
$G\{ z \}:= GL_n(\IC\{ z\})$ of germs of holomorphic bundle
automorphisms does act however and we wish to study the quotient
$\Syst(A^0)/G\{ z\}$ which is by definition the set of isomorphism
classes of germs of meromorphic connections formally equivalent to
$A^0$. Note that any {\em generic} connection is formally
equivalent to some such $A^0$.

In the Abelian case ($n=1$) and in the simple pole case ($k=1$)
$\Syst(A^0)/G\{ z\}$ is just a point; the notions of formal and
holomorphic equivalence coincide.
In the non-Abelian, irregular case ($n,k\ge 2$)
however,  
$\Syst(A^0)/G\{ z\}$ is non-trivial and we will explain how to
describe it explicitly in terms of Stokes matrices.

It is useful to consider spaces slightly larger than $\Syst(A^0)$:
\begin{defn}

$\bullet$
Let $\wh\Syst_{\rm cf}(A^0)$ be the set of compatibly framed connection
germs with both irregular and formal type $A^0$.

$\bullet$
Let $\wh\Syst_{\rm mp}(A^0):=\{
(A,\wh F)\ \bigl\vert\ A\in\Syst(A^0), 
\wh F\in G\flb z \frb, A = \wh F[A^0] \},$
be the set of {\em marked pairs}.
\end{defn}
Thus $\wh\Syst_{\rm cf}(A^0)$ is the set of pairs $(A,g_0)$ with
$A\in\Syst(A^0)$ and $g_0\in G$, such that $g_0[A]$ and $A^0$  
have the same leading term.
Clearly the projection to the first factor is a surjection 
$\wh\Syst_{\rm cf}(A^0)\twoheadrightarrow \Syst(A^0)$ and the  fibres are the
orbits of the torus action $t(A,g_0)= (A,t\cdot g_0)$ 
(where $t\in T\cong (\IC^*)^n$).

\begin{lem} \label{lem: get F}
There is a canonical isomorphism 
$\wh\Syst_{\rm mp}(A^0)\cong\wh\Syst_{\rm cf}(A^0)$: 
For each compatibly framed connection germ $(A,g_0)\in\wh\Syst_{\rm cf}(A^0)$
there is a {\em unique} formal isomorphism $\wh F\in G\flb z \frb$
with $A = \wh F[A^0]$ and $\wh F(0) = g_0^{-1}$.
\end{lem}
\pf
It is sufficient to prove that if $\wh F[A^0]=A^0$ then $\wh F\in T$,
since this implies the map $(A,\wh F)\mapsto (A,g_0)$ with 
$g_0:=\wh F(0)^{-1}$, is bijective.
Now if $\wh F[A^0]=A^0$ then the $(ij)$ matrix entry $f$ of $\wh F$
is a power series solution to 
$df=(d(q_i-q_j) + ({\lambda_i-\lambda_j})dz/z)f$, where
$\lambda_i=(\Lambda^0)_{ii}$.
It follows (using Definition \ref{generic defn} if $k=1$)
that $f=0$ unless $i=j$ when it is a constant. 
\epf
See for example \cite{BV83} for an algorithm to determine $\wh F$ from $g_0$.
Below we will use `$\wh\Syst(A^0)$' to denote either of these two sets.
Heuristically the action $g(A,\wh F)=(g[A],g\circ\wh F)$ of $G\{ z\}$ 
on the marked pairs is free 
and so one expects the quotient
$$\cH(A^0) := \wh{\Syst}(A^0)/G\{ z \}$$
to be in some sense nice (as is indeed the case). 
Moreover the actions of $T$ and 
$G\{ z\}$ on $\wh{\Syst}(A^0)$ commute so 
$\Syst(A^0)/G\{ z \}\cong \cH(A^0)/T$.

The fundamental technical result we need to quote in order 
to describe $\cH(A^0)$
is the following theorem.
First we set up a labelling convention, 
that will behave well when we vary $A^0$ in later sections.
Choose a point $p$ in one of the $r$ sectors at $0$ bounded by
anti-Stokes rays.
Label the first anti-Stokes ray when turning in a positive sense from
$p$ as $d_1$ and label the subsequent rays $d_2,\ldots,d_r$ in turn.
Write $\Sect_i := \Sect(d_i,d_{i+1})$; the `$i$th sector'
(indices are taken modulo $r$).
Note $p\in\Sect_r=\Sect_0$; the `last sector'.
Also define the `$i$th supersector' to be
$\Ssect_i := \Sect\left(d_i-\frac{\pi}{2k-2},d_{i+1}+
\frac{\pi}{2k-2}\right).$
This is a sector containing the $i$th sector
symmetrically (the same direction bisects both)
and has opening greater than $\pi/(k-1)$.
(The rays bounding these supersectors are usually referred to as
`Stokes rays'.)   
\begin{thm} 	\label{thm: multisums}
Suppose $\wh F\in G\flb z \frb$ is a formal transformation such that
$A:=\wh F[A^0]$ has convergent entries.
Set the radius of the sectors $\Sect_i, \Ssect_i$ to be less than 
the radius of convergence of $A$.
Then the following hold:

$1)$
On each sector $\Sect_i$
there is a canonical way to choose an  invertible 
$n\times n$ matrix of holomorphic functions $\Sif$ such that 
$\Sif [A^0] =  A$.

$2)$
$\Sif$ can be analytically continued to the supersector $\Ssect_i$ 
and then $\Sif$ is asymptotic to
$\wh F$ at $0$ within $\Ssect_i$.

$3)$
If $g\in G\{z\}$ and $t\in T$ then 
$\Sigma_i( g \circ \wh F \circ t^{-1})=g\circ \Sif \circ t^{-1}. $
\end{thm}
The point is that on a narrow sector there are generally many
holomorphic {isomorphisms} between $A^0$ and $A$ which are 
asymptotic to $\wh F$ and one is being chosen in a
canonical way; $\Sigma_i(\wh F)$ is in
fact uniquely characterised by property $2)$. There are various ways to
construct $\Sigma_i(\wh F)$, although the details will not be needed here.
In particular the series $\wh F$ is `$(k-1)$-summable' on
$\text{Sect}_i$, with sum $\Sigma_i(\wh F)$---see \cite{BBRS91,MalR92,MR91}.
Other approaches appear in \cite{BJL79,L-R94}.
See also \cite{Mal95,Was76} regarding asymptotic expansions on sectors.

Functions on the quotient $\cH(A^0)$ are now obtained as follows.
Let $(A,g_0)\in\wh{\Syst}(A^0)$
be a compatibly framed connection germ
and let $\wh F\in G\flb z \frb$ be the associated  formal isomorphism
from Lemma \ref{lem: get F}.
The sums of $\wh F$ on the two sectors adjacent to some anti-Stokes
ray $d_i\in\IA$ may be analytically continued across $d_i$ and they
will generally be different on the overlap. 
Thus for each anti-Stokes ray $d_i$ there is a 
matrix of holomorphic functions 
$\kappa_{i}:= \Sif^{-1}\circ\Sigma_{i-1}(\wh F)$
asymptotic to $1$ on a
sectorial neighbourhood of $d_i$.
Moreover clearly $\kappa_i[A^0]=A^0$; it is an automorphism of $A^0$.
A concrete description of $\kappa_i$ is obtained by choosing a basis
of solutions of $A^0$, which is made via a choice of branch of
$\log(z)$.

Thus choose a branch of $\log(z)$ along $d_1$ and extend it in a
positive sense across $\Sect_1,d_2,\Sect_2,d_3,\ldots,\Sect_r=\Sect_0$ 
in turn. 
In particular we get a lift $\wt p$ of the point $p\in \Sect_0$ to the
universal cover of the punctured disc $\ID\setminus\{0\}$ and we
will say that these $\log(z)$ choices are {\em associated to $\wt p$}.
\begin{defn}	\label{defn: sfs1}
Fix data $(A^0,z,\wt p)$ as above. 
The {\em Stokes factors} of a compatibly framed connection
$(A,g_0)\in\wh\Syst(A^0)$ are:
$$K_i:=e^{-Q}z^{-\Lo} \cdot\kappa_{i}\cdot z^\Lo e^Q,\qquad 
i=1,\ldots,r=\#\IA$$
using the choice of $\log(z)$ along $d_i$, where 
$\kappa_{i}:= \Sif^{-1}\circ\Sigma_{i-1}(\wh F)$.
\end{defn}
Since $z^\Lo e^Q$ is a fundamental solution of $A^0$ (i.e. its columns
are a basis of solutions) we have $d(K_i)=0$; the Stokes factors are
constant invertible matrices.
By part $3)$ of Theorem \ref{thm: multisums}, $K_i$ only depends on
the $G\{z\}$ orbit of $(A,g_0)$ and so matrix entries of $K_i$
are functions on $\cH(A^0)$.
A useful equivalent definition is:
\begin{defn}	\label{def: sf2}
Fix data $(A^0,z,\wt p)$ and choose $(A,g_0)\in\wh\Syst(A^0)$.

$\bullet$
The {\em canonical fundamental solution} of $A$ on the $i$th sector is
$\Phi_i:=\Sif z^\Lo e^Q$
where $z^\Lo$ uses the choice (determined by $\wt p$)
of $\log(z)$ on $\Sect_i$. (Note $\Phi_{i+r}=\Phi_i$.)

$\bullet$
If $\Phi_i$ is continued across the
anti-Stokes ray $d_{i+1}$ then on $\Sect_{i+1}$ we have:
$ K_{i+1} := \Phi_{i+1}^{-1}\circ\Phi_i
\text{ for all $i$ except }
 K_{1} := \Phi_{i+1}^{-1}\circ\Phi_i \circ M^{-1}_0
\text{ for $i=r$},  $
where $M_0:=e^{2\pi\sqrt{-1}\cdot\Lo}$ is the so-called `formal monodromy'.
\end{defn}
Taking care to use the right $\log(z)$ choices it is straightforward to
prove the equivalence of these two definitions of the Stokes factors.
The basic fact then is:
\begin{lem}	\label{lem: sf in sfg}
The Stokes factor $K_i$ is in the group $\ISto_{d_i}(A^0)$.
\end{lem}
\pf
From Theorem \ref{thm: multisums}, 
$\Sigma_{j}(\wh F)$
is asymptotic to 
$\wh F$ at $0$  when continued within the supersector
$\Ssect_j$, for each $j$.
Thus (if $i\ne 1$)
$z^\Lo e^Q K_i e^{-Q}z^{-\Lo} = 
 \Sigma_{i}(\wh F)^{-1}\Sigma_{i-1}(\wh F)$
is asymptotic to $1$ within the intersection
$\Ssect_i\cap\Ssect_{i-1}$.
As $K_i$ is constant we must therefore have
$(K_i)_{ab}=\delta_{ab}$ unless 
$e^{q_a-q_b}\to 0$ as $z\to\ 0$ along any ray in 
$\Ssect_i\cap\Ssect_{i-1}$.
It is straightforward to check this is equivalent to 
$(ab)$ being a root of $d_i$. (The $i=1$ case is similar.)
\epf
Thus as in Lemma \ref{lem: sf to sm} we can define the Stokes 
matrices of $(A,g_0)\in \wh \Syst(A^0)$:
$$S_i:=P^{-1}K_{il}\cdots K_{(i-1)l+1}P\in U_{+/-}$$
if $i$ is odd/even, where $i=1,\ldots,2k-2$ and $P$ is the permutation
matrix associated to the half-period $(d_1,\ldots,d_l)$.
To go directly from the canonical solutions to the Stokes matrices,
simply observe that
if $\Phi_{il}$ is continued in a positive sense across
all the anti-Stokes rays $d_{il+1},\ldots,d_{(i+1)l}$ and 
onto $\Sect_{(i+1)l}$ we have:
$\Phi_{il}=\Phi_{(i+1)l}PS_{i+1}P^{-1}$
for $i=1,\ldots,2k-3$, and
$\Phi_{il}=\Phi_lPS_1P^{-1}M_0$ 
for $i=2k-2 = r/l$
where $M_0=e^{2\pi\sqrt{-1}\Lo}$.
The main fact we need is then:

\begin{thm}[Balser, Jurkat, Lutz \cite{BJL79}] \label{thm: bjl}
Fix the data $(A^0,z,\wt p)$ as above.
Then the `local monodromy map' taking the Stokes matrices induces a
{\em bijection}
$$\cH(A^0)\mapright{\cong} (U_+\times U_-)^{k-1}; \qquad 
[(A,g_0)]\longmapsto (S_1,\ldots,S_{2k-2}).$$
In particular $\cH(A^0)$ is isomorphic to the vector space 
$\IC^{(k-1)n(n-1)}$.
\end{thm}
\noindent {\bf Sketch.\  }
For injectivity, suppose two compatibly framed systems in 
$\wh \Syst(A^0)$ have the same Stokes matrices.
Let $\wh F_1, \wh F_2$ be their associated formal isomorphisms 
(from Lemma \ref{lem: get F}).
Since the Stokes matrices (and therefore the Stokes factors and the
automorphisms $\kappa_j$) are equal, the holomorphic matrix 
$\Sigma_i(\wh F_2)\circ\Sigma_i(\wh F_1)^{-1}$ has no monodromy
around $0$ and does not depend on $i$.
Thus on any sector it has asymptotic expansion $\wh F_2\circ \wh F_1^{-1}$
and so (by Riemann's removable singularity theorem) we deduce 
the power series $\wh F_2\circ \wh F_1^{-1}$ is actually convergent
with the function $\Sigma_i(\wh F_2)\circ\Sigma_i(\wh F_1)^{-1}$
as sum. This gives an isomorphism between the systems we began
with: they represent the same point in $\cH(A^0)$.
Surjectivity follows from a result of Sibuya: See \cite{BJL79} Section 6.\hfill$\square$\\
\begin{rmk}
The set $\cH(A^0)$ is also described (by the Malgrange-Sibuya
isomorphism) as the first cohomology of a sheaf of non-Abelian
unipotent groups over the circle $S^1$, explaining our notation.
However we will not use this viewpoint: the sums  
$\Sif$ lead to canonical
choices of representatives of the cohomology classes that occur.
See \cite{BV89, L-R94, MR91} and the survey \cite{Var96}.
\end{rmk}
Finally two (by now easy) facts that we will need are:
\begin{cor} \label{cor: t and m}

$\bullet$
The torus action on $\cH(A^0)$ changing the compatible framing
corresponds to the conjugation action $t({\bf S})=
(tS_1t^{-1},\ldots,tS_{2k-2}t^{-1})$ on the Stokes matrices,
and so there is a bijection $\Syst(A^0)/G\{z\}\cong (U_+\times U_-)^{k-1}/T$ 
between the set of isomorphism
classes of germs of meromorphic connections formally equivalent to
$A^0$ and the set of $T$-orbits of Stokes matrices. 

$\bullet$
If $\Phi_0$ is continued once around $0$ in a positive sense, then on
return to $\Sect_0$ it will become 
$$\Phi_0\cdot PS_{2k-2}\cdots S_2 S_1 P^{-1}M_0$$
where $M_0=e^{2\pi\sqrt{-1}\cdot\Lo}$ is the formal monodromy.
\end{cor}
\pf
The first part is immediate from Theorem \ref{thm: multisums} statement
$3)$. 
For the second part, from Definition \ref{def: sf2} we see $\Phi_0$
becomes $\Phi_i\cdot K_i\cdots K_2K_1M_0$ when continued to 
$\Sect_i$. Then observe $K_r\cdots K_1=PS_{2k-2}\cdots S_1P^{-1}$.
\epf
\subsection*{Global Monodromy}

Recall we have fixed the data $\ba$ of 
a divisor $D=\sum k_i(a_i)$ on $\IP^1$ and connection germs
$d-\iAo$ at each $a_i$.
Now also choose $m$ disjoint open discs $D_i$ on $\IP^1$ with
$a_i\in D_i$ and, for each $i$, a coordinate $z_i$ on $D_i$ vanishing
at $a_i$.
Thus the local picture above is repeated on each such disc.
Abstractly the monodromy manifolds will be defined as spaces of
representations of the following groupoid $\wt \Ga$.

Choose a base-point $p_0\in\IP^1\setminus\{a_1,\ldots,a_m\}$
and
a point 
$b_\xi$ in each of the sectors bounded by anti-Stokes
directions at each pole $a_i$,
where $\xi$ ranges over some finite set indexing these sectors.
Let $\wt B_i$ denote the (discrete) subset of 
points of the universal cover of the punctured disc $D_i\setminus \{ a_i\}$, 
which are above one of the $b_\xi$'s.
Let $\wt B := \{ p_0 \}\cup\wt B_1\cup\cdots\cup\wt B_m$. 
If $\wt p\in\wt B$ we will write $p$ for the point of $\IP^1$
underlying $\wt p$ (namely $p_0$ or one of the $b_\xi$'s).  

\begin{defn}

$1)$ The set of objects of the groupoid $\wt \Ga$ is the set $\wt B$.

$2)$ If $\wt p_1,\wt p_2\in \wt B$, 
the set of morphisms of $\wt \Ga$ from $\wt p_1$ to $\wt p_2$ is the
set of homotopy classes of paths
$\ga:[0,1]\to\IP^1\setminus\{a_1,\ldots,a_m\}$ 
from  $p_1$ to $p_2$.
\end{defn}
This is clearly a groupoid with multiplication (of composable
morphisms) defined by path composition.

Now let $(V,\nabla,{\bf g})$ be a compatibly framed meromorphic
connection with irregular type $\iAo$ at $a_i$ for each $i$.
(Thus, if $V$ is trivial, $(V,\nabla,{\bf g})$  represents a point
of the extended moduli space $\wt \M^*(\ba)$.)
For each choice of basis of the fibre $V_{p_0}$ of $V$ at $p_0$ such
$(V,\nabla,{\bf g})$ naturally determines a representation of 
the groupoid $\wt \Ga$ in the group $G=GL_n(\IC)$, as follows.

Suppose $[\ga_{\wt p_2\wt p_1}]$ is a morphism in $\wt \Ga$, 
represented by a path $\ga_{\wt p_2\wt p_1}$ in the punctured sphere
from $p_1$ to $p_2$.
Then from Definition \ref{def: sf2} (with $\wt p=\wt p_i$)
we obtain a canonical choice of basis
$\Phi_i:\IC^n\to V$ of $\nabla$-horizontal sections 
of $V$ in a neighbourhood
of $p_i$ for $i=1,2$.
(First use any local trivialisation of $V$, and then
observe the basis obtained is independent of this choice. In the case
$\wt p_i=p_0$, use the choice of basis of $V_{p_0}$ to determine $\Phi_i$.) 
Now both bases extend uniquely (as solutions of $\nabla$) 
along the track 
$\ga_{\wt p_2\wt p_1}([0,1])$ of the path $\ga_{\wt p_2\wt p_1}$.
Since they are both $\nabla$-horizontal bases we have 
$\Phi_1=\Phi_2\cdot C$
on the track of $\ga_{\wt p_2\wt p_1}$, for some constant invertible matrix 
$C\in G$. 
The representation $\rho$ of $\wt \Ga$ is defined by setting
\begin{equation}	\label{connection matrices}
\rho(\ga_{\wt p_2\wt p_1}) := C = {\Phi_2}^{\!\!-1}\Phi_1.
\end{equation}
Clearly $C$ only depends on the homotopy class of the path in
$\IP^1\setminus\{a_1,\ldots,a_m\}$ and it is easy to check this is indeed
a representation. (For example $\rho$ maps contractible loops to $1$ 
and has composition property 
$\rho(\ga_{\wt p_3\wt p_2}\cdot\ga_{\wt p_2\wt p_1}) = 
\rho(\ga_{\wt p_3\wt p_2})\cdot \rho(\ga_{\wt p_2\wt p_1}).$)

Thus $\rho$ encodes all possible `connection matrices' between sectors at
different poles as well as all the Stokes factors and Stokes matrices at each
pole. 
To characterise the representations of $\wt \Ga$ that arise in this way we
observe: 

\begin{lem} \label{lem: stokes conds}
The representation $\rho$ has the following two properties:

$(${\em SR}$1)$
For any $i$, if $\wt p_1\in\wt B_i$ and $\wt p_2$ is the next element
of $\wt B_i$ after $\wt p_1$ 
in a positive sense and $\ga_{\wt p_2\wt p_1}$ is a
small arc in $D_i$ from $p_1$ to $p_2$  then 
$\rho(\ga_{\wt p_2\wt p_1})\in\ISto_d(\iAo)$, 
where $d$ is the unique anti-Stokes ray that $\ga_{\wt p_2\wt p_1}$ crosses.

$(${\em SR}$2)$
For each $i$ there is a diagonal matrix $\iL$
(which has distinct eigenvalues mod $\IZ$ if $k_i=1$) such that
for any $\wt p_1\in \wt B_i$, $\wt p_2\in \wt B$ and morphism
$\ga_{\wt p_2\wt p_1}$:
$$\rho(\ga_{\wt p_2(\wt p_1+2\pi)}) = 
\rho(\ga_{\wt p_2\wt p_1})\cdot\exp({2\pi\sqrt{-1}\cdot\iL})$$
where $\ga_{\wt p_2(\wt p_1+2\pi)}=\ga_{\wt p_2\wt p_1}$ as paths, but 
$(\wt p_1+2\pi)$ is the next point of $\wt B_i$ after $\wt p_1$ 
(in a positive sense) which is also above $p_1$.
\end{lem}
\pf
The first part  is immediate from Definition \ref{defn: sfs1} and
Lemma \ref{lem: sf in sfg} 
whilst the second is clear from the definition of the
canonical solutions, with $\iL$ the exponent of formal monodromy 
of $(V,\nabla,{\bf g})$ at $a_i$.
\epf
\begin{defn}	\label{dfn: sto rep}

$\bullet$
A {\em Stokes representation} $\rho$ is a representation of the
groupoid $\wt \Ga$ into $G$ 
together with a choice of $m$ diagonal matrices $\iL$ such that
(SR1) and (SR2) hold.
The set of Stokes representations will be denoted $\Hom_\IS(\wt\Ga,G)$. 

$\bullet$
The matrices $\iL$ associated to a Stokes representation $\rho$
will be called the {\em exponents of formal monodromy} of $\rho$
and the 
number $\deg(\rho):= \sum_i\tr(\iL)$ is the {\em degree} of $\rho$.

$\bullet$
Two Stokes representations are {\em isomorphic} if they are in the same orbit
of the following $G$ action on $\Hom_\IS(\wt \Ga,G)$:
if $\wt p_1,\wt p_2\in\wt B\setminus \{p_0\}$, $g\in G$ define
$$
\begin{array}{lcl}
(g\cdot\rho)(\ga_{p_0 p_0}) = g \rho(\ga_{p_0 p_0})g^{-1}, &\ &
(g\cdot\rho)(\ga_{p_0 \wt p_1}) = g \rho(\ga_{p_0 \wt p_1}), \\
(g\cdot\rho)(\ga_{\wt p_2 p_0}) = \rho(\ga_{\wt p_2 p_0})g^{-1}, &\ &
(g\cdot\rho)(\ga_{\wt p_2\wt p_1}) = \rho(\ga_{\wt p_2 \wt p_1}).
\end{array}
$$

$\bullet$
The {\em extended monodromy manifold} 
$\wt M(\ba):= \Hom_{\IS}(\wt \Ga, G)/G$ is the set of isomorphism classes of
Stokes representations.
\end{defn}
Observe that 
this $G$ action on $\Hom_{\IS}(\wt \Ga, G)$ corresponds to the choice of basis
of the fibre $V_{p_0}$ made above, and 
so a compatibly framed meromorphic connection $(V,\nabla,{\bf g})$ canonically
determines a point of $\wt M(\ba)$. 
(Also $\wt M(\ba)$ does not depend on the choices of the base-points
$p_0,b_\xi$ that were used to define the groupoid $\wt \Ga$.)
\begin{prop}[see \cite{JMU81}]	\label{prop: im iff im}
Two compatibly framed meromorphic connections are isomorphic if and only if 
they have
isomorphic Stokes representations. 
\end{prop}
\pf
Suppose $(V^1,\nabla^1,{\bf g^1}), (V^2,\nabla^2,{\bf g^2})$ are
isomorphic and both have irregular type $\iAo$ at each $a_i$.
Thus there is a vector bundle isomorphism $\varphi:V^1\to V^2$ which
relates the connections and the framings.
It is easy to check now that, for each $i$,
$\varphi$ also relates the canonical
bases ${\Phi^{1,2}}(z_i):\IC^n\to V^{1,2}$
of solutions on each sector at $a_i$, associated to any point 
$\wt p\in\wt B_i$. 
This implies the Stokes representations are isomorphic.
Conversely if the Stokes representations are isomorphic
the local isomorphisms
${\Phi^2}\circ ({\Phi^1})^{-1}:V^1\to V^2$ extend to $\IP^1$
to give the desired isomorphism $\varphi$, as in the proof of
Theorem \ref{thm: bjl}. 
\epf
Thus on restricting attention to connections on trivial bundles we get a
well-defined injective map
$\wt\nu:\wt\M^*(\ba)\to \wt M(\ba)$
from the extended moduli space $\wt\M^*(\ba)$ of Section \ref{sn: triv}.
This is the (extended) monodromy map and is the key ingredient in the
whole isomonodromy story.
It is a map between complex manifolds of the same dimension (see 
Lemma \ref{lem: sdim} and Proposition \ref{prop: alg im}) and 
moreover results of Sibuya and Hsieh \cite{Sib68,HsSib66,Sib62}
imply it is {\em holomorphic}.
(They prove each canonical
fundamental solution varies holomorphically with parameters 
and therefore so does all the
monodromy data---see also \cite{JMU81} Proposition 3.2.)
It follows immediately that $\wt\nu$ is surjective on tangent vectors and
biholomorphic onto its image (since any injective holomorphic map between 
equi-dimensional complex manifolds has these properties---see for example 
\cite{Ran86} Theorem 2.14).
We will see in Section \ref{sn: imds} that the image of $\wt \nu$ is
the complement of a divisor in the degree zero component of $\wt M(\ba)$.

Now we wish to describe the monodromy manifold 
$\wt M(\ba)$ more explicitly and this requires the
following choices:

\begin{defn}	\label{def: tentacles}
A choice of {\em tentacles} $\T$ is a choice of:

$1)$ A point $p_i$ in some sector at $a_i$ between two anti-Stokes rays
($i=1,\ldots,m$). 

$2)$ A lift $\wt p_i$ of each $p_i$ to the universal cover of the
punctured disc $D_i\setminus\{a_i\}$.

$3)$ A base-point  $p_0\in\IP^1\setminus\{a_1,\ldots,a_m\}$.

$4)$ A path $\ga_i:[0,1]\to \IP^1\setminus\{a_1,\ldots,a_m\}$ in the
punctured sphere, from $p_0$ to $p_i$ for $i=1,\ldots,m$, such that the loop
\begin{equation}	\label{contractible loop}
(\ga_m^{-1}    \cdot \be_m	\cdot\ga_m) 		\cdots\cdots
(\ga_2^{-1}    \cdot \be_2	\cdot\ga_2) 		\cdot 
(\ga_1^{-1}    \cdot \be_1	\cdot\ga_1)
\end{equation}
based at $p_0$ is contractible in $\IP^1\setminus\{a_1,\ldots,a_m\}$,
where $\be_i$ is any loop in $D_i\setminus\{a_i\}$ based at $p_i$ encircling
$a_i$ once in a positive sense.
\end{defn}
\begin{prop}	\label{prop: alg im}
For each choice of tentacles $\T$ there is an explicit algebraic 
isomorphism 
$\wt\varphi_\T: \wt M(\ba)\to \wcc_1\times\cdots\times\wcc_m\spq G$
from the extended monodromy manifold to the `explicit monodromy
manifold' of Definition \ref{def: gcc} and (\ref{eqn: msq}). 
\end{prop}
\pf
The choice $\T$ determines an isomorphism 
$\Hom_\IS(\wt \Ga,G)\mapright{\cong} {\bf{\rho}}^{-1}(1)\subset 
\wcc_1\times\cdots\times\wcc_m$ as follows.
Recall, using the convention used before, that the chosen point 
$\wt p_i\in \T$
determines a labelling of, and a $\log(z_i)$ choice on, each sector
and anti-Stokes ray at $a_i$.
Let $\iwbj$ be the element of $\wt B_i$ lying in the corresponding
lift of the $j$th sector $^i\Sect_j$ at $a_i$ to the universal cover
of the punctured disc $D_i\setminus\{a_i\}$.
Without loss of generality we assume that 
$\iwbo = \wt p_i$ and that the base-point $p_0$ of $\wt \Ga$ and $\T$
is the same.
Also the labelling determines a permutation matrix $P_i$ associated
to each $a_i$
(see Lemma \ref{lem: sf to sm}). (If $k_i=1$ set $P_i=1$.)
Let $\ga_{\wt p_ip_0}$ be the morphism of $\wt \Ga$ from $p_0$ to
$\wt p_i$ corresponding to the path $\ga_i$ and define
$C_i:=P_i^{-1}\rho(\ga_{\wt p_ip_0})\in G$ for $i=1,\ldots,m$.
Next let $\isj$ be the morphism from 
$^i\wt b_{(j-1)\cdot l}$ to $^i\wt b_{j\cdot l}$ with underlying path
a simple arc in $D_i\setminus\{a_i\}$ from 
$^ib_{(j-1)\cdot l}$ to $^ib_{j\cdot l}$ in a positive sense (where
$l=l_i=r_i/(2k_i-2)$ and  $r_i=\#{^i\!\IA}$). 
Then define the Stokes matrices (as explained before 
Theorem \ref{thm: bjl}) by the formulae:
${^i\!S_j}:= P^{-1}_i\rho(\isj)P_i$ for 
$j=2,\ldots,2k_i-2$ and 
${^i\!S_1}:= P^{-1}_i\rho(\iso)\cdot{^i\!M_0^{-1}}\cdot P_i$.
Finally set $\iLp=P_i^{-1}\iL P_i$ 
where $\iL$ is the $i$th exponent of formal monodromy of $\rho$
(Definition \ref{dfn: sto rep}).
Thus a Stokes representation $\rho$ determines a point 
$({\bf C},\bs,{\bf \Lambda'})$
of the product $\wcc_1\times\cdots\times\wcc_m$,
where ${\bf C}=(C_1, C_2,\ldots,C_m)$, $\bs=({^1\bs,\ldots,{^m\bs}})$, 
$^i\bs:=({^i\!S_1},\ldots,{^i\!S_{2k_i-2}})$ and 
${\bf \Lambda'} = ({^1\!\Lambda'},\ldots,{^m\!\Lambda'})$.
Now observe that the value 
$\rho(\ga_i^{-1}\cdot \be_i\cdot\ga_i)$ of the representation $\rho$
on the loop $\ga_i^{-1}\cdot \be_i\cdot\ga_i$ based at $p_0$ is equal
to the value $\rho_i(C_i,{^i\bs},\iLp)$ of the map $\rho_i:\wcc\to G$,
and so the contractibility of the loop (\ref{contractible loop})
implies the monodromy data 
$({\bf C},\bs,{\bf \Lambda'})$ satisfies the constraint
$\rho_m\cdots\rho_1 = 1$. 
This defines the map
$\Hom_\IS(\wt \Ga,G)\to {\bf{\rho}}^{-1}(1)$ and it is
straightforward to see it is an isomorphism 
(using Lemma \ref{lem: sf to sm} and knowledge of the fundamental
group of the punctured sphere).
This map is $G$-equivariant and so descends to give 
$\wt\varphi_\T$.
\epf
Now we turn to the non-extended version.
First, taking the exponents of formal monodromy 
${\bf \Lambda} = ({^1\!\Lambda},\ldots,{^m\!\Lambda})$
of any Stokes representation $\rho$ induces a map
$$\mu_{T^m}:\wt M(\ba)\longrightarrow \lt^m;\qquad
\rho\longmapsto {\bf \Lambda}.$$
Also for each pole $a_i$ there is a torus action on 
$\Hom_{\IS}(\wt \Ga,G)$ defined by the formulae
\begin{equation}	\label{eq: t action on monod}
\begin{array}{lcl}
(t\cdot\rho)(\ga_{\wt p_2 \wt p_1}) = 
t\rho(\ga_{\wt p_2 \wt p_1})t^{-1} &\ &
(t\cdot\rho)(\ga_{\wt q_2 \wt p_1}) = 
\rho(\ga_{\wt q_2 \wt p_1})t^{-1} \\
(t\cdot\rho)(\ga_{\wt q_2 \wt q_1}) 
= \rho(\ga_{\wt q_2 \wt q_1}) &\ &
(t\cdot\rho)(\ga_{\wt p_2 \wt q_1}) = t\rho(\ga_{\wt p_2 \wt q_1})
\end{array}
\end{equation}
for any $\wt p_1,\wt p_2\in\wt B_i$ and 
$\wt q_1,\wt q_2\in\wt B\setminus \wt B_i$, where $t\in T$.
\begin{defn}

$\bullet$
The (non-extended) space of monodromy data $M(\ba)$ 
is the set of $T^m$ orbits in $\wt M(\ba)$ which have exponents of formal
monodromy equal to 
${\bf \Lambda\!^0}= ({^1\!\Lambda\!^0},\ldots,{^m\!\Lambda\!^0})$, where 
$\iLo= \res_{a_i}(\iAo)$:
$$M(\ba):= \mu_{T^m}^{-1}({\bf \Lambda\!^0})/T^m.$$

$\bullet$
The {\em monodromy map} is the map $\nu:\M^*(\ba)\to  M(\ba)$ induced from the
extended monodromy map.
\end{defn}
The monodromy map is well-defined since part $3)$ of 
Theorem \ref{thm: multisums}
implies that the extended monodromy
map  is $T^m$-equivariant and also since it is clear that
$\mu_{T^m}\circ \wt \nu$ is the moment map for the $T^m$ 
action on the extended moduli space $\wt \M^*(\ba)$ 
(defined in 
Proposition \ref{prop: main pp}).

\begin{cor}
For each choice of tentacles $\T$ there is an explicit algebraic 
isomorphism 
$\varphi_\T:  M(\ba)\to \cC_1\times\cdots\times\cC_m\spq G$
from the monodromy space to the explicit set of monodromy
data from Definition \ref{def: gcc} (with $\cC_i$ depending on $\T$).
\end{cor}
\pf
The choice of tentacles determines a permutation matrix $P_i$ for
$i=1,\ldots,m$.
Then define $\iLp:=P_i^{-1}\cdot\iLo\cdot P_i$ and use this value to define 
$\cC_i$. 
The rest now follows from Proposition \ref{prop: alg im} since
$\wt\varphi_\T$ is $T^m$-equivariant, where 
$(t_1,\ldots,t_m)$ acts on $\wcc_i$ via 
$(P_i^{-1}t_iP_i)\in T$ and the $T$-action of
Definition \ref{def: gcc}.
\epf
\begin{rmk}(Degree.)	\label{rmk: degree}
If $\rho$ is a Stokes representation having exponents of formal monodromy 
${\bf \Lambda}$ then the degree $\deg(\rho)= \sum_i\tr(\iL)$ of $\rho$ is an
integer.
One way to see this is to choose some tentacles so $\rho$ determines 
(via Proposition \ref{prop: alg im}) 
a point 
$({\bf C},\bs,{\bf \Lambda'})$ of $\wcc_1\times\cdots\times\wcc_m$ satisfying
the constraint $\rho_m\cdots\rho_1=1$. 
By taking the determinant of this constraint we see that 
$\sum_i\tr(\iLp)=\sum_i\tr(\iL)\in\IZ$.
(It is also clear that 
the fixed-degree components $\wt M_d(\ba)$ of the extended
monodromy manifolds are pairwise isomorphic.)
On the other hand suppose $(V,\nabla,{\bf g})$ is a compatibly framed
meromorphic connection on a holomorphic vector bundle $V\to\IP^1$
with irregular type $\ba$ and exponents of formal monodromy ${\bf \Lambda}$.
Then by considering the induced connection on the determinant line bundle 
$\Lambda^n V$ of $V$ one finds  that $\sum_i\tr(\iL)$ is equal to the degree 
of the vector bundle $V$.  
The only point we need to make here is that the germs
$\iAo$ must be chosen such that $\sum_i\tr(\iLo)=0$,
if the moduli spaces $\M^*(\ba)$ are to be non-empty, 
and so we will tacitly assume this throughout.
\end{rmk}
To end this section we describe the dependence on the local coordinate choices
$z_i$ that were made right at the start.
Let $\ba$ be a choice of divisor $D$ and connection germs $d-\iAo$ as above.
This determines all the spaces $\wt\M^*(\ba), \wt M(\ba), \M^*(\ba)$
and $M(\ba)$.

\begin{prop}	\label{prop: coord depce}

$1)$ The extended monodromy map $\wt \nu:\wt\M^*(\ba)\to\wt M(\ba)$ depends
(only) on the choice of a $k_i$-jet of a coordinate $z_i$ at each $a_i$.

$2)$ This coordinate dependence is only within the $T^m$ orbits:
The monodromy map $\nu:\M^*(\ba)\to M(\ba)$ is completely intrinsic.
\end{prop}
\pf
The key point is to see how a fundamental solution 
$\Phi=\Sjf z^\Lambda e^Q$ changes when  the coordinate $z$ is changed.
Here $A^0=dQ+\Lambda dz/z$ is fixed and $Q$ is determined by $(A^0,z)$ by
requiring it to have zero constant term in its Laurent expansion with respect
to $z$.
Suppose $z'=ze^f$ is a new coordinate choice, 
for some local holomorphic function $f$. 
One finds $Q'=Q-\Lambda f+\Lambda f(0)-\res_0(Qdf)$ 
(as meromorphic functions near $z=0$), since then 
$\res_0(Q'dz'/z')=0$.
In turn $\Phi'=\Phi\cdot t^{-1}$ where $t=\exp(\res_0(Qdf)-\Lambda f(0))\in T$.
(The function $\Sjf$ is intrinsic.)
Then observe: $1)$ If $f=O(z^k)$ then $t=1$, since $Q$ has a pole of order
$k-1$, and $2)$ This action of $t\in T$ corresponds to the torus action we
have defined. 
\epf
\begin{rmk}	\label{rmk: pp dept}
One should also note that all the spaces 
$\wt\M^*(\ba), \wt M(\ba), \M^*(\ba)$ and 
$M(\ba)$ only depend on the principal
part of each germ $\iAo$.
For the monodromy manifolds this is immediate and for the moduli
spaces $\wt\M^*(\ba), \M^*(\ba)$ it is because all $\iAo$ with the
same principal part are formally equivalent via a transformation with
constant term $1$ (as explained in Section \ref{sn: triv}).
\end{rmk}

%% file: smid-smooth.tex
\section{${C^\infty}$  Approach to Meromorphic Connections} \label{sn: smooth}

This section gives a third viewpoint on meromorphic connections: a
$C^\infty$ approach.
Although we work exclusively with `generic' connections over $\IP^1$
(as we wish to study isomonodromic deformations of such connections)
we remark that this
$C^\infty$ approach works over arbitrary compact Riemann
surfaces (maybe with boundary)
and the generic hypothesis is also superfluous
(see Remark \ref{rmk: non-generic}).
\subsection*{Singular Connections: ${C^\infty}$ Connections with Poles}

Let $D = k_1(a_1)+\cdots+k_m(a_m)$ be an effective divisor on $\IP^1$
as usual and choose $m$ disjoint discs $D_i\subset \IP^1$
with $a_i\in\ D_i$ and a coordinate $z_i$ on $D_i$
vanishing at $a_i$.
Define the sheaf of `smooth functions with poles on $D$' to be
the sheaf of $C^\infty$ sections of the holomorphic line bundle
associated to the divisor $D$:	\label{sing forms}
$$C^\infty[D]:=\cO[D]\otimes_\cO C^\infty$$   
where $\cO$ is the sheaf of holomorphic functions
and $C^\infty$ the infinitely differentiable complex functions. 
Any local section of
$C^\infty[D]$ near $a_i$ is of the form 
$f/z_i^{k_i}$
for a $C^\infty$ function $f$.
Similarly define sheaves $\Omega^{r}[D]$ of $C^\infty$
$r$-forms with poles on $D$ 
(so in particular $\Omega^{0}[D]=C^\infty[D]$).
A basic feature is that `$C^\infty$-Laurent expansions' can be taken at each
$a_i$. 
This gives a map
\begin{equation}	\label{laurent}
L_i:\Omega^*[D](\IP^1)\to\IC\flb z_i,\bar z_i\frb z_i^{-k_i}
\otimes {\bigwedge}^{\!\!*}\IC^2
\end{equation}
where $\IC^2=\IC dz_i\oplus\IC d\bar z_i$.
For example if $f$ is a $C^\infty$ function defined in a neighbourhood
of $a_i$ then 
$L_i\left(f/z_i^{k_i}\right)=L_i(f)/z_i^{k_i}$
where $L_i(f)$ is
the Taylor expansion of $f$ at $a_i$.

The Laurent map $L_i$ has nice morphism properties,
for example
$L_i(\omega_1\wedge\omega_2)=L_i(\omega_1)\wedge L_i(\omega_2)$
and $L_i$ commutes with the exterior derivative $d$, where $d$ is
defined on the right-hand side of (\ref{laurent}) in the obvious way 
($d(z_i^{-1})=-dz_i/z_i^2$).

We will repeatedly make use of the fact that the kernel of $L_i$ consists
of nonsingular forms, that is: if $L_i(\omega)=0$ then 
$\omega$ is nonsingular at $a_i$.
This apparently innocuous statement 
is surprisingly tricky to prove directly, but since it is crucial for
us we remark it follows from the following:

\begin{lem}[Division] \ 	\label{lem: division}
Let $\ID\subset \IC$ be a disk containing the origin.
Suppose $f\in C^\infty(\ID)$ 
and that the Taylor expansion of $f$ at $0$ is in the ideal in 
$\IC\flb z,\bar z\frb$ generated by $z$.
Then $f/z\in C^\infty(\ID)$.
\end{lem}
\pf
This is a special case of a much more general result
of Malgrange \cite{Mal66}.
The particular instance here is discussed by Martinet 
\cite{Mart82} p115. \epf

Another fact we will use is that the $C^\infty$ Laurent expansion map
$L_i$ in (\ref{laurent}) is {\em surjective} for each $i$.
This is due to a classical result of E.Borel which we quote  here 
in the relative case that will be needed later:

\begin{thm} [E. Borel] \label{thm: borel1}
Suppose $U$ is a differentiable manifold, 
$I$ is a compact neighbourhood of the origin in $\IR$ and
$\widehat f\in \IC\flb x,y \frb\otimes C^\infty(U)$
(where $x,y$ are real coordinates on $\IC\cong\IR^2$).
Then
there exists a smooth function
$f\in C^\infty(U\times I\times I)$
such that
the Taylor expansion of $f$ at $x=y=0$ is $\widehat f$.
\end{thm}
\pf
This is easily deduced, via partitions of unity,
by using two applications of the version of
Borel's theorem proved on p16 of H\"ormander's book \cite{Horm83}.
\epf

Now let $V\to\IP^1$ be a rank $n$, $C^\infty$ complex vector bundle.

\begin{defn}	\label{dfn: sing conn}
A {\em $C^\infty$ singular connection
$\nabla$ on $V$ with poles on $D$} is a map
$\nabla: V\longrightarrow V\otimes \Omega^1[D]$
from the sheaf of ($C^\infty$) 
sections of $V$ to the sheaf of sections of 
$V\otimes \Omega^1[D]$,
satisfying the Leibniz rule:
$\nabla(fv)=(df)\otimes v+ f\nabla v$
where $v$ is a local section of $V$ and $f$ is a local $C^\infty$
function.
\end{defn}
Concretely in terms of the local coordinate $z_i$ on $\IP^1$ vanishing
at $a_i$ and a local trivialisation of $V$, 
$\nabla$ has the form:
$\nabla=d-\iA/z_i^{k_i}$
where ${\iA}$ is an $n\times n$ matrix of $C^\infty$ one-forms.
In this paper, to study the Jimbo-Miwa-Ueno isomonodromy equations, 
we need only to consider the case when $V$ is the trivial rank $n$, 
$C^\infty$ vector bundle over $\IP^1$. (Recall any degree zero 
vector bundle over $\IP^1$ is $C^\infty$ trivial.) 

\begin{defn} \ 	\label{dfn: D sing conns}

$\bullet$
Let $\A[D]$ denote the set of $C^\infty$ 
singular connections with poles on $D$ on
the trivial $C^\infty$ rank $n$ vector bundle:
$\A[D]:=\bigl\{ d-\al\ \ \bigl\vert \ 
\al\in \End_n\bigl(\Omega^{1}[D](\IP^1)\bigr)\bigl\}$
where $\Omega^{1}[D]$ is the sheaf of $C^\infty$ one-forms with poles
on $D$.

$\bullet$
The {\em gauge group} of $C^\infty$ bundle automorphisms is 
$\G:= GL_n(C^\infty(\IP^1)).$

$\bullet$ The {\em curvature} of a singular connection
$d-\al\in\A[D]$  is the matrix of singular two-forms 
$\F(\al) := (d-\al)^2 = -d\al + \al^2 \in 
\End_n\bigl(\Omega^2[2D](\IP^1)\bigr).$

$\bullet$ The {\em flat} connections are those with zero
curvature and the subset of flat singular connections
 will be denoted $\A_\fl[D]$.
\end{defn}
\begin{rmk}
Occasionally one comes across notions of curvature of singular 
connections involving distributional derivatives. 
For example a meromorphic connection on a Riemann surface 
is sometimes said to have a $\delta$-function singularity
in its curvature at the pole, to account for the monodromy around the pole.
The definition above of curvature
does {\em not} involve distributional derivatives, and so, for us,
{\em any} meromorphic connection over a Riemann surface is flat.
\end{rmk}
The group $\G$ of bundle automorphisms clearly acts 
on the singular connections $\A[D]$ 
and explicitly this is given by the formula
$g[\al]=g \al g^{-1} + (dg)g^{-1}.$
This restricts to an action on $\A_\fl[D]$
since
$\F(g[\al]) = g(\F(\al))g^{-1}$ for $g\in \G$. 

Now choose a generic diagonal connection germ $d-\iAo$ at $a_i$ 
for each $i$ 
and let $\ba$ denote this $m$-tuple of germs and the divisor $D$,
as usual.
Since $d-\al\in\A[D]$ is on the trivial vector bundle, and 
$d-\iAo$ is a germ of a connection on the trivial bundle, we
can compare the Laurent expansion of $\al$ at $a_i$ with $\iAo$.
In particular the following definition makes sense:

\begin{defn}	\label{dfn: fixed sing conns}

$\bullet$
Let $\A(\ba)$ be the set of singular connections with fixed Laurent
expansions:
$$\A(\ba):= \bigl\{ d-\al\in \A[D] \ \bigl\vert\ 
L_i(\al) = \iAo \text{ for each $i$} \bigr\}.$$

$\bullet$
Let $\wA(\ba)$ be the following {\em extended} set of singular connections
with fixed Laurent expansions:
$$\wA(\ba):= \bigl\{ d-\al\in \A[D] \ \bigl\vert\ 
L_i(\al) = \iAo + (\iL-\iLo)\frac{dz_i}{z_i} 
\text{ for some $\iL\in\lt_i$} \bigr\}$$
where $\iLo= \res_0(\iAo)$,
$\lt_i=\lt$ if $k_i\ge 2$ and  $\lt_i=\lt'$ if $k_i=1$. 

$\bullet$
Let $\G_T$ and $\G_1$ denote the subgroups of $\G$ of 
elements
having Taylor expansion equal to a constant diagonal matrix or the
identity, respectively, at each $a_i$.
\end{defn}
The basic motivation for this definition is 
Corollary \ref{cor: germ bijection} below.
Note that $\A(\ba)$ is an affine space and that 
if $d-\al\in\A(\ba)$ then  
(from the division lemma above)  
the $(0,1)$ part of $\al$ is {\em nonsingular} over all of $\IP^1$.

\subsection*{Smooth Local Picture}\label{sn: slp}

Now we will give a $C^\infty$ description of the sets $\cH(A^0)$
and the local analytic classes $\Syst(A^0)/G\{z\}$ 
defined in Section \ref{sn: gm}.

We begin with a straightforward observation.
Let $z$ be a complex coordinate on the unit disc $\ID\subset\IC$. 
From Borel's theorem we have an exact sequence of groups:
$$1\mapright{}\og_1\mapright{}\og\mapright{L_0}
GL_n(\IC\flb z,\bar z \frb)\mapright{}1.$$
where $\og$ is the group of germs at $0$ of gauge transformations
$g\in\G$ 
and
$\og_1:=\ker(L_0)$ is the subgroup of germs with Taylor expansion $1$.

Fix a generic diagonal connection germ $d-A^0$ with an order $k$ pole
at $z=0$. 
By projecting a marked pair $(A,\wh F)$ 
onto its second factor we obtain an injection
$\wh\Syst(A^0)\hookrightarrow G\flb z\frb$ 
(see Lemma \ref{lem: get F})
and so in turn may regard $\wh\Syst(A^0)$ as a subset of 
$GL_n(\IC\flb z,\bar z \frb)$.
Define $\wh\cS (A^0):= L_0^{-1}(\wh\Syst(A^0))$ to be the lift of this
subset to $\og$.
Also lift the stabiliser torus $T\cong (\IC^*)^n$, 
that is define $\og_T :=L_0^{-1}(T).$
\begin{lem} 		\label{smooth lifts}
Taking Taylor series at $0$ induces isomorphisms:
\begin{equation}	\label{smooth im 1}
G\{ z\} \backslash \wh\cS(A^0)/ \og_1
\ \cong \ \cH(A^0) 
\end{equation}
and (by considering the residual action of $\og_T/\og_1\cong T$):
\begin{equation}	\label{smooth im 2}
G\{ z\} \backslash \wh\cS(A^0)/ {\og_T}
\ \cong \ \Syst(A^0)/G\{z\}.
\end{equation}
\end{lem}
\pf
First observe $L_0$ induces isomorphisms 
$\wh\cS(A^0)/\og_1 \cong \wh\Syst(A^0)$
and
$\og_T/\og_1 \cong T$. 
Then recall $\cH(A^0):= G\{z\}\backslash\wh\Syst(A^0)$
\epf

Having lifted things up into a smooth context 
a new interpretation of the smooth quotients above  
will be given.
In particular it is desirable to remove the groups $G\{z\}$ occurring
on the left-hand sides in  (\ref{smooth im 1}) and 
(\ref{smooth im 2}).

Let $\oa[k]=\oa[k(0)]$ \label{dfn: oak}
denote the set of germs at $0$ of $C^\infty$
singular connections
on the trivial bundle, with poles of order at most $k$.

Now given $g\in \wh\cS(A^0)$ we can apply the formal transformation
$L_0(g)$ to $A^0$ to obtain a meromorphic connection 
$A:=L_0(g)[A^0]$.
Now apply the $C^\infty$ gauge transformation $g^{-1}$ to $A$ to
define a singular connection
$\si(g):= g^{-1}[A]=g^{-1}[L_0(g)[A^0]]$.
This defines a map $\si:\wh\cS(A^0)\rightarrow\oa[k]$.
Observe that 

$\bullet$
$\si(g)$ has Laurent expansion $A^0$
(from the morphism properties of $L_0$),

$\bullet$
If $h\in G\{ z\}$ is holomorphic then  $\si(hg)=\si(g)$, 
as $h[A] = L_0(h)[A]$, and

$\bullet$ $\si(g)$ is a {\em flat} singular connection,
since it is $C^\infty$ gauge equivalent to the meromorphic
connection $A$.

Thus $\si$ gives a map into the flat connection germs with Laurent
expansion $A^0$, i.e. into $\oa_\fl(A^0)$. 
In fact it is surjective and its fibres are precisely the $G\{z\}$ orbits:

\begin{prop}		\label{big local bijection}
The map $\si$ defined above 
induces an isomorphism
$$G\{ z\}\backslash \wh\cS(A^0)\mapright{\cong} \oa_{\fl}(A^0)$$
onto the set of flat singular connection germs with Laurent
expansion $A^0$. 
\end{prop}
\pf
We have seen the induced map is well defined and now  show it is bijective.
For surjectivity, suppose $d-\al\in\oa_{\fl}(A^0)$ is a flat singular
connection with Laurent expansion $A^0$.
Thus the $d\bar z$ component
$\al^{0,1}$ of $\al$ has zero Laurent expansion at $0$ and so in
particular is nonsingular.
It follows (see \cite{AB82} p555 or \cite{BV89} p67)
that there exists $g\in \og$
with $(\bar\partial g)g^{-1} = \al^{0,1}$ and so
$A:=g^{-1}[\al]$ is still flat and has no $(0,1)$ part.
By writing $A=\gamma dz/z^k$ for smooth $\gamma$ observe that flatness 
implies $\bar\partial\gamma = 0$ and so $A$ is meromorphic.
We claim now that $A$ is formally equivalent to $A^0$, and that
$L_0(g)$ is a formal isomorphism between them.
Firstly $L_0(g)$ has no terms containing $\bar z$ because
$L_0(\bar\partial g)=L_0(\al^{0,1}g)=0$ since $L_0(\al^{0,1})=0$.
Secondly just observe
$$L_0(g^{-1})[A^0]=L_0(g^{-1})[L_0(\al)]=L_0(g^{-1}[\al])=L_0(A)=A$$
and so the claim follows.
In particular $g^{-1}\in\wh\cS(A^0)$ and by construction
$\si(g^{-1})=\al$ and so $\si$ is onto.
Finally if $g_1[A]=g_2[B]$ with $A,B$ meromorphic then
$h[A]=B$ with $h:=g_2^{-1}g_1$.
Looking at $(0,1)$ parts 
gives $(\bar\partial h)h^{-1}=0$ and so $h$ is
holomorphic. This proves injectivity.
\epf

Combining this with Lemma \ref{smooth lifts} immediately yields the
main local result:
\begin{cor}		\label{cor: germ bijection}
There are canonical isomorphisms:
$$\oa_{\fl}(A^0)/\og_1 \cong \cH(A^0)\quad\text{and}\quad
\oa_{\fl}(A^0)/\og_T \cong \Syst(A^0)/G\{z\}$$
between the $\og_1$ orbits of flat singular connection germs with Laurent
expansion $A^0$ and the set of analytic equivalence classes of
compatibly framed systems with formal type $A^0$, and
between the $\og_T$ orbits of flat singular connection germs with Laurent
expansion $A^0$ and the set of analytic equivalence classes of
connection germs formally equivalent to $A^0$.
\end{cor}
\pf
This follows directly by substituting $\oa_\fl(A^0)$ for
$G\{z\}\backslash \wh\cS(A^0)$ in Lemma \ref{smooth lifts}.
In summary: to go from a flat singular connection
$d-\al\in\oa_\fl(A^0)$ to $\cH(A^0)$ just solve 
$(\bar \partial g)g^{-1}=\al^{0,1}$ and take the $G\{z\}$ orbit of 
$L_0(g^{-1})\in G\flb z\frb$ to give an element of $\cH(A^0)$
(see the proof of Proposition \ref{big local bijection}).
Conversely, given $\wh F\in G\flb z\frb$ such that $A:=\wh F[A^0]$ is
convergent, 
use E.Borel's theorem to find $g\in\og$ such
that $L_0(g)=\wh F^{-1}$. 
Then set $\al=g[A]$ to give $\al\in\oa_\fl(A^0)$.
\epf

Thus the analytic equivalence classes may be encoded
in an entirely $C^\infty$ way.
These bijections can be thought of as relating the two 
distinguished types of elements (the meromorphic connections
and the connections with fixed Laurent expansion)
within the $^0\G$ orbits in $\oa_{\fl}[k]$.
That is, they relate the conditions $\al\in\Syst(A^0)$ and 
$\al\in\oa_{\fl}(A^0)$ on 
$\al\in\oa_{\fl}[k]$ by moving within $\al$'s $\og$ orbit.

\begin{rmk}	\label{rmk: non-generic}
Corollary \ref{cor: germ bijection}
easily extends to the general
(non-generic) case, with the same proof.
The precise statement is as follows (but won't be needed elsewhere in
this paper).
Let $d-A$ be {\em any} meromorphic connection germ and let 
$\og_{\text{\rm\footnotesize Stab}(A)}$ be the subgroup of $\og$
consisting of elements  $g$ whose Taylor expansion stabilises
$A$ (i.e. $L_0(g)[A]=A$).
Then the set of analytic isomorphism 
classes of meromorphic connection germs
formally equivalent to $A$ is canonically isomorphic to the set of
$\og_{\text{\rm\footnotesize Stab}(A)}$-orbits 
of flat singular connection germs with Laurent
expansion $A$: 
$\Syst(A)/G\{z\}\cong \oa_\fl(A)/\og_{\text{\rm\footnotesize Stab}(A)}.$
Similarly $\cH(A):=\wh \Syst_{\rm mp}(A)/G\{z\}$
is canonically isomorphic to $\oa_\fl(A)/\og_1$
(but in general this cannot be interpreted in terms of 
compatibly framed systems, only in terms of marked pairs).
\end{rmk}
\subsection*{Globalisation}      \label{globalisation}

Recall we have fixed the data $\ba$ 
(of a
divisor $D=\sum k_i(a_i)$ on $\IP^1$ and 
connection germs $d-\iAo$)
and defined $\A(\ba)$ to be the
set of singular connections on the trivial rank $n$ vector bundle on
$\IP^1$ having Laurent expansion $\iAo$ at $a_i$ for each $i$.
Following the results of the last section we are led to consider such 
connections which are {\em flat}.
The main result is:
\begin{prop}	\label{prop: global bijn 1}
There is a canonical bijection between the set of $\G_T$ orbits of
flat $C^\infty$ singular connections with fixed Laurent
expansions $\ba$ and the set of isomorphism classes of meromorphic
connections with
formal type $\ba$ on {\em degree zero} holomorphic bundles over $\IP^1$:
$$\M(\ba)\cong \A_\fl(\ba)/\G_T.$$
\end{prop}
\pf
Suppose $(V,\nabla)$ represents an isomorphism class in $\M(\ba)$.
The meromorphic connection $\nabla$ is in particular a $C^\infty$
singular connection, according to Definition \ref{dfn: sing conn}.
Since $V$ is degree zero it is $C^\infty$ trivial so, by choosing a
trivialisation, $(V,\nabla)$ determines a singular connection $d-\al$
on the trivial bundle over $\IP^1$.

From the local picture just described, 
since $\nabla$ is formally equivalent to $\iAo$
at $a_i$, we can choose $g\in\G$ such that $g[\al]$ has Laurent
expansion $\iAo$ at $a_i$ for all $i$.
This gives an element $g[\al]$
of $\A_\fl(\ba)$ and we take the $\G_T$ orbit through it to define the
required map. 
We need to check this $\G_T$ orbit only depends on the isomorphism
class of $(V,\nabla)$ and that the map is bijective.

Suppose we have two such pairs $(V,\nabla)$ and $(V',\nabla')$ and we
choose $C^\infty$ trivialisations of $V$ and $V'$ so that $\nabla,
\nabla'$ give singular connections $d-\al_1$, $d-\al_2$  
respectively.
Now a standard $\bar \partial$-operator argument implies
$(V,\nabla)\cong(V',\nabla')$ 
if and only if $\al_1$
and $\al_2$ are in the same $\G$ orbit.
Thus an isomorphism class $[(V,\nabla)]$ of
meromorphic connections determines (and is determined by) a $\G$ orbit
of singular connections on the trivial bundle.
This $\G$ orbit has a subset of singular connections having Laurent
expansion $\iAo$ at $a_i$ for each $i$.
This subset is a $\G_T$ orbit of singular connections (since $T$ is the
stabiliser of $\iAo$) and is the element of 
$\A_\fl(\ba)/\G_T$ corresponding to
$[(V,\nabla)]$. 
\epf

\begin{cor}	\label{cor: global bijn 2}
The set
$\wt \M(\ba)$ of isomorphism classes of triples
$(V,\nabla,{\bf g})$ consisting of a generic meromorphic connection
$\nabla$ (with poles on $D$) 
on a {\em degree zero} holomorphic vector bundle $V$ over $\IP^1$
with compatible framings ${\bf g}$
such that
$(V,\nabla,{\bf g})$ has irregular type $\ba$
is canonically isomorphic to the set of $\G_1$ orbits of
flat connections in $\wt\A(\ba)$:
$$\wt \M(\ba)\cong \wt\A_\fl(\ba)/\G_1.$$
\end{cor}
\pf
As in Corollary \ref{cor: germ bijection}, replacing
$\G_T$ by $\G_1$ in Proposition \ref{prop: global bijn 1} corresponds
to incorporating a compatible framing as required for $\wt \M(\ba)$.
The desired isomorphism is then obtained by simply repeating the proof
of Proposition \ref{prop: global bijn 1} for each possible set of
choices of exponents of formal monodromy 
${\bf\Lambda}$.
\epf

\subsection*{Monodromy of Flat Singular Connections}
Having related the $C^\infty$ approach to meromorphic connections we
now relate it to the monodromy approach of Section \ref{sn: gm}.
The key step is to define the generalised monodromy data of 
flat $C^\infty$ 
singular connections with fixed Laurent expansions, 
but this is easy since they too have
canonical solutions on sectors: 
\begin{lem}	\label{C^infty canonical solutions}
Suppose $\al\in\oaf(A^0)$. 
For each choice of $\log(z)$ 
there is a canonical
choice $\Phi_i$ of fundamental solution of $\al$ on $\Sect_i$, given by:
$$\Phi_i:=g\Sigma_i(L_0(g^{-1}))z^\Lo e^Q$$
for any $g\in\og$ solving $(\bar\partial g)g^{-1}=\al^{0,1}$.
\end{lem}
\pf
From the proof of Proposition \ref{big local bijection},
 such $g$
is unique up to right
multiplication by $h\in G\{z\}$ and
$A:=L_0(g^{-1})[A^0]=g^{-1}[\al]$ is a convergent meromorphic
 connection germ.
Theorem \ref{thm: multisums} then provides an analytic isomorphism  
$\Sigma_i(L_0(g^{-1}))$ between $A^0$ and $A$ on $\Sect_i$.
It follows that
$g\Sigma_i(L_0(g^{-1}))$ is an isomorphism between 
$A^0$ and $\al$ which is independent of the choice of $g$. 
Composing this with the fundamental solution $z^\Lo e^Q$ of $A^0$
gives the result.
\epf

Thus, exactly as in Section \ref{sn: gm}, a singular connection
$d-\al\in\wt\A(\ba)$ determines a Stokes representation $\rho$ 
(upto isomorphism).
This gives a map (which will be referred to as the $C^\infty$
monodromy map):
$$\wt\nu:\wt\A_{\fl}(\ba)\longrightarrow \wt M(\ba);\qquad
\al\longmapsto [\rho].$$
Since the connections in $\wt\A_{\fl}(\ba)$ are on a degree zero
bundle, the image of $\wt \nu$ is in the degree zero component 
$\wt M_0(\ba)$ of the extended monodromy manifold 
(see Remark \ref{rmk: degree}).
The main result is then:
\begin{prop}		\label{prop: gm props}
The $C^\infty$ monodromy map
$\wt\nu:\wt\A_{\fl}(\ba)\longrightarrow \wt M_0(\ba)$
is surjective and has precisely the $\G_1$ orbits in $\wt\A_{\fl}(\ba)$ as
fibres, so that
$$\wt\A_{\fl}(\ba)/\G_1 \cong \wt M_0(\ba).$$
Moreover $\wt \nu$ 
intertwines the $\G_T$ action on $\wt\A_{\fl}(\ba)$ with
the $\G_T$ action on $\wt M_0(\ba)$ defined via the evaluation map
$\G_T\to T^m$ and the
torus actions of (\ref{eq: t action on monod}).  
\end{prop}
Before proving this we deduce what the monodromy data corresponds to
in the meromorphic world:
\begin{cor}		\label{gm ims}
Taking monodromy data induces bijections:
$$\wt\M(\ba)\cong \wt M_0(\ba)\quad\text{ and }\quad \M(\ba)\cong M(\ba)$$
between the spaces of meromorphic connections on degree zero bundles
and the corresponding spaces of monodromy data.
In particular $\wt \M(\ba)$ inherits the structure of a
complex manifold from $\wt M_0(\ba)$.
\end{cor}
\pf 
The first bijection follows directly from 
Propositions \ref{prop: global bijn 1} and \ref{prop: gm props}. 
The second bijection follows from the first by fixing the exponents of
formal monodromy and quotienting by the $T^m$ action (using the
intertwining property of $\wt \nu$).
\epf

\pfms [of Proposition \ref{prop: gm props}].\ 
Choose some tentacles $\T$ and a thickening 
$\overline\ga_i:[0,1]\times [0,1]\to\IP^1\setminus\{a_1,\ldots,a_m\}$
of each path $\ga_i$ (i.e. a ribbon following $\ga_i$ whose track 
$\vert\overline\ga_i\vert$ is a closed tubular neighbourhood of the
track of $\ga_i$).
Let $D_0$ be a disc in $\IP^1$ containing $p_0$ and disjoint from each
disc $D_1,\ldots,D_m$.
Let $\vert\T\vert:= 
\overline D_0\cup\bigcup_{i=1}^m(\overline D_i\cup
\vert\overline\ga_i\vert)\subset \IP^1$
be the union of the $m+1$ discs $\overline D_i$ and the 
$m$ ribbons $\vert\overline\ga_i\vert$.
We will suppose (as is clearly possible) that 
$\T,\vert\overline\ga_i\vert, D_0$ have been chosen such that:
1) if $i\ne j$ then $\vert\overline\ga_i\vert$ only intersects 
$\vert\overline\ga_j\vert$ inside $D_0$ and 
2) that $\vert\T\vert$ is homeomorphic to a (closed) disc.

For surjectivity, let $\rho$ be any degree zero Stokes representation.
From Proposition \ref{prop: alg im}, 
specifying $\rho$ is equivalent to
specifying
the matrices $({\bf C},\bs,{\bf \Lambda'})=\wt\varphi_\T(\rho)$.
(Also as in Proposition \ref{prop: alg im},
$P_i$ will denote the permutation matrix associated to $a_i$ 
via the tentacle choice.)
Since the Stokes matrices classify germs of singular
connections up to $C^\infty$ 
gauge transformations with Taylor expansion $1$,
germs $\al_i\in{^i\!\wt\A_{\fl}(\iAo)}$ may be obtained having any given
Stokes matrices and residue for each $i=1,\ldots,m$
(combine Theorem \ref{thm: bjl} with Corollary \ref{cor: germ bijection}).
It is straightforward to extend $\al_i$ arbitrarily to $\overline D_i$. 
Next the  $\al_i$'s are patched together along the ribbons
$\vert\overline\ga_i\vert$.
Let $^i\Phi_{0}$ be the canonical solution of $\al_i$ on
$^i\Sect_{0}$ from Lemma \ref{C^infty canonical solutions}.
Since $G=GL_n(\IC)$ is path connected it is possible to choose a smooth map 
$\chi_i:\vert\overline\ga_i\vert\rightarrow G$ such that 
$\chi_i=1$ on $\overline D_0\cap \vert\overline\ga_i\vert$ and
$\chi_i={^i\Phi_{0}}P_iC_i$ on $^i\Sect_{0}\cap \vert\overline\ga_i\vert$ for
$i=1,\ldots,m$.
Define $\al$ over $\vert\T\vert$ as follows:
$\al\vert_{\overline D_0}=0$ and for $i=1,\ldots,m$
$\al\vert_{\overline D_i}={\al_i}$
and
$\al\vert_{\vert\overline\ga_i\vert}=(d\chi_i)\chi_i^{-1}$.
It is easy to check these definitions 
agree on the overlaps and that when
the basis $^i\Phi_{0}$ is extended over $\vert\overline\ga_i\vert$ as
a solution of $\al$ then 
$\rho(\ga_{\wt p_ip_0}) = {^i\Phi_{0}^{-1}}\cdot\Phi$
on $\vert\overline\ga_i\vert$, where $\Phi$ is the basis which 
equals $1$ on $D_0$.

Now we must extend $\al$ to the rest of $\IP^1$.
First the product relation $\rho_m\cdots\rho_1=1$ ensures
that $\al$ has no monodromy around the boundary circle 
$\partial \vert\T\vert\cong S^1$, 
so that any local fundamental solution $\Psi$ extends to give a map
$\Psi:\partial \vert\T\vert\to G$.
Then the hypothesis that $\deg(\rho)=0$ 
ensures that this
loop $\Psi$ in $G$ is contractible.
To see this, firstly recall that the determinant map
$\det:G\to \IC^*$ expresses $G$ as a fibre bundle over
$\IC^*$, with fibres diffeomorphic to $SL_n(\IC)$, and that $SL_n(\IC)$
is simply connected.
Then, from the homotopy long exact sequence for fibrations, 
it follows that $\det$ induces
an isomorphism of fundamental groups:
$\pi_1(G)\cong\pi_1(\IC^*)\cong\IZ$.
Thus we need to see that the loop 
$\psi:= \det(\Psi):\partial \vert\T\vert\to \IC^*$
in the punctured complex plane does not wind around zero.
But the winding number of $\psi$ is 
$$\frac{1}{2\pi i}\oint_{\partial \vert\T\vert}\frac{d\psi}{\psi} = 
\frac{1}{2\pi i}\oint_{\partial \vert\T\vert}\tr(\al)$$
and the $C^\infty$ version of Cauchy's integral theorem
(see Lemma \ref{lem: CIT})
implies this is equal to $\sum \tr(\iL)=\deg(\rho)$ 
(using the flatness of $\al$ to deduce $d\tr(\al)=0$).

Thus the loop $\Psi$ in $G$ may be extended  
to a smooth map from the 
complement of $\vert\T\vert$ in $\IP^1$ to  $G$.
We then define $\al = (d\Psi)\Psi^{-1}$ on 
this complement and thereby obtain 
$\al\in\wt\A_{\fl}(\ba)$ having the desired monodromy data.
Hence the $C^\infty$ monodromy map is indeed surjective.

Next observe (from Theorem
\ref{thm: multisums} and Lemma \ref{C^infty canonical solutions}) that if 
$h\in\G_T$ and $\al'=h[\al]$ then the canonical
solutions of $\al$ and $\al'$ are related by:
$^i\Phi'_j=h\cdot{^i\Phi_j}\cdot t_i^{-1}$
where $t_i=h(a_i)\in T$.
The intertwining property and the fact that the $\G_1$ orbits are
contained in the fibres of $\wt \nu$ are then immediate from 
the definition of the Stokes representation 
in terms of the canonical solutions

The proof that the fibres are precisely the $\G_1$ orbits is much like
the proof of Proposition \ref{prop: im iff im}:
Suppose
$\al,\al'\in\wt\A_{\fl}(\ba)$ have the same monodromy data.
Let 
$\varphi:={^1\Phi'_{0}}({^1\Phi_{0}})^{-1}$
be the induced isomorphism between $\al$ and $\al'$ on $^1\Sect_{0}$.
Then $\varphi$ is single-valued when extended to 
$\IP^1\setminus\{a_1,\ldots,a_m\}$ as a solution of the induced
connection $\Hom(\al,\al')$ on the trivial bundle $\End(\IC^n)$.
(When $\varphi$ is extended around any loop 
$\ga$ based at $p_1$ it has no monodromy
since, when extended around this loop,
${^1\Phi'_{0}}$ and ${^1\Phi_{0}}$ are both multiplied on the
right by the same constant matrix.) 
Also, since the monodromy data encodes the transitions between the
various canonical fundamental solutions it follows that
$\varphi={^i\Phi'_j}({^i\Phi_j})^{-1}$ for any $i,j$.
Now observe (from Theorem
\ref{thm: multisums} and Lemma \ref{C^infty canonical solutions}) that
${^i\Phi'_j}({^i\Phi_j})^{-1}$ is asymptotic to $1$ at $a_i$ on some
sectorial neighbourhood of $^i\Sect_j$ ($j=1,\ldots,r_i,\ i=1,\ldots,m$).
It follows that $\varphi$ extends smoothly to $\IP^1$ and has Taylor
expansion $1$ at each $a_i$.
By construction $\al'=\varphi[\al]$ so $\al$ and $\al'$ are in the same
$\G_1$ orbit.
\epfms

\section{Symplectic Structure and Moment Map}  \label{sn: ss+mm}

In this section we observe that the well known Atiyah-Bott symplectic
structure on nonsingular connections naturally generalises to the
singular case we have been studying.
Moreover, as in the nonsingular case we find that
the curvature is a moment map for the action of the gauge group. 
Thus the moduli spaces of {\em flat}
singular connections with fixed expansions arise 
as infinite dimensional symplectic quotients.

The main technical difficulty here is that standard Sobolev/Banach
space methods cannot be used since we want to fix infinite-jets of
derivatives at the singular points $a_i\in\IP^1$.
Instead the infinite dimensional spaces here are naturally Fr\'echet
manifolds.
We will not use any deep properties of Fr\'echet spaces but do need a
topology and differential structure. 
(The explicitness of our situation
means we can get by without using an implicit function
theorem---the
monodromy description gives $\wt\A_\fl(\ba)/\G_1$
the structure of a complex manifold 
and local slices 
for this $\G_1$ action will be constructed directly.)
The reference used for Fr\'echet spaces is Treves \cite{Tre} and 
for Fr\'echet manifolds or Lie groups  see Hamilton \cite{Ham82} 
and Milnor \cite{Mil83}; we will give direct references to these works
rather than full details here.

\subsection*{The Atiyah-Bott Symplectic Structure on $\wt\A(\ba)$}

Let $E$ denote the trivial rank $n$ complex vector bundle over $\IP^1$
and
consider the complex vector space $\Omega^1[D](\IP^1,\End(E))$
of $n\times n$ matrices of global 
$C^\infty$ singular one-forms on $\IP^1$ with poles on $D$
(see Section \ref{sn: smooth}).
This is the space of sections of a $C^\infty$ vector bundle and so can
be given a  Fr\'echet topology in a standard way
(\cite{Ham82} p68). Now define $W$ to be the vector subspace 
$$W:= \left\{
\phi\in \Omega^1[D](\IP^1,\End(E)) \ \bigl\vert \ 
L_i(\phi)\in\lt\frac{dz_i}{z_i} \text{ for $i=1,\ldots,m$ } \right\}$$
of $\Omega^1[D](\IP^1,\End(E))$
of elements having Laurent expansion zero at each $i$, except for a
possibly nonzero, diagonal residue term.
This is a closed subspace\footnote{
since the $C^\infty$ Laurent expansion maps $L_i$ are {\em continuous} 
(if we put the topology of simple convergence of
coefficients on the formal power series ring 
which is the image of the Laurent expansion map $L_i$);
see \cite{Tre} p390, where this fact is used to prove 
E.Borel's theorem on the surjectivity of $L_i$.}
and so inherits a Fr\'echet topology.

\begin{lem}	\label{lem: Frechet mfd of conns}
The extended space $\wt\A(\ba)$ of singular connections 
is a complex Fr\'echet manifold and if $\al\in\wt\A(\ba)$ then
the tangent space to $\wt\A(\ba)$ at $\al$ is canonically isomorphic
to the complex Fr\'echet space $W$ defined above:
$T_\al\wt\A(\ba)\cong W.$
\end{lem}
\pf
If all $k_i\ge 2$ then
$\wt\A(\ba)$ is an affine space modelled on $W$:
If $\al_0\in\wt\A(\ba)$ then
$\wt\A(\ba)= \{ \al_0+\phi \ \bigl\vert \ \phi\in W \}.$
In general (some $k_i=1$), $\wt\A(\ba)$ is identified in this way
with an open subset of $W$ (recall the residues must be regular mod
$\IZ$):
if $\al_0\in\wt\A(\ba)$ then the map
$ \{ \al_0+\phi \ \bigl\vert \ \phi\in W \} \to \lt^m;$
$\al\mapsto \bigl(\res_iL_i(\al)\bigr)_{i=1}^m$
taking the residues is continuous and $\wt\A(\ba)$ is the inverse
image of an open subset of $\lt^m$.
Thus $\wt\A(\ba)$ is identified with an open subset of $W$;
it is thus a Fr\'echet manifold (with just one coordinate chart)
and the tangent spaces are canonically identified with $W$ as in the
finite dimensional case (see discussion \cite{Mil83} p1030).
\epf

Thus following Atiyah-Bott \cite{AB82} we can define a two-form
\begin{equation}	\label{AB form}
\omega_\al(\phi,\psi) := \frac{1}{2\pi i}\ip \tr(\phi\wedge\psi)
\end{equation}
on $\wt\A(\ba)$,
where $\al\in \wt\A(\ba)$ and $ \phi,\psi\in T_\al\wt\A(\ba)$.
This integral is well defined since the two-form $\tr(\phi\wedge\psi)$
on $\IP^1$ is nonsingular;
its Laurent expansion at $a_i$ is a $(2,0)$ form and so zero. 
Thus $\omega_\al$ is a skew-symmetric complex 
bilinear form on the tangent space $T_\al\wt\A(\ba)$.
It is nondegenerate in the sense that 
if $\omega_\al(\phi,\psi)=0$ for all $\psi$ then $\phi=0$
(if $\phi \ne 0$ then $\phi$ is nonzero at some point 
$p \ne a_1,\ldots,a_m$ 
and it is easy then to construct $\psi$ vanishing
outside a neighbourhood of
$p$ and such that  $\omega_\al(\phi,\psi)\ne 0$).
Also $\omega_\al$ is continuous as a map $W\times W \to \IC$, 
since it is continuous in each factor, and 
(for Fr\'echet spaces) such `separately continuous' bilinear maps
are continuous (\cite{Tre} p354).
Finally the right-hand side of (\ref{AB form}) is independent of
$\al$, so $\omega$ is a constant two-form on $\wt\A(\ba)$ and in
particular it is closed.
Owing to these properties we will say 
$\omega$ is a complex {\em symplectic} form on $\wt\A(\ba)$.
(See for example Kobayashi \cite{Kob} for a discussion of the more
well-known theory of symplectic {\em Banach} manifolds.) 

\subsection*{Group Actions}

First, the full gauge group 
$\G:= GL_n(C^\infty(\IP^1))$
is a Fr\'echet Lie group; that is, it is a Fr\'echet manifold such
that the group operations $g,h\mapsto g\cdot h$ and $g\mapsto g^{-1}$
are $C^\infty$ maps (see \cite{Mil83} Example 1.3).
$\G$ is locally modelled on the Fr\'echet space 
$\Lia(\G) := C^\infty(\IP^1,\gl)$
and has a complex analytic structure coming from the
exponential map
$\exp : \Lia(\G)\to\G;$
$x\mapsto \exp(x)$
which is defined pointwise in terms of
the exponential map for $G$.
This implies $\Lia(\G)$ is
a canonical coordinate chart for $\G$ in a neighbourhood  of the
identity since $\exp$ has a local inverse $g\mapsto \log(g)$ (also
defined pointwise). 
In particular $\Lia(\G)$ is so identified with the tangent space to
$\G$ at the identity; the Lie algebra of $\G$.

The group we are really interested in here is $\G_1$, the subgroup of
$\G$ consisting of elements $g\in \G$ having Taylor expansion 
$1$ at each $a_i\in\IP^1$.
As above, the Taylor expansion maps are continuous and so 
$\G_1$ (the intersection of their kernels) is a closed subgroup of
$\G$.
It follows that $\G_1$ is a complex Fr\'echet Lie group with Lie algebra
$$\Lia(\G_1) := \{ x\in\Lia(\G) \ \bigl\vert \ 
L_i(x)=0 \text{ for $i=1,\ldots, m$} \}$$
where $L_i$ is the Taylor expansion map at $a_i$.
(The same statements also hold for $\G_T$ except now
$\Lia(\G_T) := \{ x\in\Lia(\G) \ \bigl\vert \ 
L_i(x)\in\lt \text{ for $i=1,\ldots, m$} \}$.)

\begin{lem}	\label{inf fvf}
The groups $\G_1$ and $\G_T$ act holomorphically on $\wt\A(\ba)$ by gauge
transformations and the fundamental vector field of $x\in \Lia(\G_T)$
takes the value
$ -d_\al x \in T_\al\wt\A(\ba)$
at $\al\in\wt\A(\ba)$, where 
$d_\al$
is the singular connection on $\End(E)$ induced from $\al$.
\end{lem}
\pf
First the action map 
$\G_T\times\wt\A(\ba)\rightarrow \wt\A(\ba);$
$(g,\al)\mapsto g\al g^{-1}+ (dg)g^{-1}$
can be factored into simpler maps each of which is 
holomorphic (see \cite{Ham82}).
By convention the fundamental vector field is minus the tangent field to the
flow generated by $x$, which may be calculated using the exponential map
for $\G_T$.  
\epf

\subsection*{The Curvature is a Moment Map}

It is clear that the action of $\G_T$ on $\wt\A(\ba)$ preserves
the symplectic form $\omega$: If $g\in \G_T$ 
and $\al\in \wt\A(\ba)$ then the derivative of the action of $g$
is simply conjugation:
$$(g[\cdot])_*:T_\al\wt\A(\ba)\to T_{g[\al]}\wt\A(\ba);\quad
\phi\mapsto g\phi g^{-1}.$$ 
and so $\omega$ is preserved as
$\tr(\phi\wedge \psi) = \tr(g\phi g^{-1}\wedge g\psi g^{-1})$
for any $\phi,\psi\in T_\al\wt\A(\ba)$.

More interestingly, this action is Hamiltonian.
If we firstly look at the smaller group $\G_1$, then,
as observed by Atiyah and Bott in the nonsingular case,
the curvature is a moment map.
To start with observe:
\begin{lem}	\label{inf dif F}
The curvature map 
$\F:\wt\A(\ba)\longrightarrow \Omega^2(\IP^1,\End(E))$
is an infinitely differentiable (even holomorphic) map to
the Fr\'echet space of $\End(E)$ valued {\em nonsingular} two-forms on
$\IP^1$.
The derivative of $\F$ at $\al\in\wt\A(\ba)$ is
$$(d\F)_\al : T_\al\wt\A(\ba)\to \Omega^2(\IP^1,\End(E));\qquad
\phi\longmapsto -d_\al\phi$$
where $\phi\in T_\al\wt\A(\ba)=W$ and
$d_\al : \Omega^1[D](\IP^1,\End(E)) \to \Omega^2[2D](\IP^1,\End(E))$
is the operator naturally induced from the singular connection $\al$.
\end{lem} 
\pf
Recall the curvature is given explicitly by 
$\F(\al)= -d\al+\al\wedge\al$
and observe (by looking at Laurent expansions and using the division lemma) 
that this is a matrix of {\em nonsingular} two-forms.
That $\F$ is $C^\infty$ with the stated derivative follows from basic
facts about calculus on Fr\'echet spaces
(see \cite{Ham82} Part I).
\epf

Next there is a natural inclusion from $\Omega^2(\IP^1,\End(E))$ to
the dual of the Lie algebra of $\G_1$ given by taking the trace and
integrating: 
$$\iota:\Omega^2(\IP^1,\End(E)) \rightarrow \Lia(\G_1)^*;\qquad
\F(\al)\mapsto \left(x\mapsto \frac{1}{2\pi i}\ip \tr(\F(\al)x)\right)$$
where $x\in\Lia(\G_1)$ is a matrix of functions on $\IP^1$.
Using this inclusion we will regard $\F$ as a map into the dual of the
Lie algebra of the group.
We then have

\begin{prop}	\label{prop: F is mmap}
The curvature
$\F : \wt\A(\ba)\longrightarrow \Lia(\G_1)^*$
is 
an equivariant moment map for the $\G_1$ action on the extended space
$\wt\A(\ba)$ of singular connections.
\end{prop}
\pf
Everything has been set up so that the arguments from the nonsingular
case still work, as we will now show. 
Given $x\in\Lia(\G_1)$, define a (Hamiltonian) 
function $H_x$ on $\wt\A(\ba)$ to be the $x$ component of 
$\F$:
$$H_x:=\langle \F , x \rangle:\wt\A(\ba)\to \IC;\qquad
H_x(\al)=\frac{1}{2\pi i}\ip \tr(\F(\al)x)$$
where the angled brackets denote the natural pairing between
$\Lia(\G_1)$ and its dual.
We need to show that the fundamental vector field of $x$ is the
Hamiltonian vector field of $H_x$, i.e. that
$(dH_x)_\al= \omega_\al(\cdot,-d_\al x)$
as elements of $T^*_\al\wt\A(\ba)$.
Now if $\phi\in T_\al\wt\A(\ba)$ then 
\begin{equation}	\label{curvm2}
(dH_x)_\al(\phi)= -\frac{1}{2\pi i}\ip\tr ((d_\al \phi)x)
\end{equation}
from Lemma \ref{inf dif F} and the chain rule.
Now observe that $\tr(\phi x)$ is a {\em nonsingular} one-form on
$\IP^1$ (as $L_i(x)=0$ for all $i$).
Therefore Stokes' theorem implies
$d\tr(\phi x) = \tr((d_\al \phi)x) - \tr(\phi\wedge d_\al x)$ 
integrates to zero over $\IP^1$.
Hence (\ref{curvm2}) becomes 
$$(dH_x)_\al(\phi) = -\frac{1}{2\pi i}\ip\tr(\phi \wedge d_\al x) =
\omega_\al(\phi,-d_\al x)$$
proving that the curvature is indeed a moment map.

The equivariance follows directly from the definition of the coadjoint
action of $\G_1$: If $x\in\Lia(\G_1)$ then 
$
\langle \Ad^*_g(\F(\al)), x\rangle  := 
\langle \F(\al), \Ad_{g^{-1}}(x) \rangle =
\langle \F(g[\al]), x \rangle 
$
using the fact that 
$\tr(\F(\al) g^{-1}xg)  = \tr( \F(g[\al]) x)$.
\epf

Thus the subset of flat connections is the preimage of zero under the
moment map:
$\wt\A_{\fl}(\ba) = \F^{-1}(0).$
Therefore, at least in a formal sense, the moduli space is a
symplectic quotient:
$\wt\A_{\fl}(\ba)/\G_1 = \F^{-1}(0)/\G_1.$
(Recall $\wt\A_{\fl}(\ba)/\G_1$ was identified in Section
\ref{sn: smooth} with the space $\wt M_0(\ba)$ of monodromy data,
analogously to the non-singular case.)
In the next section we will show that this prescription does
define a genuine symplectic structure on at least the dense open
subset of $\wt M_0(\ba)$ which is the image of the extended monodromy
map $\wt \nu$.

\subsection*{Torus Actions}

To end this section we consider the action of the larger group $\G_T$
on the extended space of singular connections $\wt\A(\ba)$.
This action is also Hamiltonian:

\begin{prop}
Let $\mu : \wt\A(\ba)\longrightarrow \Lia(\G_T)^*$ be the map 
given by taking the curvature together with the residue at each $a_i$:
If $x\in\Lia(\G_T)$ and $\al\in \wt\A(\ba)$
$$\langle \mu(\al) , x \rangle :=
\frac{1}{(2\pi \sqrt{-1})}\ip \tr\left(\F(\al) x\right) - 
\sum_{i=1}^m\res_iL_i(\tr (\al x)).$$
Then $\mu$ is an equivariant moment map for the $\G_T$ action on 
$\wt\A(\ba)$.
\end{prop}
\pf
For any $x\in\Lia(\G_T)$ 
define the function $H_x:\wt\A(\ba)\to\IC$ to be the $x$ component of
$\mu$:
$H_x(\al):= \langle \mu(\al) , x \rangle.$
If $\phi\in T_\al\wt\A(\ba)$ then
\begin{equation}	\label{dH}
(dH_x)_\al(\phi)= 
-\frac{1}{(2\pi \sqrt{-1})}\ip\tr ((d_\al \phi)x)- 
\sum_{i}\res_iL_i(\tr (\phi x)).
\end{equation}
Our task is to show $\omega_\al(\phi,-d_\al x)=(dH_x)_\al(\phi)$.
We do this by using the $C^\infty$ Cauchy integral theorem
(see Lemma \ref{lem: CIT}).
Recall $\phi$ is a matrix of $C^\infty$ one-forms on 
$\IP^1$ with (at worst)
first order poles in its $(1,0)$ part at each $a_i$.
Also $x\in\Lia(\G_T)$ is a matrix of functions on $\IP^1$ and has
Taylor expansion equal to a constant diagonal matrix at each $a_i$.
Thus for each $i$ we can choose a $C^\infty$ function 
$f_i:\IP^1\to\IC$ which vanishes outside $D_i$, such that
$\tr(\phi x) = \theta + f_1{dz_1/z_1}+\cdots f_m{dz_m/z_m}$
for some {\em nonsingular} one-form $\theta$ on $\IP^1$.
Thus
$d\tr(\phi x) = d\theta - \sum_i 
\frac{\partial f_i}{\partial \bar z_i}
\frac{dz_i\wedge d\bar z_i}{z_i}$
and so by Stokes' theorem and Cauchy's integral theorem:
$$\ip d\tr(\phi x) = \sum_i \int_{\overline D_i}
\frac{\partial f_i}{\partial \bar z_i}
\frac{dz_i\wedge d\bar z_i}{z_i}=-(2\pi\sqrt{-1})\sum_i f_i(a_i).$$
(Note $f_i(a_i)=\res_iL_i(\tr (\phi x))$.)
Then the equality $\omega_\al(\phi,-d_\al x)=(dH_x)_\al(\phi)$ follows
from the fact that $ d\tr(\phi x) = \tr(d_\al(\phi x)) = 
\tr((d_\al\phi)x) - \tr(\phi\wedge d_\al x)$.
The equivariance follows exactly as before 
since $\G_T/\G_1\cong T^m$ is Abelian.\ \  
\epf

Instead we could do the reduction in stages, and consider the $T^m$
action on $\wt\A_\fl(\ba)/\G_1$. This matches up with the
Hamiltonian $T^m$ actions considered in Section \ref{sn: triv},
since the residues above are the exponents of formal monodromy 
$\iL$.

%% file: smid-mmsp.tex
\section{The Monodromy Map is Symplectic} \label{sn: mm is sp} \ 

Most of the story so far can be summarised in the commutative
diagram:
\begin{equation}	\label{main diagram}
\begin{array}{cccccc}
&&  \wt\M(\ba) & \mapright{\cong} & \wt\A_{\fl}(\ba)/\G_1 &\\
&&  \bigcup && \mapdown{\cong}&\\
\wt O_1\times\cdots\times\wt O_m\spq G &\cong&
  \wt\M^*(\ba) & \mapright{\wt\nu} & \wt M_0(\ba). &\hphantom{\longrightarrow}
\end{array}
\end{equation}
The extended moduli space $\wt\M^*(\ba)$ was defined in
Section \ref{sn: triv} to be the set of isomorphism 
classes of compatibly framed meromorphic connections on trivial rank
$n$ vector bundles with irregular type $\ba$.
It was given an intrinsic
complex symplectic structure  explicitly in terms of (finite dimensional)
coadjoint orbits and cotangent bundles.
The extended monodromy manifold $\wt M(\ba)$ was defined as
the set of isomorphism classes of Stokes representations 
 and looks like a multiplicative version of
$\wt\M^*(\ba)$ (when both are described explicitly).
$\wt M_0(\ba)$ is the degree zero component of $\wt M(\ba)$ and
was identified with the set of $\G_1$ orbits in the extended space 
$\wt\A_{\fl}(\ba)$ of flat $C^\infty$ singular connections. 
Moreover the curvature was shown to be a moment map for the action
of the gauge
group $\G_1$ on the symplectic Fr\'echet manifold $\wt\A(\ba)$,
so that (formally) $\wt\A_\fl(\ba)/\G_1$ is a complex
symplectic quotient.
$\wt \M(\ba)$ has the same definition as $\wt\M^*(\ba)$ except
with the word `trivial' replaced by `degree zero'.  
The act of taking monodromy data defines both the right-hand
isomorphism 
in the diagram and the monodromy map $\wt\nu$,
which is a
biholomorphic map onto its image (a dense open submanifold of
$\wt M_0(\ba)$).

Basically the bottom line appears in the work 
\cite{JMU81} of Jimbo, Miwa and Ueno
but the symplectic structures and the rest of the diagram do not.
The torus $T^m\cong(\IC^*)^{nm}$ acts on each space in (\ref{main diagram})
and these actions are intertwined by all the maps. The non-extended
picture arises by taking the symplectic quotient 
(fixing the exponents of formal monodromy and 
quotienting by $T^m$).
We then obtain another commutative diagram as above but with all
the tildes removed and $\G_1$ replaced by $\G_T$.

In this section we show that the symplectic
structure on $\wt\A(\ba)$ does induce a symplectic structure
on (at least) the dense open submanifold of $\wt M_0(\ba)$ that is
the image of the monodromy map $\wt \nu$, and
that this symplectic form pulls back 
along $\wt \nu$ to the explicit symplectic
form on $\wt\M^*(\ba)$.
In other words we will prove:
\begin{thm}	\label{thm mmap is sp}
The monodromy map $\wt\nu$ is symplectic.
\end{thm}
Analogous results have been proved in the 
simple pole case independently by Hitchin \cite{Hit95} and by
Iwasaki \cite{Iwa91, Iwa92}.
(Note that Iwasaki considers only
$PSL_2(\IC)$ Fuchsian equations, but he does so
over (fixed) arbitrary
genus Riemann surfaces.)
\subsection*{Factorising the Monodromy Map}

Recall from Proposition \ref{prop: global bijn 1} how the 
isomorphism at the top of the above diagram 
arose: a meromorphic connection gives rise to a $\G$ orbit
of $C^\infty$ 
singular connections and we consider the subset with fixed
Laurent expansion at each $a_i$ to define the map.
In other words we can choose $g\in\G$ to `straighten' a meromorphic
connection to have fixed $C^\infty$ Laurent expansions
and thereby
specify an element of $\wt\A_\fl(\ba)$.
Here we show that this straightening procedure can be carried out
for a family of connections all at the same time, and so the 
monodromy map factorises through $\wt\A_\fl(\ba)$.

As usual we fix the data $\ba$ 
consisting of an effective divisor $D=\sum k_i(a_i)$ and diagonal
generic connections germs $d-\iAo$. 
Also choose a coordinate $z$ to identify $\IP^1$ with $\IC\cup\infty$
such that each $a_i$ is finite and 
let $D_1,\ldots,D_m\subset \IP^1$ be disjoint open disks
with $a_i\in D_i$, so that $z_i:=z-a_i$ is a coordinate on $D_i$. 
\begin{prop}	\label{prop: straightening}
Let $U\subset\wt\M^*(\ba)$ be an open subset.
Then
there exists a universal family 
$d_{\IP^1}-A$ of meromorphic connections
on the trivial bundle over
$\IP^1$ (with compatible framings 
${\bf g}= ({^1\!g_0},\ldots,{^m\!g_0})$)
parameterised by $u\in U$
and a family of smooth bundle automorphisms
$g\in GL_n(C^\infty(U\times \IP^1))$
such that for each $u\in U$ and each $i=1,2,\ldots,m$:

$\bullet$ $g(u,a_i)\in GL_n(\IC)$ is the compatible 
framing ${^i\!g_0}(u)$ at $a_i$ specified by 
$u\in U$,

$\bullet$ The singular connection $\al(u):=g(u)[A(u)]$ on $\IP^1$
has Laurent expansion $\iAo+{^i\!R(u)}$ at $a_i\in\IP^1$, 
where 
${^i\!R}(u)=(\iL-\iLo)dz/z_i$, $\iLo=\res_i(\iAo)$ and 
$\iL$ is the exponent of formal
monodromy of $(d-A(u),{\bf g})$ at $a_i$,

$\bullet$
If $z\in\IP^1\setminus(D_1\cup\cdots\cup D_m)$ then $g(u,z)=1$.
\end{prop}

\pf
The construction of the universal family is immediate from the proof
of Proposition \ref{prop: main pp}:
Using the coordinate choices, $\wt\M^*(\ba)$ can be identified with
the submanifold of 
$\mu_G^{-1}(0)\subset\wt O_1\times\cdots\times\wt O_m$
which has ${^1\!g_0}=1$.
This subset was identified as a set of matrices of meromorphic
one-forms $A(u)$ on $\IP^1$, together with compatible framings 
$\bf g$.
(Although we do not need this fact, it is easy to check that the
family $(d-A,{\bf g})$ of compatibly
framed connections on the trivial bundle has the appropriate
universal property for $\wt \M^*(\ba)$; it is a {\em fine} moduli
space.)

Now consider the Laurent expansion $L_i(A)\in
\End_n\bigl(\IC\{z_i\}\otimes\cO(U)\bigr)dz/z_i^\mki$
of $A$ at
$a_i\in\IP^1$,
where the coefficients are now holomorphic functions on $U$.
(If $u\in U$ then  
$L_i(A)(u)= L_i(A(u))$ as elements of $\End_n(\IC\{z_i\})dz/z_i^{k_i}$.)
Recall from Lemma \ref{lem: get F} that the compatible framings
determine formal isomorphisms: In the relative case here this means
that, for each $i$, there is a unique invertible matrix 
${^i\widehat g}\in GL_n(\IC\flb z_i\frb\otimes\cO(U))$
of formal power series with coefficients in $\cO(U)$
which agrees with the compatible framing at $a_i$
and for each $u\in U$ satisfies:
\begin{equation} 	\label{formal ims}
{^i\widehat g}(u)[A(u)]= {^i\!A^0} + {^i\!R}(u)\in 
\End_n\bigl(\IC\flb z_i \frb\bigr)\frac{dz}{z_i^\mki}
\end{equation}
with ${^i\!R}(u)$ as in the statement of the proposition.
(The algorithm to construct such ${^i\widehat g}$'s is as before;
it works with coefficients in $\cO(U)$.)

The crucial step is to now use E. Borel's result 
that the Taylor expansion map is surjective
(Theorem \ref{thm: borel1} above).
Applying this to
each matrix entry of each ${^i\widehat g}$ in turn for
$i=1,\ldots,m$
gives matrices of functions
${^i\!g}\in \End_n(C^\infty(U\times \overline D_i))$ 
such that for each $u\in U$
the Taylor expansion of ${^i\!g}$ at $a_i$ is
${^i\widehat g(u)}$.
Since $\det{^i\!g}(u,a_i)=\det{^i\!g_0}(u)$ 
is nonzero for all $u\in U$, there is a 
neighbourhood of $U\times\{a_i\}\subset U\times \IP^1$ 
throughout which $\det({^i\!g})$ is nonzero. 
It follows (as $GL_n(\IC)$ is connected) that there is a smooth 
bundle automorphism
$g\in GL_n(C^\infty(U\times\IP^1))$
that equals ${^i\!g}$ in some neighbourhood of
$U\times\{a_i\}\subset U\times \IP^1$
for each $i$ and equals $1$ outside $U\times(D_1\cup\cdots\cup D_m)$.
In particular $g$ has the desired Taylor expansions at each $a_i$ 
so that $\al=g[A]$ has the desired $C^\infty$ Laurent expansions by
construction.
\epf

\begin{cor}\footnote{This statement is corrected from the published version. Thanks are due to B. Malgrange for querying this. This  does not affect the rest of the article.}	\label{cor: factorisation of mm}
The monodromy map $\wt \nu$ factorises through
$\wt\A_\fl(\ba)$:
It is possible to choose a 
map $\wh\nu$ from the extended moduli space $\wt\M^*(\ba)$
to the extended space of flat singular connections
determined by $\ba$
such that the diagram:
$$
\begin{array}{ccccc}
  \wt\M^*(\ba) & \mapright{\wh \nu} & \wt\A_\fl(\ba) &
\hookmapright{i} & \wt\A(\ba) \\
  \Bigl\| && \mapdown{/\G_1} && \\
  \wt\M^*(\ba) & \mapright{\wt \nu} & \wt M_0(\ba)&&
\end{array}
$$
commutes and the composition $i\circ \wh \nu$ into the Fr\'echet
manifold $\wt\A(\ba)$ is $C^\infty$.
\end{cor}

\pf
Construct $g$ as in Proposition \ref{prop: straightening}
with $U=\wt\M^*(\ba)$ and then define $\wh \nu(u)= g(u)[A(u)]$ for all
$u\in\wt\M^*(\ba)$. 
All that remains is to see that the composition
$i\circ\wh\nu$ is $C^\infty$.
Recall (from Lemma \ref{lem: Frechet mfd of conns}) that by choosing a
basepoint $\wt\A(\ba)$ is identified with a Fr\'echet submanifold of
the Fr\'echet space $\Omega^1[D](\IP^1,\End(E))$ of matrices of
$C^\infty$ one-forms with poles on the divisor $D$.
Thus we must prove that the map
$\wt\M^*(\ba)\rightarrow \Omega^1[D](\IP^1,\End(E));$
$u\mapsto \al(u):=g(u)[A(u)]$
is $C^\infty$.
Now if $u_0\in\wt\M^*(\ba)$ and 
$W_0\in T_{u_0}\wt\M^*(\ba)$ is a tangent vector at $u_0$,
then we will denote the partial derivative of $\al$ along $W_0$ by
$W_0(\al)\in \Omega^1[D](\IP^1,\End(E)).$
Here we think of $\al$ as a section of the $C^\infty$ vector bundle
$\pi^*(\End_n(\Omega^1[D]))$ over $\IP^1\times U$ 
(where $\pi:\IP^1\times U\to\IP^1$ is the obvious projection).
This vector bundle is trivial in the $U$ directions so the partial
derivative makes sense. 
(Concretely, local sections are of the form $\sum h_i\theta_i$ for
$C^\infty$ functions $h_i$ on $U$ and sections $\theta_i$ of 
$\End_n(\Omega^1[D])$.
Then $W_0$ differentiates just the $h_i$'s: $W_0(\sum h_i\theta_i)=
\sum W_0(h_i)\theta_i$.)
It then follows from basic facts about calculus on Fr\'echet spaces
that the map $i\circ\wh\nu$ is $C^\infty$ and has derivative $W_0(\al)$
along $W_0$ at $u_0$.
(This can be deduced from Examples 3.1.6 and 3.1.7 in \cite{Ham82}.)
\epf

\subsection*{Main Proof}\

\pfms [of Theorem \ref{thm mmap is sp}].\  
Choose $g$ as in Corollary \ref{cor: factorisation of mm} above and
let $\wh \nu:\wt\M^*(\ba)\to \wt\A_\fl(\ba)$
be the corresponding lift of the monodromy map.
It is sufficient for us to prove that the composite map 
$i\circ\wh \nu:\wt\M^*(\ba)\to \wt\A(\ba)$ is symplectic.
This is because the symplectic form on $\wt M_0(\ba)$ is defined
locally as $(i\circ s)^*\omega_{\wt\A(\ba)}$ for any local slice
$s:\wt M_0(\ba)\to\wt\A_\fl(\ba)$ of the $\G_1$ action.
But, over the subset $\wt\nu(\wt\M^*(\ba))$, such a slice is given by
$\wh\nu\circ\wt\nu^{-1}$. 
Thus $\wt\nu^*\omega_{\wt M_0(\ba)} = (i\circ\wh\nu)^*\omega_{\wt\A(\ba)}$
and, if $i\circ\wh\nu$ is symplectic, this is $\omega_{\wt\M^*(\ba)}$.

Now choose a point $u_0\in\wt\M^*(\ba)$ and two tangent vectors
$W_1,W_2\in T_{u_0}\wt\M^*(\ba)$.
Define two matrices of singular one-forms on $\IP^1$,
$\phi_j:=W_j(\al)\in\End_n(\Omega^1[D](\IP^1))$ $(j=1,2)$, to  be the 
corresponding partial derivatives of $\al(u):=g(u)[A(u)]$.
As in the proof of Corollary \ref{cor: factorisation of mm}, 
$\phi_j$ is the derivative $(i\circ\wh\nu)_*(W_j)$ of the map $i\circ\wh\nu$
along $W_j$.
Therefore what we must prove is:
\begin{equation}	\label{eq: wwmp}
\frac{1}{2\pi i}
\ip\tr(\phi_1\wedge\phi_2)=\omega_{\wt\M^*(\ba)}(W_1,W_2).
\end{equation}
The first step is to obtain a formula for the right-hand side in terms
of $g$.
This comes from Lemma \ref{lem: taction} since, by construction, the
first $k_i$ terms of the Taylor expansion of $g$ at $a_i$ give a
section of the $i$th `winding map' $w$.
For $j=1,2$ 
define $\idLj= W_j(\iL)\in\lt$ where $\iL$ is the $i$th exponent of
formal monodromy (which is regarded as a $\lt$-valued function on
$\wt\M^*(\ba)$). 
Let  $\idRj := \idLj dz/z_i$
and denote the derivatives of $g$ as
$\dot g_j := W_j(g)\in \End_n(C^\infty(\IP^1))$.
Then according to Lemma \ref{lem: taction}, if we define
${^i\!X_j}\in\g_{k_i}$ to be the first $k_i$ terms in the Taylor
expansion of $g(u_0)^{-1}\dot g_j$ at $a_i$ then
\begin{equation}	\label{eqn: 1st exprn}
\omega_{\wt\M^*(\ba)}(W_1,W_2)=
\sum_{i=1}^m \left(
  \langle {^i\!\dot R_1},{^i\!\wt X_2} \rangle
- \langle {^i\!\dot R_2},{^i\!\wt X_1} \rangle 
+ \langle {^i\!A(u_0)},[{^i\!X_1},{^i\!X_2}] \rangle \right)\ \ \ 
\end{equation}
where 
${^i\!A}$ is the Laurent expansion of $A$ at $a_i$ and
${^i\!\wt X_j}=
{^i\!g_0(u_0)}\cdot{^i\! X_j}\cdot {^i\!g_0(u_0)}^{-1}\in\g_{k_i}$
for $j=1,2$ and  $i=1,\ldots,m$.

Now we will calculate the left-hand side of (\ref{eq: wwmp}).
First observe that the two-form $\tr(\phi_1\wedge\phi_2)$ on $\IP^1$
is non-singular. Indeed the $C^\infty$
Laurent expansion of $\phi_j$ at $a_i$ is
$\idRj$, and so the expansion of $\tr(\phi_1\wedge\phi_2)$ is a $(2,0)$
form and so zero. 
Then the division lemma implies $\tr(\phi_1\wedge\phi_2)$ is non-singular.

Next, by differentiating the expression $g(u)[A(u)]$ for $\al$ along
$W_j$ we find 
\begin{equation}	\label{eqn: phi tilde}
\phi_j=g(u_0)\cdot  \wt\phi_j \cdot g(u_0)^{-1},\qquad 
\wt\phi_j := \dot A_j+ d_{A(u_0)}\left(g(u_0)^{-1}\dot g_j\right)
\end{equation}
for $j=1,2$, where $\dot A_j:= W_j(A(u))$ and $g(u)=g(u,\cdot)\in\G$.
(Note that this formula is the basic reason why the `straightening' procedure
makes the Atiyah-Bott formula (\ref{AB form}) non-trivial in this situation.)
In particular we have $\tr(\phi_1\wedge\phi_2)=\tr(\wt\phi_1\wedge\wt\phi_2)$.
Observe this two-form is zero outside of the disks $D_i$,
since $\dot g_j$ is zero there and each $\dot A_j$ has type $(1,0)$.
It follows that the integral splits up into integrals over the closed
disks:
\begin{equation}	\label{int splits}
\ip \tr(\phi_1\wedge\phi_2)=
\sum_{i=1}^m\int_{\overline D_i}\tr(\wt\phi_1\wedge\wt\phi_2).
\end{equation}
We break each term in this sum into two pieces, 
using the definition (\ref{eqn: phi tilde}) of $\wt\phi_2$:
$\tr(\wt\phi_1\wedge\wt\phi_2)=
\tr(\wt\phi_1\wedge\dot A_2)+
\tr\left(\wt\phi_1\wedge d_{A(u_0)}(g(u_0)^{-1}\dot g_2)\right).$
Therefore by comparing with the expression (\ref{eqn: 1st exprn}),  
the theorem 
now follows immediately from:
\begin{claim}

$1)$ 
$
\frac{1}{(2\pi\sqrt{-1})}\int_{\overline D_i}\tr(\wt\phi_1\wedge\dot A_2)=
\langle {^i\!A(u_0)},[{^i\!X_1},{^i\!X_2}] \rangle
- \langle  {^i\!\dot R_2},{^i\!\wt X_1} \rangle,
$

$2)$
$\frac{1}{(2\pi\sqrt{-1})}\int_{\overline D_i}
\tr(\wt\phi_1\wedge d_{A(u_0)}(g(u_0)^{-1}\dot g_2))=
\langle {^i\!\dot R_1},{^i\!\wt X_2} \rangle.$
\end{claim}
The basic tool we will use to evaluate these integrals is:
\begin{lem} (Modified $C^\infty$ Cauchy Integral Theorem). \label{lem: CIT}
Let $k$ be a nonnegative integer, $a\in\IC$ a complex number and 
$D_a$ a disk in $\IC$ containing the point $a$.
Suppose $f\in C^\infty(\overline D_a)$ and 
$\left(
\frac{\partial f}{\partial \bar z}
\right)/(z-a)^k\in C^\infty(\overline D_a)$.
Then
$\left(
\frac{\partial f}{\partial \bar z}
\right)          
\frac{dz\wedge d\bar z}{(z-a)^{k+1}}$
is absolutely integrable over $\overline D_a$ and
$$\frac{(2\pi i)}{k!}\frac{\partial^k f}{\partial z^k}(a)=
\oint_{\partial \overline D_a}\frac{f(z)dz}{(z-a)^{k+1}}+
\int_{\overline D_a}\frac{\partial f}{\partial \bar z}
\frac{dz\wedge d\bar z}{(z-a)^{k+1}}
$$
where the line integral is taken in an anti-clockwise direction.
\end{lem} 
\pf
The $k=0$ case is the usual $C^\infty$ Cauchy integral theorem, see
\cite{GH94} p2.
Differentiating with respect to $a$ then gives the above result:  
we may reorder the integration and
differentiation due to the 
absolute integrability. 
\epf

Part $1)$ of the claim arises as follows.
Since $\dot A_2$ is a matrix of meromorphic one-forms  we have 
$\tr(\wt\phi_1\wedge\dot A_2)=\tr(\wt\phi_1^{(0,1)}\wedge\dot A_2),$
and from (\ref{eqn: phi tilde}):
$$\wt\phi_1^{(0,1)}= \bar\partial ( g(u_0)^{-1}\dot g_1)=
\frac{\partial( g(u_0)^{-1}\dot g_1)}{\partial\bar z}d\bar z.$$
Also $\dot A_2$ has a pole of order at most $\mki$ at $a_i$ and so
we can define a smooth function on $\overline D_i$, 
$f\in C^\infty(\overline D_i)$, 
by the prescription
$f dz = (z-a_i)^\mki\cdot \tr \bigl(g(u_0)^{-1}\dot g_1 
\dot A_2\bigr)$ on $\overline D_i.$
By taking the exterior derivative of both sides of this
and dividing through by $(z-a_i)^\mki$ we deduce
$$\tr(\wt \phi_1\wedge\dot A_2)=
-\frac{\partial f}{\partial\bar z}\frac{dz\wedge d\bar z}
{(z-a_i)^\mki}\qquad\text{on $\overline D_i$},$$
where the minus sign occurs since we have reversed the order of $dz$
and $d\bar z$.
Observe that the Taylor expansion of $f dz$ at $a_i$ has no 
terms containing $\bar z_i$.
Thus $\partial f/\partial \bar z$ has zero Taylor expansion at $a_i$
and in particular using the division lemma we see 
$f$ satisfies the conditions in Lemma \ref{lem: CIT}.
Also $f$ is zero on the boundary $\partial\overline D_i$ since $\dot g_1$
is zero there.
Therefore Cauchy's integral theorem gives
\begin{equation}	\label{eqn: CIT use}
\frac{1}{(2\pi\sqrt{-1})}\int_{\overline D_i}\tr(\wt\phi_1\wedge\dot A_2)=
-\frac{1}{k!}
\frac{\partial^k f}{\partial z^k}(a_i)\qquad
\text{with $k=\mki-1$.}
\end{equation}
This value is just 
$-\Res_i(f dz/(z-a_i)^\mki)=
-\res_i(\tr \bigl(g(u_0)^{-1}\dot g_1\dot A_2\bigr))$, where `residue'
just means taking 
the coefficient of $dz/z_i$ in the ($C^\infty$) Laurent expansion.
This last expression only involves the principal part of $\dot A_2$ at $a_i$
and the first $\mki$ terms of the Taylor expansion of 
$g(u_0)^{-1}\dot g_1$.
By definition these first $\mki$ terms are
given by ${^i\!X_1}$.
Also, by construction, 
the principal part of $A$ at $a_i$ is the same as the principal
part of 
$g(u)^{-1}\bigl({^i\!A^0}+(\iL-\iLo)dz/z_i\bigr)g(u)$.
It follows directly that
$$\PP_i(\dot A_2)=\PP_i W_2(A(u))=
[{^i\!A(u_0)},{^i\!X_2}]+{^i\!g_0(u_0)}^{-1}\cdot{^i\!\dot R_2}\cdot
{^i\!g_0(u_0)}.$$
Statement $1)$ of the claim is now immediate, upon substituting this
and ${^i\!X_1}$ into the expression 
$-\res_i(\tr \bigl(g(u_0)^{-1}\dot g_1\dot A_2\bigr))$ for the
integral (\ref{eqn: CIT use}).

Now for part $2)$ of the claim.
First observe that $d_{A(u_0)}\wt\phi_1=0$ as a matrix of two-forms on
$\IP^1$. This is equivalent to
$d_{\al(u_0)}\phi_1=0$
(since $\wt \phi_1=g(u_0)^{-1}\phi_1g(u_0)$ and $\al(u)=g(u)[A(u)]$),
which follows immediately by differentiating the equation
$d(\al(u))=\al(u)\wedge \al(u)$ for the
flatness of $\al$ along $W_1$.

Therefore, by Leibniz 
$\tr(\wt\phi_1\wedge d_{A(u_0)}(g(u_0)^{-1}\dot g_2))
=-d\tr(\wt\phi_1 g(u_0)^{-1}\dot g_2)$.
Now, the Laurent expansion of $\phi_1$ at $a_i$ is
just ${^i\!\dot R_1}$, so that 
$\phi_1 ={^i\psi_1}+ {^i\!\dot R_1}$ on
$\overline D_i$ for some matrix of non-singular one-forms
${^i\psi_1}$.
Thus the integrand in $2)$ is 
$$
-d\tr (g(u_0)^{-1}\cdot{^i\psi_1}\cdot\dot g_2)
 -d\tr (g(u_0)^{-1}\cdot{^i\!\dot R_1}\cdot\dot g_2).
$$
The first term integrates to zero over the disk 
by Stokes' theorem 
(the boundary term is zero as $\dot g_2$ vanishes on $\partial\overline D_i$). 
Now from the definition  
${^i\!\dot R_1}={^i\!\dot\Lambda_1}\frac{dz}{z-a_i}$
we find that the second term is 
$\frac{\partial f}{\partial \bar z}
\frac{dz\wedge d\bar z}{z-a_i}$
where $f:=\tr (g(u_0)^{-1}{^i\!\dot\Lambda_1}\dot g_2)$.
This smooth function $f$ vanishes on $\partial\overline D_i$
and so part $2)$ of the claim follows 
using Cauchy's integral theorem since 
$f(a_i) = \langle {^i\!R_1},{^i\!\wt X_2} \rangle$.
This completes the proof of Theorem \ref{thm mmap is sp}.
\epfms

%% file: smid-imds.tex
\section{Isomonodromic Deformations and Symplectic Fibrations}	
\label{sn: imds}

Now we will consider smoothly varying the data $\ba$ that was
previously held fixed (consisting of the pole positions and the choices 
of generic connection germs $d-\iAo$ at the poles)---all 
that is now fixed throughout is the
rank $n$ of the bundles, the number $m$ of distinct poles and the
multiplicities $k_1,\ldots,k_m$ of the poles.
This leads naturally to the notion of `isomonodromic deformations' of
meromorphic connections.
Our aim is to explain, and then prove, the following:
\begin{thm}	\label{thm2}
The Jimbo-Miwa-Ueno isomonodromic deformation equations are 
equivalent to a flat symplectic Ehresmann connection on a
symplectic fibre bundle, having the moduli spaces $\wt \M^*(\ba)$ as
fibre.
\end{thm}

\subsection*{The Betti Approach to Isomonodromy}

A choice of data $\ba$ determines all the spaces 
$\wt\M^*(\ba), \wt M(\ba), \M^*(\ba)$ and 
$M(\ba)$.
Note however that the extended spaces $\wt\M^*(\ba)$ and $\wt M(\ba)$ 
only depend on the principal part of each diagonal matrix $d({^iQ})$ of
meromorphic one-forms, where $\iAo= d({^iQ}) + \iLo dz_i/z_i$.
(cf. Remark \ref{rmk: pp dept}.)
Thus if $a\in\IP^1$ it is useful to define the set $X_k(a)$ of
`order $k$ irregular types at $a$', to be the set of such principal
parts. 
Upon choosing a local coordinate $z$ vanishing at $a$ we have an
isomorphism 
\beq	\label{eqn: irtypes}
X_k(a)\cong (\IC^n\setminus\text{diagonals})\times(\IC^n)^{k-2}
\eeq
obtained by taking the coefficients of $dz/z^j$ of the 
Laurent expansion in $z$ of $A^0$, for $j=k,k-1,\ldots,2$. 
(If $k=1$ define $X_k(a):=(\text{point})$.) 

For the rest of this section we will change notation slightly,
and let $\ba$ denote data 
$(a_1,a_2,\ldots, a_m, {^1\!A^0},\ldots, {^m\!A^0})$
where ${^i\!A^0}\in X_{k_i}(a_i)$ and the $a_i$ are pairwise distinct
points of $\IP^1$.
Thus such $\ba$ determines the extended spaces 
$\wt\M^*(\ba)$ and $\wt M(\ba)$ (although we need to further specify 
exponents of formal monodromy $\iLo$ to define 
$\M^*(\ba)$ and $M(\ba)$). 

There are three manifolds of deformation parameters we will consider:
\begin{defn}

$\bullet$
The {\em basic manifold of deformation parameters} $X$ is simply the set of such $\ba$.

$\bullet$
The {\em extended manifold of deformation parameters} $\wt X$ 
is the set of such
$\ba$ together with the choice of a $k_i$-jet of a coordinate $z_i$ at
each $a_i$.

$\bullet$
If $z$ is a fixed coordinate identifying $\IP^1$ with
$\IC\cup\infty$,
the {\em Jimbo-Miwa-Ueno manifold of deformation parameters} 
$X_{\text{\rm\footnotesize JMU}}$
is the set of all such $\ba$ having $a_1=\infty$.
\end{defn}
It is easy to see these are complex manifolds, with 
$\dim(X)=\dim(X_{\text{\rm\footnotesize JMU}})+1= \dim(\wt X) - \sum k_i= 
m-mn+n\sum k_i$. 
There is an obvious embedding 
$X_{\text{\rm\footnotesize JMU}}\hookrightarrow X$ and a projection
$\wt X\twoheadrightarrow X$ forgetting the jets of local coordinates.
Moreover using the chosen coordinate $z$ there is an embedding
$X_{\text{\rm\footnotesize JMU}}\hookrightarrow \wt X$ obtained by
using the jets of the coordinates $z_i := z-a_i$ for $i=2,\ldots,m$ and
$z_1:= 1/z$.
$X_{\text{\rm\footnotesize JMU}}$ can be described very explicitly:
via (\ref{eqn: irtypes})
these coordinates identify it with
$$(\IC^{m-1}\setminus\text{diagonals}) \times
(\IC^n\setminus\text{diagonals})^{m-l}
\times(\IC^n)^{l+\sum (k_i-2)}$$
where $l=\#\{ i \ \bigl\vert \ k_i=1 \}$ is the number of simple
poles. However our aim here is more to understand the intrinsic geometry
of isomonodromic deformations, than seek explicitness, and so we will
mainly use $X$ and $\wt X$.

Now we move on to the construction of bundles over these parameter
spaces.
\begin{defn}	\label{dfn: eppb}
The bundle of extended moduli spaces $\wt\M^*$ 
is the set of isomorphism classes of data
$(V,\nabla,{\bf g},{\bf a})$ 
consisting of a {\em generic} meromorphic connection
$\nabla$ (with compatible framings ${\bf g}$)
on a {\em trivial} rank $n$ holomorphic vector bundle $V$ over a fixed
copy of $\IP^1$
such that
$\nabla$ has $m$ poles which are labelled 
$a_1,\ldots,a_m$ and the order of the pole at $a_i$ is $k_i$.
\end{defn}
It is clear from the discussion in Section \ref{sn: triv} that a
generic compatibly framed connection determines an irregular type at
each pole and it follows that there is a natural projection
$\wt\M^*\twoheadrightarrow X$ onto 
the manifold $X$ of deformation parameters, taking
the pole positions and the irregular types.
The fibre of this projection over a point $\ba\in X$ is
the extended moduli space $\wt \M^*(\ba)$. 
The results of Section \ref{sn: triv} now yield the following, which
will amount to half of Theorem \ref{thm2}:
\begin{prop}	\label{prop: sp fibrn}
The bundle $\wt \M^*$ of extended moduli spaces 
is a complex manifold and the projection above expresses it as a 
locally trivial {\em symplectic}
fibre bundle over $X$.
In particular $\wt\M^*$ has an intrinsic complex Poisson structure,
its foliation by symplectic leaves is fibrating and the space of
leaves is $X$.
\end{prop}
\pf
The only non-trivial part left is to see that $\wt \M^*$ is locally trivial
as a bundle of symplectic manifolds.
The decoupling lemma from Section \ref{sn: triv} is useful here.
Choose $m$ disjoint open disks $D_i\subset \IP^1$ and choose a
coordinate $z$ on $\IP^1$ which is non-singular on all the $D_i$'s.
Restrict to the open subset $X'$ of $X$ having $a_i\in D_i$ for each
$i$. Let $z_i:=z-a_i$.
Now, from Proposition \ref{prop: main pp}, over $X'$
any fibre $\wt\M^*(\ba)$ can be identified (using the coordinates
$z_i)$ 
with a symplectic submanifold of $\wt O_1\times\cdots\times\wt O_m$
(e.g. as the subset of $\mu_G^{-1}(0)$ which has $^1\!g_0=1$).
In turn, using Lemma \ref{lem: decoupling}, $\wt\M^*(\ba)$
is identified (if all $k_i\ge 2$)
with the symplectic manifold
$$(T^*G)^{m-1}\times{^1O_{\!B}}\times\cdots\times {^mO_{\!B}}$$
where ${^iO_{\!B}}$ is the $B_{k_i}$-coadjoint orbit through
the element of $\lb_{k_i}^*$ determined by $\iAo$ (on expanding $\iAo$
with respect to $z_i$ and replacing $z_i$ by $\ze$).
Thus the dependence of $\wt\M^*(\ba)$ on $\ba$ is clear: as $\ba$
varies, the orbit ${^iO_{\!B}}$ moves around in $\lb_{k_i}^*$.
The key fact now is that $B_{k_i}$ is a nilpotent Lie group:
coadjoint orbits of nilpotent Lie groups are diffeomorphic to vector
spaces and admit {\em global} Darboux coordinates.
Indeed M.Vergne \cite{Vergne72} shows how to find $\dim({^iO_{\!B}})$ functions
on $\lb_{k_i}^*$ which restrict to global Darboux coordinates on any 
$^iO_{\!B}$ that arises as $\iAo$ varies ($^iO_{\!B}$ is always a 
generic orbit).
Such coordinates immediately give a symplectic trivialisation of 
$\wt \M^*$ over $X'$.
(If $k_i=1$ for some $i$ then $\wt O_i$ is a fixed symplectic
submanifold of $T^*G$; there is no ${^iO_{\!B}}$ factor to worry
about.)
\epf

Similarly there is a fibre bundle $\wt M$ over $X$ whose fibres are
the extended monodromy manifolds $\wt M(\ba)$.
The key feature of the bundle $\wt M$ is that it has a canonical
complete flat Ehresmann
connection on it---in other words there is a canonical isomorphism
between nearby fibres.
In essence this connection arises by `keeping the monodromy data
constant' so we will call it the {\em isomonodromy connection}.
There is a subtlety however because it is the Stokes matrices which are
held 
constant locally, rather than the Stokes factors:
For example any anti-Stokes direction with multiplicity greater than
one can break up into distinct anti-Stokes directions under
arbitrarily small deformations of the data $\ba$, and the dimensions
of the groups of Stokes factors jump accordingly (so the 
notion of keeping the Stokes factors constant 
makes no sense directly, in general).
A precise description of the isomonodromy connection is as follows.

Suppose $\ba\in X$ is a choice of pole positions $(a_1,\ldots,a_m)$
and irregular types.
Choose disjoint open discs $D_i\subset \IP^1$ with $a_i\in\IP^1$,
together with a coordinate on each disc (so directions at $a_i$ can be
drawn as lines on $D_i$). 
If we choose a set of tentacles $\cT$ (see Definition 
\ref{def: tentacles}) then there is, from Proposition 
\ref{prop: alg im}, an isomorphism 
$\wt\varphi_\T: \wt M(\ba)\to \wcc_1\times\cdots\times\wcc_m\spq G$
to the explicit monodromy manifold (which is completely independent of
$\ba$). 
The point is that, by continuity, 
there is a small open neighbourhood $U_\ba$ of $\ba$ in $X$ such
that if $\bap$ moves around in $U_\ba$, then none of the anti-Stokes
directions at $a'_i$ cross over the base-point $p_i\in D_i$ chosen as
part of the tentacles.
Thus using the maps $\wt\varphi_\cT$ (with $\cT$ fixed and $\bap$ varying) 
we get a local trivialisation of $\wt M$ over $U_\ba$.
Repeating this process gives an open cover of $X$ with a choice of
trivialisation of $\wt M$ over each patch.
This describes the bundle $\wt M$ explicitly with clutching
functions of the form
$\wt\varphi_{\cT_1}\circ\wt\varphi_{\cT_2}^{-1}$.
Now the fact that these clutching functions are constant with respect
to the parameters $\ba\in X$ means that we have a well
defined flat connection on $\wt M$ (the local horizontal sections of which  
have constant 
explicit monodromy data 
$({\bf C},\bs,{\bf \Lambda'})\in\wcc_1\times\cdots\times\wcc_m\spq G$). 
This is the isomonodromy connection.

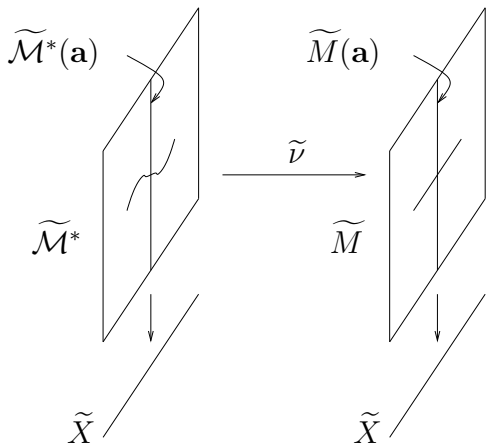
\begin{figure}[h] 	
\centerline{\input{extimdstex.pstex_t}}	
\caption{Isomonodromic Deformations} \label{imd fig}
\end{figure}

Now we want to define the relative version of the extended monodromy map.
However recall from Proposition \ref{prop: coord depce}, that this
requires a choices of coordinate jets.
Thus we first pull both bundles $\wt\M^*$ and $\wt M$ back to the
extended manifold $\wt X$ of deformation parameters along the
projection $\wt X\twoheadrightarrow X$. (These bundles over $\wt X$
will also be denoted $\wt\M^*$ and $\wt M$ but this should not lead to
confusion.) 
Then, using the jets of coordinates encoded in $\wt X$, 
the fibrewise monodromy maps
fit together to define a holomorphic bundle map,
$\wt\nu:\wt\M^*\to \wt M$
between the bundles over $\wt X$.
(As before this is holomorphic since the
canonical solutions depend holomorphically on parameters.)
\begin{defn}
The {\em isomonodromy connection on $\wt\M^*$} is the pull-back of the
isomonodromy connection on $\wt M$ along $\wt \nu$.
\end{defn}
See Figure \ref{imd fig}.
The point is that $\wt\nu$ is a highly nonlinear map with respect to the
explicit descriptions of the bundles $\wt\M^*$ and $\wt M$;
whilst being trivial on $\wt M$, the isomonodromy connection defines
interesting nonlinear differential equations on $\wt \M^*$, 
such as the Painlev\'e or Schlesinger equations
(indicated by a wavy line in the figure).

Equivalently one may view $\wt \nu$ as a kind of nonlinear Fourier-Laplace
transform (the `monodromy transform'), converting hard nonlinear
equations on the left-hand side into trivial equations on the right.
The image of $\wt \nu$ is a subset of the degree zero component 
$\wt M_0$ and as before, for dimensional reasons (since it is injective
and holomorphic)
$\wt \nu$ is  biholomorphic onto its image.
Moreover Miwa \cite{Miw} has proved that the inverse 
$\wt \nu^{-1}:\wt M_0\to \wt \M^*$ is {\em meromorphic}, so that local
horizontal sections of the isomonodromy connection on $\wt\M^*$ will
develop at worst poles when extended around $\wt X$: this is the
Painlev\'e property of the equations.
In particular this implies the image of $\wt \nu$ is the complement of
a divisor in $\wt M_0$.

Note that the isomonodromy connections are equivariant under the 
$PSL_2(\IC)$ action on the bundles $\wt \M^*, \wt M$, induced from
automorphisms of $\IP^1$.

To be precise, the isomonodromy equations of Jimbo, Miwa and Ueno are
the equations for horizontal sections of the restriction of the
isomonodromy connection on $\wt \M^*$ to $X_\jmu\hookrightarrow \wt X$,
as we will explain in the appendix.
The key idea required to actually write down equations for such
horizontal sections is the following recharacterisation of the
isomonodromy connection, which will also be 
very useful in the proof of Theorem \ref{thm2}.

\begin{rmk}	\label{new features}
Observe that in the order two pole case $k_i=2$, the $B_{k_i}$ 
coadjoint orbit $^iO_{\!B}$ in  Proposition \ref{prop: sp fibrn}
is just a point. Thus if there are no poles of order three or more,
the bundle $\wt\M^*$ has a {\em canonical} symplectic trivialisation
(which is not directly related to the isomonodromy connection) and so
the isomonodromy equations will be naturally identified with 
time-dependent flows on a fixed symplectic manifold. 
In general however, 
choices are needed in the use of Vergne's theorem in  Proposition
\ref{prop: sp fibrn}, so we 
do not know a natural way to make such an identification.
In particular, 
one must find/choose such a symplectic trivialisation
before the notion of time-dependent Hamiltonians for isomonodromy 
even makes sense.
This is a question we hope to return to in the future.  
(One suspects such a trivialisation arises naturally by requiring
Hamiltonians to come from the logarithmic derivative of the
Jimbo-Miwa-Ueno $\tau$ function.)
\end{rmk}

\subsection*{De$\!$ Rham Approach to Isomonodromy}

Suppose $\pi:Y\to X$ is some fibration over $X$, with manifolds $Y_t$
as fibres.
Replacing each $Y_t$ by its cohomology $\hh^\bullet(Y_t,\IC)$ yields a
vector bundle $\hh_\rel^\bullet(Y,\IC)\to X$.
This vector bundle has a natural flat connection on it: the
Gauss-Manin connection. 
One way to see this is from the homotopy invariance of cohomology: if
$\Delta\subset X$ is an open ball then $Y\vert_\Delta$ is homotopy
equivalent to any fibre $Y_t\subset Y\vert_\Delta$ 
so there is a
canonical isomorphism $\hh^\bullet(Y_t,\IC)\cong \hh^\bullet(Y_s,\IC)$
for any $s,t \in\Delta$.
Alternatively there is a de$\!$ Rham approach as follows. Given a closed
differential form $\theta_t$ on a fibre $Y_t$, choose any closed form
$\theta$ on $Y\vert_\Delta$ extending $\theta_t$,
and let $\theta_s$ be the restriction of
$\theta$ to $Y_s$. The cohomology class of $\theta$ in
$\hh^\bullet(Y\vert_\Delta,\IC)$ is
uniquely determined by the cohomology class of $\theta_t$ or of
$\theta_s$: this process defines the isomorphism 
$\hh^\bullet(Y_t,\IC)\cong \hh^\bullet(Y_s,\IC)$ over $\Delta$.

At least for $\hh^1$, this generalises to non-Abelian cohomology,
replacing $\IC$ by $G=GL_n(\IC)$.
Topologically $\hh^1(Y_t,G)= \Hom(\pi_1(Y_t),G)/G$ is the set of
conjugacy classes of fundamental group representations.
These fit together into a (non-linear) fibre bundle 
$\hh_\rel^1(Y,G)\to X$, which again clearly has a natural 
flat (Ehresmann) connection on it,
due to the homotopy invariance of the fundamental group: the
Gauss-Manin connection in non-Abelian cohomology. 
Simpson \cite{Sim94ab} refers to this as the Betti approach and
studies the corresponding de$\!$ Rham version. 
In the non-Abelian case, one-forms are replaced by connections 
on vector bundles, closedness is replaced by flatness, and
the notion of differing by an exact form is replaced by gauge equivalence.
Thus the de$\!$ Rham version of $\hh^1(Y_t,G)$ is the set of isomorphism
classes of flat connections on rank $n$ vector bundles over $Y_t$, and
the isomorphism $\hh^1(Y_t,G)\cong \hh^1(Y_s,G)$ arises by extending a flat
connection over a fibre $Y_t$ to a flat connection over the family 
$Y\vert_\Delta$ and then
restricting to $Y_s$.

The main realisation now is that one can very usefully 
view the isomonodromy
connection described above as the analogue in the meromorphic case 
of this non-Abelian
Gauss-Manin connection.
This emphasises the basic geometrical nature of isomonodromy and
suggests many generalisations (we are, after all, working over $\IP^1$
with $G=GL_n(\IC)$).
Note however, the necessity of having explicit descriptions of the 
moduli spaces in order to have explicit equations: the distinction
between $\wt\M^*(\ba)$ and $\wt\M(\ba)$ is important.

Thus, in the de$\!$ Rham approach,
horizontal sections of the isomonodromy connection on $\wt\M^*$ over 
some ball $\Delta\subset \wt X$ are   
related to flat meromorphic connections on vector bundles over
$\IP^1\times \Delta$.
This alternative approach was one of the main results of
\cite{JMU81} (although not expressed in
these terms).
More precisely, in the extended case, the following holds:
\begin{thm}[see \cite{JMU81}]	\label{JMUs main thm}
Let $\Delta\subset \wt X$ be an open ball.
Then there is a canonical one to one correspondence between horizontal
sections of the isomonodromy connection on $\wt\M^*$ over $\Delta$ and
isomorphism classes of triples $(V,\nabla,{\bf g})$ consisting of 
flat meromorphic connections $\nabla$ on vector bundles $V$ over
$\IP^1\times \Delta$ with {\em good} compatible framings $\bf g$, 
such that for any $t\in \Delta$ the restriction of 
$(V,\nabla,{\bf g})$ to the projective line $\IP^1\times\{ t\}$
represents an element in the fibre $\wt \M^*_t$.
\end{thm}
\noindent {\bf Sketch.}\  
See the appendix for more details, and in particular 
for the definition of `good' compatible framings.
To go from such triples $(V,\nabla,{\bf g})$ to sections of 
$\wt\M^*$ over $\Delta$, simply restrict to the $\IP^1$ fibres. 
Lemma \ref{lem: can hoz solns} shows
why the flatness of $\nabla$ implies the isomonodromicity of this
family of connections over $\IP^1$. 
Conversely, suppose we have a horizontal
section of the isomonodromy connection on $\wt\M^*$ over $\Delta$, or
equivalently a compatibly framed 
isomonodromic family $\nabla_t$ of meromorphic connections over $\IP^1$,
parameterised by $t\in\Delta$.
Then for each fixed $t$ we have a canonical basis of horizontal solutions of
$\nabla_t$ on each sector at each pole on $\IP^1\times\{ t\}$.
The key idea is that as  $t$ varies, these bases (where defined) vary
holomorphically with $t$ and $\nabla$ is defined by declaring all of
these bases to be horizontal sections of it. 
The isomonodromicity of the original family implies this $\nabla$ is
well-defined and flat.
Moreover one can deduce that $\nabla$ is meromorphic and,
by summing its principal parts, write down an algebraic expression 
for $\nabla$ in terms of the original horizontal section.
This leads directly to the explicit deformation equations.
\hfill$\square$\\

The non-extended version can easily be deduced from the above result, by
forgetting the framings, and is closer in spirit to the non-singular
(Gauss-Manin) case.
First choose an $m$-tuple $\bL$ 
of diagonal $n\times n$ matrices. 
Then define bundles $\M^*(\bL)$ and $M(\bL)$ over the space $X$ of
deformation parameters, by restricting the bundles
$\wt\M^*\to X$ and $\wt M\to X$ to the
subsets which have 
exponents of formal monodromy $\bL$ and quotienting by the action of
$T^m\cong(\IC^*)^{nm}$. 
The isomonodromy connection
on $\wt M$ 
descends to
induce a canonical isomorphism between nearby fibres of $M(\bL)\to X$,
and in turn we obtain a well-defined notion of local horizontal sections of the
isomonodromy connection on $\M^*(\bL)\to X$. 
Immediately we obtain the following (see Malgrange \cite{Mal-imd1,Mal-imd2}
for some global statements
along these lines): 
\begin{cor}
Horizontal
sections of the isomonodromy connection on $\M^*(\bL)$ 
over $\Delta\subset X$ correspond canonically
to isomorphism classes of pairs $(V,\nabla)$ consisting of 
flat meromorphic connections $\nabla$ on vector bundles $V$ over
$\IP^1\times \Delta$, 
such that for any $t\in \Delta$ the restriction of 
$(V,\nabla)$ to the projective line $\IP^1\times\{ t\}$
represents an element in the fibre $\M^*(\bL)_t$.
\end{cor}
\subsection*{Isomonodromic Deformations are Symplectic}
\label{ssn: imds are sp}
Now we will establish the second part of Theorem \ref{thm2},
thereby revealing the symplectic nature of the full family of
Jimbo-Miwa-Ueno isomonodromic deformation equations:
\begin{thm}	\label{thm: imds are sp}
The isomonodromy connection on the bundle $\wt\M^*\to \wt X$ of extended
moduli spaces, is a symplectic connection.
In other words, the local analytic diffeomorphisms induced by the isomonodromy
connection
between the fibres of $\wt \M^*$ are symplectic diffeomorphisms.
\end{thm}

\pf
We will show that arbitrary, small, isomonodromic deformations
induce symplectomorphisms.
Let $u_0$ be any point of $\wt\M^*$ and
let $x_0$ be the image of $u_0$ in $\wt X$.
Let $\gamma$ be any holomorphic map from the open unit disk
$\ID\subset\IC$
into $\wt X$ such that $\gamma(0)=x_0$.
For $t\in \ID$, let
$\wt\M^*_t$ denote the (symplectic) 
extended moduli space which is the fibre of
$\wt\M^*$ over $\gamma(t)$.
The standard vector field $\partial/\partial t$ on $\ID$
gives a vector field on $\gamma(\ID)\subset \wt X$ which we lift 
to a vector field $V$ on 
$\wt\M^*\vert_{\gamma(\ID)}$, transverse to the
fibres $\wt\M^*_t$,  using the isomonodromy connection. 
This lifted vector field may be integrated throughout a neighbourhood of 
$u_0$ in $\wt \M^*\vert_{\gamma(\ID)}$.
Concretely, this means that there is a contractible 
neighbourhood $U$ of $u_0$ in $\wt\M^*_0$, 
a neighbourhood $\Delta\subset\ID$ of $0$ in $\IC$
and a holomorphic map
$F:U\times\Delta\to\wt\M^*\vert_{\gamma(\Delta)}$
such that for all $u\in U$ and $t\in\Delta$:
$$F(u,t)\in\wt\M^*_t ,\quad F(u,0)=u\in\wt\M^*_0 \quad \text{and}\quad
\frac{\partial F}{\partial t}(u,t)=V_{F(u,t)}.$$
In particular for each $t\in\Delta$ we have a symplectic form 
$\omega_t:= (F\vert_t)^*(\omega_{\wt\M^*_t})$ on $U$,
where $\omega_{\wt\M^*_t}$ is the symplectic form defined on the
extended moduli space $\wt\M^*_t$ in Section \ref{sn: triv} and 
$F\vert_t=F(\cdot,t):U\to\wt\M^*_t$.
Now, given any two tangent vectors $W_1,W_2$ to $U$ at $u_0$,
it is sufficient for us to show that the function 
$\omega_t(W_1,W_2)$ of $t$
is {\em constant} in some neighbourhood of $0\in\Delta$.

First,
as in Proposition \ref{prop: straightening} it is easy to construct 
a local universal family over the image of $F$ in $\wt\M^*$.
Pulling back along $F$ yields a family of meromorphic connections on
the trivial bundle over $\IP^1$ parameterised by $U\times\Delta$.
For each fixed $u\in U$ we get an isomonodromic family parameterised 
by $\Delta$, that is, a `vertical' meromorphic connection
on the trivial bundle over
$\Delta\times\IP^1$ (where $\IP^1$ is the vertical direction), such that
each connection on $\IP^1$ has the `same' monodromy data.
The result of Jimbo, Miwa and Ueno (Theorem \ref{JMUs main thm} above) 
then tells us how to extend this vertical connection 
to a full {\em flat} connection over $\Delta\times\IP^1$. 
From the algebraic formula (\ref{omega defn 2}) 
for this extension it is clear that
this process behaves well as we vary $u\in U$:
for each $u\in U$ we obtain a flat meromorphic connection, which we will
denote $\nabla_u$,
on the trivial bundle over $\Delta\times\IP^1$,
that depends holomorphically on $u$.
The poles of $\nabla_u$ will be denoted by
$a_1(t),\ldots,a_m(t)$ and the polar divisor in $\Delta\times\IP^1$
of $\nabla_u$ by
$\wt D=\sum\mki\Delta_i$ (these are all independent of $u\in U$).
Shrinking  $\Delta$ if necessary, choose disjoint open
discs $D_i$ in $\IP^1$ such that $a_i(t)\in D_i$ for all $t\in\Delta$.
For each $i$ let $z_i:D_i\times\Delta\to \IC$ 
be a function which, for each fixed
$t\in\Delta$ is a coordinate on $D_i$, vanishing at
$a_i(t)$ and having the $k_i$-jet at $a_i(t)$ as specified by the point of the
base $\wt X$ below $\ga(t)$.

The next step is to push everything over to the $C^\infty$ picture
where the symplectic forms are expressed simply as integrals.
To do this we choose a smooth bundle automorphism:
$g\in GL_n\bigl(C^\infty(U\times\Delta\times\IP^1)\bigr)$
which `straightens' the whole family of connections $\nabla_u$
at the same time, as in Section \ref{sn: mm is sp}.
The map $F$ into $\wt\M^*$ specifies a family of good compatible framings
$\igo:U\times\Delta_i\to GL_n(\IC)$
of $\nabla_u$ along $\Delta_i$ for each $i$ and all $u\in U$.
Use the coordinate $z_i$ to define uniquely a  family 
$\iAo := d_{\IP^1}({^iQ})$ 
of diagonal matrices of meromorphic one-forms on $D_i$, 
parameterised by $U\times \Delta$.
(Recall only the principal part of ${^iQ}$ is specified by $\wt X$:
declare the other terms are zero in its Laurent expansion with respect
to $z_i$.)
As in Proposition \ref{prop: straightening} 
the framings extend uniquely to formal isomorphisms 
$\iwg\in GL_n\bigl(\IC\flb z_i\frb\otimes\cO(U\times\Delta_i)\bigr)$
to (uniquely determined) diagonal connections
$d_{\IP^1} - d_{\IP^1}({^iQ}) - \iL(u) d_{\IP^1}z_i/z_i$.  
By definition, that the framings are good, means
$\iwg$ satisfies a stronger condition: it
transforms the Laurent expansion of $\nabla_u$ along $\Delta_i$
into a standard full connection associated to the 
normal forms for each $u$:
\begin{equation}	\label{eqn: std full conn}
\iwg[L_i(\nabla_u)]=d-d({^i\!Q(t)})-{^i\!\Lambda}(u)\frac{d(z_i)}{z_i},
\end{equation}
where $d$ denotes the exterior derivative on the product
$\Delta\times\IP^1$, rather than just $\IP^1$. 
The automorphism $g$ is now constructed using Borel's theorem, as in
Proposition \ref{prop: straightening} 
to have Taylor expansion at $a_i(t)$ equal to $\iwg$ for all
$t\in\Delta$
and for all $u\in U$.

Thus we can use $g$ to straighten the whole family $\nabla_u$ at the
same time.
Define two families of $C^\infty$ singular connections. First a family
$\wt \nabla_u:=g[\nabla_u]$
on $\Delta\times \IP^1$ parameterised by $U$, and second
$d_\al=d_{\IP^1}-\al:=\wt\nabla_u\vert_{\IP^1}$
on $\IP^1$ parameterised by $U\times \Delta$.
By construction 
the $C^\infty$ Laurent expansion of $\wt \nabla_u$ at $a_i$ 
is given by (\ref{eqn: std full conn}).
It follows, for all $u\in U$ and $t\in\Delta$, that 
$d_\al$ is an element of the extended space
$\wt\A_{\fl}(\ba_t)\subset\wt\A(\ba_t)$ of 
flat singular connections associated to $\ba_t := \{\iAo\}$.

Now differentiate $\wt \nabla_u$ and $d_\al$ with respect to $u$
along both $W_1$ and $W_2$ at $u=u_0$.
Define these derivatives to be $\Psi_j:=W_j(\wt\nabla_u)$ and 
$\psi_j:=W_j(d_\al)=\Psi_j\vert_{\IP^1}$ respectively,
for $j=1,2$.
Each $\Psi_j$ is a matrix of singular one-forms on $\Delta\times\IP^1$
and each $\psi_j$ is a matrix of singular one-forms on $\IP^1$
parameterised by $\Delta$.
Clearly $\tr(\psi_1\wedge\psi_2)=\tr(\Psi_1\wedge\Psi_2)\vert_{\IP^1}.$
Also since the Laurent expansion of $\wt \nabla_u$ is given by
(\ref{eqn: std full conn}) 
at each $a_i$ we can deduce what the Laurent expansions
of $\Psi_1$ and $\Psi_2$ are:
$L_i(\Psi_j)= W_j({^i\!\Lambda(u)})
{d_{\Delta\times\IP^1}(z_i)}/{z_i}$
for $j=1,2$ and $i=1,\ldots,m$.
It follows that 
$\tr(\Psi_1\wedge\Psi_2)$ is a {\em nonsingular} two-form on 
$\Delta\times\IP^1$ since
$L_i(\Psi_1\wedge\Psi_2)=L_i(\Psi_1)\wedge L_i(\Psi_2)=0$
for each $i$.

Now observe that for each $u\in U$ the flatness of $\nabla_u$ 
implies the flatness of $\wt \nabla_u$.
By differentiating the equation $\wt\nabla_u\circ\wt\nabla_u=0$ 
with respect to $u$ along $W_1$ and
$W_2$ we find
$\wt\nabla_{u_0}\Psi_1=0$
and
$\wt\nabla_{u_0}\Psi_2=0.$
In particular, by Leibniz, 
the two-form $\tr(\Psi_1\wedge\Psi_2)$ on
$\Delta\times\IP^1$ is {\em closed}.

Thus if we do the fibre integral over $\IP^1$ we obtain a zero-form
on $\Delta$ (i.e. a function of $t$):
$$\ip \tr(\Psi_1\wedge\Psi_2)=\ip\tr(\psi_1\wedge\psi_2).$$
This is a {\em closed} $0$-form (i.e. a {\em constant} function)
since integration over the fibre commutes with exterior
differentiation.
See for example Bott and Tu \cite{Bott+Tu} Proposition 6.14.1 
(it is important here that $\tr(\Psi_1\wedge\Psi_2)$ is nonsingular). 

Finally we appeal to Theorem \ref{thm mmap is sp} to see that for all
$t\in\Delta$:
$$\frac{1}{2\pi i}\ip\tr(\psi_1\wedge\psi_2)=
\omega_t(W_1,W_2)$$
and so the symplectic form is indeed independent of $t$.
\epf

\subsection*{Closing Remarks}
One upshot of Theorem \ref{thm2} is that the symplectic structure on each
monodromy manifold is independent of the choice of deformation
parameters; the isomonodromy connection on $\wt M$ is symplectic. This is the
generalisation of the `symplectic nature of the fundamental group'.
As in the non-singular case, one then wonders if there is an intrinsic
finite-dimensional/algebraic approach to this symplectic structure
(generalising the cup product in group cohomology).
This should be possible by combining the $C^\infty$ approach here
with the ideas of Alekseev, Malkin and Meinrenken \cite{AMM}.

Alternatively (or perhaps equivalently) a direct connection between
Stokes matrices and Poisson Lie groups was observed in 
\cite{Boa}, which we will briefly sketch here since it is 
quite intriguing.
Consider the case of connections on $\IP^1$ with just two poles, of orders one
and two respectively.
The choice of an irregular type at the order two pole determines the moduli
space $\wt\M^*(\ba)$ and the monodromy manifold $\wt M(\ba)$.
If we forget the framing at the order one pole, we obtain
$\wt\M^*(\ba)/T$ which 
is isomorphic as a Poisson manifold to (a covering of a dense
open subset of) $\g^*$.
Also $\wt M(\ba)/T$ is isomorphic to a covering of a dense
open subset of $U_+\times U_-\times\lt$.
The monodromy map extends to a map
$\nu:\g^*\to U_+\times U_-\times\lt$, taking the Stokes matrices and the
exponent of formal monodromy at $0$ of the connection
$d-(Udz/z^2+Vdz/z)$, where $V\in\g\cong\g^*$ and $U$ is a fixed diagonal
matrix with distinct eigenvalues.
The basic observation now is that $ U_+\times U_-\times\lt$ may be identified
with the simply connected Poisson Lie group $G^*$ dual to $GL_n(\IC)$.
We then claim that, 
under such identification, $\nu:\g^*\to G^*$ is a Poisson map,
where $\g^*$ and $G^*$ both have their standard Poisson structures.\footnote
{This has now been proved, cf. P.P.Boalch, math.DG/0011062} 
In particular, taking $V$ to be skew-symmetric, this claim yields a new
approach to the Poisson bracket on Dubrovin's local
moduli space of semisimple Frobenius manifolds.

\appendix
\section{}	\label{apx}
We will give more details regarding Theorem \ref{JMUs main thm},
relating flat connections to horizontal sections of the isomonodromy
connection.
This differs from \cite{JMU81} in that the coordinate dependence is isolated
here. 
At the end we will write down the deformation equations.

First some generalities on the local structure of meromorphic
connections in higher dimensions.
The local model is of a meromorphic connection $\nabla=d-\wt A$ on 
the trivial rank $n$
vector bundle over a product 
$\ID\times\Delta$
of the unit disc $\ID\subset \IC$ and some contractible
space of parameters $\Delta$.
We suppose, for each $t\in\Delta$ that the restriction  
$\nabla_t:=\nabla\vert_{\ID\times\{t\}}$ to the corresponding disc
has only one pole (of order $k$)
at some point  
$a(t)\in\ID$ and is formally equivalent to a generic diagonal connection 
$d_\ID-A^0(t)$ depending holomorphically on $t$.
Assume the divisor 
$\Delta_0:=\{(a(t),t)\}\subset\ID\times\Delta$ is smooth.
Let $z_0:\ID\times\Delta\to \IC$ be any holomorphic function vanishing on
$\Delta_0$ which restricts to a coordinate on 
$\ID\times\{t\}$ for each $t\in\Delta$ (only the $k$-jet of the Taylor
expansion of $z_0$ along $\Delta_0$ will be significant below).
Write $A^0 = d_\ID Q + \Lambda^0(t) d_\ID z_0/z_0$, as usual and
define the `standard full connection' to be
$d-\wt A^0$ where $\wt A^0:= d Q + \Lambda^0(t) d(z_0)/z_0$ and $d$
denotes the full exterior derivative on $\ID\times\Delta$.

If we  choose a compatible framing $g_0$ of $\nabla$ along $\Delta_0$
then, as in Proposition \ref{prop: straightening}, there is a unique
family of formal isomorphisms 
$\wh g\in GL_n(\IC\flb z_0\frb\otimes\cO(\Delta_0))$ satisfying 
$\wh g\vert_{\Delta_0} = g_0$ and
${\widehat g_t}[\nabla_t]= d_\ID- {A^0}$ for each fixed $t$ (after
possibly permuting the entries of $A^0$). 
The basic structural result is then:
\begin{lem}[see \cite{Mal-imd2}]	\label{lem: flat str}
If $\nabla$ is {\em flat} then $\Lambda^0$ is constant and there is
a diagonal matrix valued holomorphic function
$F\in\End_n(\cO(\Delta_0))$ (which is unique upto the addition of a
constant diagonal matrix) such that
$$\wh g[\nabla]_{\ID\times\Delta} = d - 
(\wt A^0 + \pi^*(d_{\Delta_0}F))$$
where $\pi:\ID\times\Delta\to\Delta_0$ is
the projection along the $\ID$ direction.
\end{lem}
\pf
Let $d_{\Delta_0} - B$ be the $\Delta_0$ component of the Laurent
expansion of  
$\wh g[\nabla]_{\ID\times\Delta}$ so that 
$\wh g[\nabla]_{\ID\times\Delta} = d_{\ID\times\Delta} - 
(A^0 + {{B}}).$
This is flat because $\nabla$ is. 
The ($\ID$-$\Delta_0$) part of the equation for this flatness is:
\begin{equation}	\label{vertizontal}
d_{\ID}{{B}} + d_{\Delta_0} A^0 = {A^0}\wedge {{B}} + {{B}}\wedge {A^0}.
\end{equation}
Since ${A^0}$ is diagonal this equation splits into two independent
pieces, the diagonal part and the off-diagonal part.
First we deduce that the off-diagonal part ${{B}}^\od$ of
${{B}}$ is zero:
Suppose ${{B}}^\od\ne 0$ and let $M/z_0^r$ be its leading term, 
$M\in \End_n^\od(\Omega^1_\hol(\Delta_0))$.
Equation (\ref{vertizontal}) implies 
$
d_{\ID}{{B}}^\od={A^0}\wedge {{B}}^\od + {{B}}^\od \wedge {A^0}.
$
Counting the pole orders we deduce
${{B}}^\od=0$ unless $k=1$.
If $k=1$, say ${A^0}=A_{1}^0dz_0/z_0$, then considering
coefficients of $dz_0/z_0^{r+1}$ we see
$(-r)M=[A_{1}^0,M]$
which implies $M=0$ (and therefore ${{B}}^\od=0$) since 
${A^0}$ is generic; the difference between any two eigenvalues of
$A^0_{1}$ is never the integer $-r$.
Thus ${{B}}$ is diagonal, and so (\ref{vertizontal}) now reads
$d_{\ID}{{B}}+ d_{\Delta_0} {A^0}=0$.
This implies $d_{\Delta_0} {\Lambda^0(t)}=0$ since
$d_{\ID}{{B}}$ will have no residue term, and so $\wt A^0$ is flat.
Thus
$d_{\ID}{{B}} =  -d_{\Delta_0} \wt A^0 = d_{\ID} \wt A^0$.
Hence
${{B}} = \wt A^0_{\Delta_0} + \phi(t)$
for some diagonal matrix of one-forms 
$\phi\in\End_n(\Omega_\hol^1(\Delta_0))$ where 
$\wt A^0_{\Delta_0}$ is the $\Delta_0$ component of $\wt A^0$.
Finally the ($\Delta_0$-$\Delta_0$) part of the equation for the flatness
of $d-{A^0}-{{B}}$ implies $d_{\Delta_0} {{B}}=0$. 
It follows that $d_{\Delta_0}(\phi(t))=0$ and
so, since $\Delta_0$ is contractible, $\phi=d_{\Delta_0} F$ for some 
diagonal $F\in \End_n(\cO(\Delta_0))$
\epf

This leads us to make the following:
\begin{defn}
If $\nabla$ is flat then 
a compatible framing $g_0$ of $\nabla$ 
along $\Delta_0$ is {\em good} 
if 
$\wh g[\nabla]_{\ID\times\Delta} = d - \wt A^0$
where $\wh g$ is the formal series associated to $g_0$.
\end{defn}
Thus an arbitrary compatible framing $g_0$ can be made good by replacing it
by $e^{-F}g_0$ where  $F$ is from Lemma \ref{lem: flat str}.
It is worth saying the same thing slightly differently.
In the convention we are using, the columns of the inverse $g_0^{-1}$
of the compatible framing are a basis of sections of
$V\vert_{\Delta_0}$,
where $V$ is the bundle that $\nabla$ is on. 
Thus, since good compatible framings are determined upto a constant,
there is a flat holomorphic connection $\nabla_0$ on 
$V\vert_{\Delta_0}$ whose horizontal sections are the columns of 
$g_0^{-1}$ for any good compatible framing $g_0$.
A direct definition is:
\begin{defn}	\label{def: ind conn}
If $g_0$ is any compatible framing of $\nabla$ along $\Delta_0$ then
the {\em induced connection along $\Delta_0$} is
$
\nabla_0= 
(\nabla + \wh g^{-1}\cdot \wt A^0 \cdot \wh g)\bigl\vert_{\Delta_0},
$
where $\wh g$ is the formal series associated to $g_0$.
\end{defn}
It is easy to check this definition is independent of the choice of
compatible framing and, if $g_0$ is good, then the columns of
$g_0^{-1}$ are horizontal. Moreover this definition makes sense for
non-flat $\nabla$, but then $\nabla_0$ may not be flat.
One may also check that $\nabla_0$ only
depends on $\nabla$ and the choice of $k$-jets of coordinates $z_0$.
(Also in the logarithmic case $k=1$, $\nabla_0$ coincides with the
usual (canonical) notion of induced connection
$\nabla\vert_{\Delta_0}$, 
provided 
$z_0$ satisfies $(dz_0/z_0)\vert_{\Delta_0}=0$.)
Thus one can alternatively define good framings to be the compatible
framings 
$g_0$ such that the columns of $g_0^{-1}$ are horizontal for $\nabla_0$.
The reason for restricting how the framings vary along $\Delta_0$ is
the following:

\begin{lem}[see \cite{JMU81} Theorem 3.3]	\label{lem: can hoz solns}
Let $\nabla$ be a full flat connection as above and 
let $g_0$ be a good compatible framing 
with corresponding formal series $\wh g$.
Fix any point $t_0\in \Delta$, choose a labelling of the sectors between the
anti-Stokes directions at 
$a(t_0)\in \ID\times\{ t_0\}$, and choose $\log(z_0)$ branches on
$\ID\times\{t_0\}$.
Let $\Delta'$ be a neighbourhood of $t_0\in\Delta$ such that the last
sector at $a(t_0)$ deforms into a unique sector at $a(t)$ for all
$t\in\Delta'$ (the last sector at $a(t)$).

Then the canonical fundamental solution 
$\Phi_0 := \Sigma_0(\wh g^{-1})z_0^{\Lambda^0} e^Q$ 
of $\nabla\vert_\vrt$
on the last sector at $a(t)\in \ID\times\{t\}$ varies holomorphically
with $t\in\Delta'$ and $\Phi_0(z,t)$ is a local fundamental solution of 
the original full connection $\nabla$. (Similarly on the other
sectors: just relabel.)
\end{lem}
\pf
Write $\nabla = d- \wt A$ and let $\Omega$ be the $\Delta$ component
of $\wt A$ so that $\wt A = A + \Omega$.
The aim is to show that $d_\Delta \Phi_0 = \Omega \Phi_0$.
From the definition of $\wh g$ we have 
$A+\Omega = \wh g^{-1}[\wt A^0]_{\ID\times\Delta}$
and this has $\Delta$ component
$\Omega = \wh g^{-1} \cdot\wt A^0_\Delta\cdot \wh g - \wh
g^{-1}d_\Delta\wh g.$
Now the key observation is that the equation 
$d_\Delta A=-d_{\ID}\Omega+A\wedge\Omega+\Omega\wedge A$
(from the flatness of $\nabla$)
implies that 
the matrix of one-forms $d_\Delta \Phi_0 - \Omega \Phi_0$
satisfies the equation
$d_{\ID}
(d_\Delta \Phi_0 - \Omega \Phi_0) = A (d_\Delta \Phi_0 - \Omega \Phi_0)$
(also using the fact that $d_{\ID}\Phi_0 = A \Phi_0$).
Then if we define a matrix 
$K:= \Phi_0^{-1}(d_\Delta \Phi_0 - \Omega \Phi_0)$
of one-forms it follows that  $d_{\ID}K = 0$ so that $K$ is constant in the $\ID$
direction.
Then using the fact that the asymptotic expansion of 
$\Phi_0$ in the last sector at $a(t)$ is $\wh g^{-1}z_0^{\Lambda^0} e^Q$, 
it follows that $K$
has  zero asymptotic expansion there.
(This uses the fact that the asymptotic expansions are uniform in
$t$ to see that $d_\Delta$ commutes with the operation of taking
the asymptotic expansion.)
It follows immediately that $K=0$ because $K$ is constant in the
$\ID$ direction, and so $d_\Delta \Phi_0 = \Omega \Phi_0$.
\epf

This is the main result needed to prove Theorem \ref{JMUs main thm} as
sketched. 
All that remains is to write down the deformation 
equations of Jimbo, Miwa and Ueno.
Restrict the parameter space to $X_\jmu\hookrightarrow \wt X$.
The bundle $\wt\M^*$ over $X_\jmu$ can be decribed explicitly:
using Proposition \ref{prop: main pp} (and removing the $G$ action by
fixing $^1\!g_0=1$) it is identified as a subbundle of the
trivial bundle over $X_\jmu$ with fibre 
$$(GL_n(\IC)\times\g^*_{k_1})\times\cdots\times (GL_n(\IC)\times \g^*_{k_m}).$$
When described in this way the bundle $\wt\M^*\to X_\jmu$ is identified as the
`manifold of singularity data' of \cite{JMU81}.
Now suppose we have a horizontal section of the isomonodromy
connection on $\wt\M^*$ over some ball $\Delta\hookrightarrow X_\jmu$.
From Section \ref{sn: triv}
this determines a family of meromorphic connections $d_{\IP^1}-A$ on
the trivial bundle over $\IP^1$ and compatible framings $^i\!g_0$ 
(the principal parts of $A$ lie in the $\g^*_{k_i}$'s 
using the coordinate choices).
As above we also get (algebraically) formal isomorphisms $\iwg$,
connection germs $d_{\IP^1}-\iAo$ and `full' connection germs
$d_{\IP^1\times\Delta}-\iwAo$ , where $d$ is the full exterior
derivative on ${\IP^1\times\Delta}$.
(The holomorphic terms in the expansion of $\iAo$ with respect to
$z_i$ are defined to be zero.)

From the sketch of the proof of Theorem \ref{JMUs main thm}, 
$d_{\IP^1}-A$ is the vertical component of a full connection 
$\nabla=d-\wt A$, where $\wt A = (d\Phi)\Phi^{-1}$ for any local canonical
fundamental solution 
$\Phi(z,t):= {^i\Sigma_j}({^i\wh g^{-1}})z_i^{^i\!\Lambda}e^{^i\!Q}$ 
on (say) the $j$th sector at the $i$th pole. 
These local definitions agree as the family $d_{\IP^1}-A$ is 
isomonodromic.
Let $\Omega$ denote the $\Delta$ component of $\wt A$, so 
$\wt A=A+\Omega$. 
From the definition 
we know the asymptotics of $\Phi(z,t)$ (uniformly) on the $j$th
(super)sector at the $i$th pole and so we can deduce the asymptotics
of $\Omega$:
\begin{equation}	\label{asexp omega}
\asexp_i(\Omega) = \bigl(\iwg^{-1}\cdot\iwAod\cdot \iwg\bigr) 
- \iwg^{-1} \cdot d_\Delta(\iwg) 
\end{equation}
A priori this only holds on some sector at the $i$th pole, but
choosing a different $\Phi$, we get the same expansion on every sector.
It follows that $\Omega$ is meromorphic, with {\em Laurent} expansion
(\ref{asexp omega}).
First, it follows immediately from this expression that the
compatible framing $\igo$ is a good compatible framing of $\nabla$.
Secondly it is clear that the $i$th principal part of $\Omega$ is 
the principal part of $\iwg^{-1}\cdot\iwAod\cdot \iwg$ and so is 
determined algebraically.
Also we need a formula for the induced connections $\nabla_i$ on
the polar divisors of $\nabla$.
Upon pulling $\nabla_i$ down to the base $\Delta$, from Definition 
\ref{def: ind conn}, one finds that $\nabla_i$ becomes $d_\Delta-\Theta_i$
where
\begin{equation}	\label{theta defn 2}
\Theta_i = 
{^i\!g_0^{-1}}(d_\Delta a_i) {^i\!g_1} + 
\const_{z_i}( \iwg^{-1} \cdot \iwAod\cdot \iwg) -
\const_{z_i}(\Omega)
\end{equation}
with $\iwg = {^i\!g_0} + {^i\!g_1}\cdot z_i + O(z_i^2)$ 
and where $\const_{z_i}$ takes the constant term in the
Laurent expansion with respect to $z_i$.
Since we are working in the trivialisation determined by the first
framing (${^1\!g_0}=1$), we have $\Theta_1=0$ and so the expression 
(\ref{theta defn 2}) determines the constant term in the expansion of
$\Omega$ at $a_1=\infty$. 
Thus $\Omega$ is completely determined by this constant and the
principal parts:
\begin{equation}	\label{omega defn 2}
\Omega = \const_{z_1}(\owg^{-1}\cdot\owAod\cdot \owg) +
\sum_{i=1}^m \PP_{z_i}\bigl(\iwg^{-1}\cdot\iwAod\cdot \iwg\bigr).
\end{equation}

Now the flatness of the full connection $\nabla$ over
$\Delta\times\IP^1$ implies two equations.
Firstly 
$d_\Delta\Omega=\Omega\wedge\Omega$,
which says that $\Omega$ is a family of
flat connections on $\Delta$ depending rationally on the `spectral
parameter' $z$; a situation that often arises in soliton theory.
Secondly 
\begin{equation}	\label{flat}
d_\Delta A=-d_{\IP^1}\Omega+A\wedge\Omega+\Omega\wedge A.
\end{equation}
Also the `goodness' of the compatible framings $\igo$ implies that 
\begin{equation}	\label{fram}
d_\Delta(\igo)=-(\igo)\Theta_i.
\end{equation}

Note that the formulae (\ref{theta defn 2}) and (\ref{omega defn 2})
for $\Omega$ and $\Theta_i$ 
make sense for an arbitrary section of the bundle $\wt\M^*$ so that the
equations (\ref{flat}) and (\ref{fram})
amount to a coupled system of nonlinear {\em algebraic} differential 
equations for horizontal
sections $s=({\bf g},{^1\!A},\ldots,{^m\!A})$
of the isomonodromy connection on $\wt\M^*$ over $X_\jmu$:
They are the Jimbo-Miwa-Ueno isomonodromic deformation equations
\cite{JMU81}.

A number of examples are given in \cite{JM81, JMU81} and
in particular the cases of the Schlesinger equations and the six
Painlev\'e equations are explained.

%% file: extimdstex.pstex_t
\begin{picture}(0,0)%
\includegraphics{extimdstex.pstex}%
\end{picture}%
\setlength{\unitlength}{1973sp}%
\begingroup\makeatletter\ifx\SetFigFont\undefined%
\gdef\SetFigFont#1#2#3#4#5{%
  \reset@font\fontsize{#1}{#2pt}%
  \fontfamily{#3}\fontseries{#4}\fontshape{#5}%
  \selectfont}%
\fi\endgroup%
\begin{picture}(6027,5610)(1186,-7159)
\put(1201,-2236){\makebox(0,0)[lb]{\smash{{\SetFigFont{12}{14.4}{\rmdefault}{\mddefault}{\updefault}{\color[rgb]{0,0,0}$\wt\M^*(\ba)$}%
}}}}
\put(4726,-3511){\makebox(0,0)[lb]{\smash{{\SetFigFont{12}{14.4}{\rmdefault}{\mddefault}{\updefault}{\color[rgb]{0,0,0}$\wt\nu$}%
}}}}
\put(1501,-4636){\makebox(0,0)[lb]{\smash{{\SetFigFont{12}{14.4}{\rmdefault}{\mddefault}{\updefault}{\color[rgb]{0,0,0}$\wt\M^*$}%
}}}}
\put(5551,-7036){\makebox(0,0)[lb]{\smash{{\SetFigFont{12}{14.4}{\rmdefault}{\mddefault}{\updefault}{\color[rgb]{0,0,0}$\wt X$}%
}}}}
\put(1951,-7036){\makebox(0,0)[lb]{\smash{{\SetFigFont{12}{14.4}{\rmdefault}{\mddefault}{\updefault}{\color[rgb]{0,0,0}$\wt X$}%
}}}}
\put(4951,-2236){\makebox(0,0)[lb]{\smash{{\SetFigFont{12}{14.4}{\rmdefault}{\mddefault}{\updefault}{\color[rgb]{0,0,0}  $\wt M(\ba)$}%
}}}}
\put(5251,-4636){\makebox(0,0)[lb]{\smash{{\SetFigFont{12}{14.4}{\rmdefault}{\mddefault}{\updefault}{\color[rgb]{0,0,0}  $\wt M$}%
}}}}
\end{picture}%